\documentclass[12pt,a4paper,reqno]{amsart}
\usepackage{amsmath} 
\usepackage{amssymb} 
\usepackage{amscd} 
\usepackage{amsxtra}
\usepackage{epic,eepic}

\numberwithin{equation}{section}

\theoremstyle{plain}
\newtheorem{theorem}{Theorem}[section]
\newtheorem{theorem-definition}[theorem]{Theorem-Definition} 
\newtheorem{proposition}[theorem]{Proposition}
\newtheorem{corollary}[theorem]{Corollary}
\newtheorem{lemma}[theorem]{Lemma}

\theoremstyle{definition} 
\newtheorem{definition}[theorem]{Definition}
\newtheorem{remark}[theorem]{Remark}
\newtheorem{condition}[theorem]{Condition} 

\newtheorem{example}[theorem]{Example}

\newtheorem*{notation*}{Notation}
\newtheorem*{remark*}{Remark}

\newcommand{\C}{{\mathbb C}}
\newcommand{\R}{{\mathbb R}}
\newcommand{\Z}{{\mathbb Z}}

\newcommand{\A}{{\mathbb A}}
\newcommand{\bL}{{\mathbb L}}
\newcommand{\T}{{\mathbb T}}

\newcommand{\Q}{{\mathbb Q}}

\newcommand{\Proj}{{\mathbb P}}

\newcommand{\bt}{{\boldsymbol t}}
\newcommand{\bc}{{\boldsymbol c}}
\newcommand{\bs}{{\boldsymbol s}}
\newcommand{\be}{{\boldsymbol e}}
\newcommand{\bof}{{\boldsymbol f}}
\newcommand{\bg}{{\boldsymbol g}}
\newcommand{\bpartial}{\boldsymbol{\partial}}

\newcommand{\Qp}{\boldsymbol{P}}
\newcommand{\bq}{\boldsymbol{q}}
\newcommand{\hQp}{\hat{\boldsymbol{P}}}

\newcommand{\bb}{\boldsymbol{b}}

\newcommand{\mV}{{\mathcal{V}}}
\newcommand{\mH}{{\mathcal{H}}}
\newcommand{\mL}{{\mathcal{L}}}
\newcommand{\mW}{{\mathcal{W}}}
\newcommand{\mA}{{\mathcal{A}}}

\newcommand{\mR}{{\mathcal{R}}}
\newcommand{\mD}{{\mathcal{D}}}
\newcommand{\hmD}{{\hat{\mathcal{D}}}}
\newcommand{\qmD}{\mathsf{D}}
\newcommand{\mP}{{\mathcal{P}}}
\newcommand{\mU}{{\mathcal{U}}}
\newcommand{\mN}{{\mathcal{N}}}
\newcommand{\mO}{{\mathcal{O}}}
\newcommand{\qmO}{\mathsf{O}}
\newcommand{\mI}{{\mathcal{I}}}
\newcommand{\mE}{{\mathcal{E}}}

\newcommand{\mEh}{\mathcal{E}^\hbar}

\newcommand{\hmEh}{\hat{\mathcal{E}}^\hbar}
\newcommand{\hmE}{\hat{\mathcal{E}}}
\newcommand{\qmEh}{\mathsf{E}^\hbar}
\newcommand{\hqmEh}{\hat{\mathsf{E}}^\hbar}
\newcommand{\mFh}{\mathcal{F}^\hbar}
\newcommand{\mF}{\mathcal{F}}

\newcommand{\mJ}{\mathcal{J}}
\newcommand{\hmJ}{\hat{\mathcal{J}}}

\newcommand{\frF}{\mathfrak{F}}
\newcommand{\frI}{\mathfrak{I}}

\newcommand{\frc}{\mathfrak{c}}

\newcommand{\hq}{\hat{q}}
\newcommand{\hatt}{\hat{t}}
\newcommand{\hPhi}{\hat{\Phi}}

\newcommand{\hH}{\mathbb{H}}
\newcommand{\hA}{\hat{A}}
\newcommand{\hAA}{\hat{\A}}
\newcommand{\hL}{\hat{L}} 
\newcommand{\hJ}{\hat{J}}
\newcommand{\htau}{\hat{\tau}}

\newcommand{\Lie}{\operatorname{Lie}}
\newcommand{\Hom}{\operatorname{Hom}}

\newcommand{\Res}{\operatorname{Res}}
\newcommand{\ad}{\operatorname{ad}}
\newcommand{\ch}{\operatorname{ch}}
\newcommand{\Euler}{\operatorname{Euler}}

\newcommand{\Map}{\operatorname{Map}} 

\newcommand{\id}{\operatorname{id}}
\newcommand{\End}{\operatorname{End}}
\newcommand{\Aut}{\operatorname{Aut}}
\newcommand{\Pic}{\operatorname{Pic}} 
\newcommand{\Loc}{\operatorname{Loc}}
\newcommand{\rank}{\operatorname{rank}}

\newcommand{\ord}{\operatorname{ord}}

\newcommand{\modif}{\operatorname{modif}}

\newcommand{\frm}{\frak{m}}
\newcommand{\bwp}{\boldsymbol{\wp}}
\newcommand{\hbwp}{\hat{\bwp}}
\newcommand{\hwp}{\hat{\wp}}
\newcommand{\smnabla}{\boldsymbol{\nabla}}
\newcommand{\Itw}{I^{\rm tw}} 
\newcommand{\Qla}{Q_\lambda} 

\newcommand{\barop}{\overline{\phantom{A}}}

\def\ov#1{\overline{#1}}
\def\corrV#1{\left\langle{#1}\right\rangle^{\mathcal{V}}}
\def\parti#1{{q^{#1}\frac{\partial}{\partial q^{#1}}}}
\def\parfrac#1#2{{\frac{\partial #1}{\partial #2}}}
\def\pair#1#2{\langle #1,#2\rangle}
\def\pairV#1#2{\langle #1,#2\rangle^{\mathcal{V}}}

\begin{document}
\title{Quantum D-modules and Generalized Mirror Transformations}
\author{Hiroshi Iritani}
\address{Department of Mathematics, Graduate School of Science, 
Kyoto University, Oiwake-cho, Kitashirakawa, Sakyo-ku, Kyoto, 606-8502, Japan.} 
\email{iritani@math.kyoto-u.ac.jp}
\begin{abstract}
In the previous paper \cite{iritani-EFC}, 
the author defined equivariant Floer cohomology  
for a complete intersection in a toric variety and  
showed that it is isomorphic to the small quantum $D$-module
after a mirror transformation 
when the first Chern class $c_1(M)$ of the tangent bundle is nef.  
In this paper, even when $c_1(M)$ is not nef, 
we show that the equivariant Floer cohomology reconstructs 
the big quantum $D$-module 
under certain conditions on the ambient toric variety. 
The proof is based on a mirror theorem by Coates and Givental 
\cite{coates-givental}. 
The reconstruction procedure here gives a 
generalized mirror transformation first observed 
by Jinzenji in low degrees \cite{jinzenji1,jinzenji2}. 
\end{abstract}
\maketitle

\tableofcontents

\section{Introduction}

The $S^1$-equivariant Floer cohomology proposed by Givental 
\cite{givental-ICM,givental-homologicalgeom} is a conjectural 
semi-infinite cohomology of the free loop space of a symplectic manifold $M$. 
The $S^1$ action here is the rotation of loops. 
Givental conjectured that this should 
have a natural $D$-module structure and 
is isomorphic to the small quantum $D$-module defined by 
quantum cohomology. 
Here, the small quantum $D$-module (henceforth small QDM) 
is a trivial $H^*(M)$-bundle over $H^2(M,\C^*)$ 
with a flat connection $\smnabla^\hbar$, called the Dubrovin connection: 
\[
\smnabla^\hbar_a = \hbar Q^a \parfrac{}{Q^a} + p_a * 
\]
where $Q^1,\dots, Q^r$ is a co-ordinate system on $H^2(M,\C^*)$ and 
$p_1*,\dots,p_r*$ denotes the quantum multiplication 
by a basis of $H^2(X)$ dual to $Q^a$. 

In \cite{iritani-EFC}, the author constructed the 
equivariant Floer cohomology $FH_{S^1}^*$ 
for toric complete intersections 
as an inductive limit of ordinary equivariant cohomology. 
He also introduced {\it abstract quantum $D$-module} 
(henceforth AQDM) which generalizes the small QDM. 
AQDM is a module over a Heisenberg algebra 
$\mD=\C[\hbar][\![Q_1,\dots,Q_r]\!]\langle \Qp_1,\dots,\Qp_r\rangle$ 
where 
\[
[\Qp_a, Q^b] = \hbar \delta_a^b Q^b, \quad 1\le a,b\le r.  
\]   
In case of the small QDM, $\Qp_a$ is given by 
the Dubrovin connection $\smnabla^\hbar_a$. 
The author showed that the equivariant Floer cohomology 
$FH_{S^1}^*$ has the structure of an AQDM. 
Using a mirror theorem by Givental 
\cite{givental-mirrorthm-projective, givental-mirrorthm-toric}, 
he also showed that $FH_{S^1}^*$ is isomorphic to the small QDM
as an AQDM if the first Chern class $c_1(M)$ of the tangent bundle is nef. 

This isomorphism of two AQDMs --- 
small QDM and equivariant Floer cohomology --- 
is the same as {\it mirror transformation} 
in the context of mirror symmetry. 
The small QDM is the A-model 
whereas the $FH_{S^1}^*$ plays the role of the B-model. 
More precisely, differential equations 
satisfied by a generator of $FH_{S^1}^*$  
coincide with the Picard-Fuchs equations arising from the B-model. 
By the mirror transformation, we can compute 
the quantum cohomology \emph{i.e.} the genus zero Gromov-Witten invariants. 

When the first Chern class $c_1(M)$ is not nef, however, 
$FH_{S^1}^*$ is not 
isomorphic to the small QDM, 
and the ordinary mirror transformation does not work.  
A {\it generalized} mirror transformation is a prescription 
to remedy the situation. 
Here we need to consider not only the small QDM but also the big one.  
The big quantum $D$-module (henceforth big QDM) 
is a trivial $H^*(M)$-bundle on the total cohomology group 
$H^*(M)$ with the flat Dubrovin connection  
and its restriction to $H^2(M)$ can be 
identified with the small QDM. 
We introduce the notion of {\it big AQDM} 
which generalizes the big QDM. 
The AQDM in the previous paper \cite{iritani-EFC} 
is called small AQDM in this paper.  
We show that the small AQDM reconstructs a big AQDM uniquely 
under certain conditions corresponding to that $H^*(M)$ is 
generated by $H^2(M)$. 
In particular, the small AQDM $FH_{S^1}^*$ reconstructs 
a certain big AQDM $\hmEh$ such that $FH_{S^1}^*$ is isomorphic 
to the restriction of $\hmEh$ to a subspace of its base space. 
Under Condition \ref{cond:ambienttoric} on the ambient toric variety, 
we will show that the reconstructed $\hmEh$ 
is isomorphic to the big QDM of $M$ (Theorem \ref{thm:main}). 
In this way, we can recover the big quantum cohomology from $FH_{S^1}^*$.  
Note that both $FH_{S^1}^*$ and the small QDM of $M$ 
are obtained as restrictions of $\hmEh$ to 
certain subspaces, but {\it the loci are  different} in general, 
\emph{i.e.} the locus of $FH_{S^1}^*$, 
which can be considered to be the B-model locus, 
may not coincide with the linear subspace $H^2(M)$ 
(see Figure \ref{fig:embedding}). 
The reconstruction from $FH_{S^1}^*$ to $\hmEh$ 
can be viewed as a generalization of Kontsevich and Manin's 
reconstruction theorem \cite{kontsevich-manin}. 
Hertling and Manin \cite{hertling-manin} also proved 
a similar reconstruction theorem for (TE) structures.

In order to show that 
the reconstructed $\hmEh$ is isomorphic to the big QDM of $M$, 
we use a mirror theorem by Coates and Givental \cite{coates-givental}. 
Coates and Givental introduced an infinite dimensional Lagrangian cone 
in the symplectic space 
$H^*(M)\otimes \C[\hbar,\hbar^{-1}]\!][\![Q]\!]$ 
encoding all the information of the genus zero Gromov-Witten theory. 
They described a relationship between the Gromov-Witten 
theory of $M$ itself and 
the twisted theory by a vector bundle $\mV$ on $M$ 
as a symplectic transformation of the corresponding Lagrangian cones. 
We interpret Coates--Givental's symplectic 
transformation in terms of $D$-modules.  
In the framework of AQDMs, a symplectic transformation 
corresponds to a change of a frame at $Q=0$ 
and and a shift of the origin 
(Proposition \ref{prop:symplectictransf_as_gaugetransf}).  

In generalized mirror transformations, the use of
(non-convergent) formal power series 
in the Novikov ring parameters is inevitable. 
In a forthcoming paper \cite{iritani-convergent}, 
we will discuss the convergence of generalized mirror transformations 
in some refined sense. 

The paper is organized as follows. 
In section 2, we review the theory of quantum $D$-modules. 
In section 3, we review the equivariant Floer cohomology for 
toric complete intersections. 
In section 4, we formulate abstract big/small quantum $D$-modules 
and prove the reconstruction theorem. 
In section 5, we give a proof of the generalized mirror transformations. 
We also include the review of Coates--Givental's theory.
In section 6, we illustrate generalized mirror transformations by examples. 

\noindent 
{\bf Acknowledgments}
Thanks are due to Professor Hiraku Nakajima for  
his encouragement and guidance. 
The author also expresses gratitude to Professor Masao Jinzenji 
for explaining his works and his computer programs. 
He is also grateful to Kazushi Ueda for valuable discussions. 
This research is partially supported by JSPS Fellows and 
Scientific Research 15-5482. 

\section{Big and small QDMs} 
\label{sec:QDM} 
In \cite{coates-givental}, Coates and Givental introduced 
a twist of Gromov-Witten theory by a vector bundle and 
a multiplicative characteristic class. 
In this section, we review the genus zero twisted theory 
from a viewpoint of the $D$-module structure. 
We explain fundamental solutions, $J$-functions 
and Kontsevich and Manin's reconstruction theorem. 

\subsection{Twisted quantum cohomology} 
Let $M$ be a smooth projective variety 
and $\mV$ be a vector bundle over $M$. 
For simplicity, we assume that the total cohomology ring of $M$ 
consists only of the even degree part, $H^*(M)=H^{\rm even}(M,\C)$. 
Following \cite{coates-givental}, we introduce 
the following general multiplicative characteristic class $\bc$ 
for a vector bundle $F$: 
\begin{equation}
\label{eq:general_c}
\bc(F) := \exp\left(\sum_{k\ge 0} s_k \ch_k(F)\right)
\end{equation}
where $s_0,s_1,s_2,\dots$ are arbitrary parameters. 
Let $\C[\![\bs]\!]$ denote the completion of  
$\C[s_0,s_1,\dots]$ with respect to the valuation 
$v\colon \C[s_0,s_1,\dots] \rightarrow \R\cup\{\infty\}$ defined by 
\begin{equation}
\label{eq:valofCs}
v(s_k)=k+1, \quad k\ge 0. 
\end{equation}
Then $\bc$ is defined over $\C[\![\bs]\!]$. 
Let $\ov{M}_{0,n}(M,d)$ be the moduli space of 
genus zero, degree $d$ stable maps to $M$ with $n$ marked points. 
We have the following structure maps: 
\[
\begin{CD}
\ov{M}_{0,n+1}(M,d) @>{e_i}>> M \\
@V{\pi_i}VV & @. \\
\ov{M}_{0,n}(M,d)
\end{CD}
\]
where $\pi_i$ is the $i$-th forgetful map and 
$e_i$ is the evaluation map at the $i$-th marked point. 
The genus zero $(\mV,\bc)$-{\it twisted} Gromov-Witten invariants 
are defined by 
\begin{equation}
\label{eq:twisted_GW_inv}
\corrV{\psi_1^{k_1} \alpha_1,\dots  
          ,\psi_n^{k_n} \alpha_n}_d
  :=\int_{[\ov{M}_{0,n}(M,d)]^{\rm virt}}
  \bc (R^\bullet{\pi_{n+1}}_*e_{n+1}^*\mV) \cup 
  \prod_{i=1}^n e_i^*(\alpha_i) \psi_i^{k_i} 
\end{equation} 
where $\psi_i$ is the first Chern class of the $i$-th cotangent line, 
$\alpha_1,\dots,\alpha_n\in H^*(M)$ 
and  $[\ov{M}_{0,n}(M,d)]^{\rm virt}$ is the virtual fundamental class. 
These correlators take values in $\C[\![\bs]\!]$. 
We also introduce the twisted Poincar\'{e} pairing on $H^*(M)$:
\begin{equation}
\label{eq:twisted_Poincare}
\pairV{\alpha_1}{\alpha_2} := 
\int_M \alpha_1 \cup \alpha_2 \cup \bc(\mV).    
\end{equation}
The twisted Gromov-Witten invariants satisfy almost all 
the formal properties of the ordinary Gromov-Witten invariants ---
the string equation, the dilaton equation, 
the topological recursion relations (TRR) and 
the divisor equation (see for instance 
\cite[Section 1]{pandharipande}, \cite{givental-symplecticgeom}) 
--- but the pairing appearing in the TRR 
should be replaced with the twisted Poincar\'{e} pairing. 
In \cite{givental-symplecticgeom}, 
Givental explained these properties by 
geometry of a Lagrangian cone 
(see also Section \ref{subsec:review_QL}).

Let $p_0$ be the unit class $1\in H^0(M,\C)$ and 
$\{p_1,\dots,p_r\}$ be an integral basis of 
$H^2(M,\Z)_{\rm free}=H^2(M,\Z)/H^2(M,\Z)_{\rm tor}$. 
We can choose a basis so that each 
$p_a$ $(1\le a\le r)$ is a nef class, \emph{i.e.} 
$p_a$ intersects all effective curve classes non-negatively.  
Let $\{p_{r+1},\dots,p_s\}$ be a basis of $H^{\ge 4}(M,\C)$. 
Let $t^0,t^1,\dots,t^r,t^{r+1},\dots,t^s$ be linear co-ordinates 
of $H^*(M)$ dual to the basis $\{p_0,p_1,\dots,p_r,p_{r+1},\dots,p_s\}$. 
Let $\tau=\sum_{i=0}^s t^i p_i$ be a general point on $H^*(M)$. 
Let $Q^a$, $1\le a\le r$ be the Novikov variable dual to 
$p_a\in H^2(M,\Z)_{\rm free}$. 
This gives the following co-ordinate on $H^2(M,\C^*)$: 
\[
Q^a\colon H^2(M,\C^*) \cong H^2(M,\C/2\pi\sqrt{-1}\Z) \ni \left[\sum_{a=1}^r t^a p_a\right] 
\longmapsto \exp(t^a) \in \C^*.  
\] 
For $d\in H_2(M,\Z)$, we denote by $Q^d$ the monomial 
$(Q^1)^{\pair{p_1}{d}}(Q^2)^{\pair{p_2}{d}}\dots (Q^r)^{\pair{p_r}{d}}$.
The $(\mV,\bc)$-twisted big and small quantum cohomology 
$QH^*_{\bc}(M,\mV)$, $SQH^*_{\bc}(M,\mV)$ 
are tensor products of the cohomology ring 
and the formal power series ring in co-ordinate variables: 
\begin{align*}
QH^*_{\bc}(M,\mV)&=H^*(M)\otimes
\C[\![\bs]\!][\![Q^1,\dots,Q^r]\!][\![t^0,t^1,\dots,t^s]\!], 
\quad &\text{(big)} 
\\
SQH^*_{\bc}(M,\mV)&
=H^*(M)\otimes \C[\![\bs]\!][\![Q^1,\dots,Q^r]\!]. 
\quad &\text{(small)} 
\end{align*}
The quantum products $*$ on these modules are 
bilinear over the formal power series ring and 
deformations of the cup product $\cup$. 
The product of $SQH^*_{\bc}(M,\mV)$ is obtained as the specialization 
$t^0=t^1=\cdots=t^s=0$ of that of $QH^*_{\bc}(M,\mV)$.  
We define the big quantum product $*$ of $(M,\mV)$ by the formula  
\begin{equation}
\label{eq:defofquantumprod}
\pairV{p_i *p_j}{p_k}
  =\sum_{d\in \Lambda}\sum_{n\ge 0}\frac{1}{n!}
    \corrV{p_i,p_j,p_k,
    \overbrace{\tau,\dots,\tau}^{\text{$n$ times }}}_d Q^d 
\end{equation} 
where $\tau=\sum_{i=0}^s t^i p_i$ and 
$\Lambda\subset H^2(M,\Z)$ denotes the Mori cone, \emph{i.e.} 
the semigroup generated by effective curves.  Because $p_a$ is nef, 
$Q^d$ in the summation does not contain negative powers of 
$Q^a$, \emph{i.e.} $Q^d\in \C[Q^1,\dots,Q^r]$. 

\begin{remark*}
By applying the divisor equation, we get 
\[
\pairV{p_i*p_j}{p_k} =
\sum_{d\in \Lambda} \sum_{n\ge 0} \frac{1}{n!} 
\corrV{p_i,p_j,p_k,
    \overbrace{\tau_{\ge 4},\dots,\tau_{\ge 4}}^{\text{$n$ times }}}_d 
    e^{\pair{\tau}{d}} Q^d,  
\]
where we write $\tau= t^0 p_0 + \sum_{a=1}^r t^a p_a + \tau_{\ge 4}$ with 
$\tau_{\ge 4} \in H^{\ge 4}(X,\C)$.  
Thus for $\tau= \sum_{a=1}^r t^a p_a \in H^2(X,\C)$, we have 
\[
\pairV{p_i*p_j}{p_k} =
\sum_{d\in \Lambda}  
\corrV{p_i,p_j,p_k}_d 
\prod_{a=1}^r (e^{t^a} Q^a)^{\pair{p_a}{d}}.   
\] 
This shows that the small quantum product (restriction of the big one 
to $t^i=0$ for $0\le i\le s$) 
is equivalent to the big quantum product restricted to $\tau\in H^2(X,\C)$. 
The big quantum product on $H^2(X,\C)$ is obtained by 
substituting $Q^a$ in the small quantum product with $e^{t^a}Q^a$. 
\end{remark*} 


As an important example of the twists, 
we consider the twist by the $S^1$-equivariant Euler class $\be$: 
\begin{equation}
\label{eq:equivariant_Eulerclass}
\be(F) := \sum_{i\ge 0} \lambda^{\rank(F)-i} c_i(F)
\end{equation}
where the $S^1$-action on $F$ is defined by  
the scalar multiplication on each fiber and 
$\lambda$ is a generator of $H_{S^1}^*({\rm pt})$.  
This corresponds to the choice of parameters: 
\[
s_0 = \log \lambda, \quad s_k = (-1)^{k-1} (k-1)!/\lambda^k.  
\]
In this case, the quantum product is defined on the following modules.  
\begin{align*}
QH^*_{\be}(M,\mV)&=H^*(M)\otimes
\C(\!(\lambda^{-1})\!)[\![Q^1,\dots,Q^r]\!][\![t^0,t^1,\dots,t^s]\!], 
\quad &\text{(big)} 
\\
SQH^*_{\be}(M,\mV)&
=H^*(M)\otimes \C(\!(\lambda^{-1})\!)[\![Q^1,\dots,Q^r]\!]. 
\quad &\text{(small)} 
\end{align*}
Now assume that the bundle $\mV$ is {\it convex}, \emph{i.e.}  
for any holomorphic map $f\colon\Proj^1\rightarrow M$,
$H^1(\Proj^1,f^*(\mV))=0$ holds. 
In this case, the structure constants of $QH_\be^*(M,\mV)$ 
take values in $\C[\lambda][\![Q^1,\dots,Q^r]\!][\![t^0,\dots,t^s]\!]$, 
\emph{i.e.} we need not invert the variable $\lambda$. 
We refer the reader to \cite{pandharipande} for details. 
Thus we can define the quantum product on the smaller modules:
\begin{align*}
QH^*_\be(M,\mV) &= 
H^*(M)\otimes \C[\lambda][\![Q^1,\dots,Q^r]\!][\![t^0,\dots,t^s]\!] 
\quad & \text{(big)}\\ 
SQH^*_\be(M,\mV) &= 
H^*(M)\otimes \C[\lambda][\![Q^1,\dots,Q^r]\!] 
\quad & \text{(small)}
\end{align*}  
In particular, we can define the non-equivariant ($\lambda=0$) 
Euler twisted quantum cohomology for a pair $(M,\mV)$. 
The non-equivariant version is closely related 
to the quantum cohomology of zero-locus $N\subset M$ 
of a transverse section of $\mV$ 
as the following theorem shows. 

\begin{theorem}[\cite{kim-kresch-pantev}] 
\label{thm:kkp}
Let $\mV$ be a convex vector bundle on $M$ and 
$\iota\colon N\hookrightarrow M$ be a smooth subvariety 
defined by a regular section of $\mV$. 
Then 
\[
[\ov{M}_{0,n}(M,d)]^{\rm virt} \cap 
\Euler(R^\bullet({\pi_{n+1}}_* e_{n+1}^* \mV)) 
= \sum_{d=\iota_*\beta} \iota_* [\ov{M}_{0,n}(N,\beta)]^{\rm virt} 
\]
for $d\in H_2(M,\Z)$. 
In particular, we have 
\[
\lim_{\lambda\to 0} \pairV{p_i*_{\mV}p_j}{p_k}(\tau)
=\pair{\iota^*p_i*_N \iota^*p_j}{\iota^*p_k}^N(\iota^*\tau)
\Bigr|_{H_2(N)\to H_2(M)}. 
\] 
Here, $*_\mV$ and $*_N$ are the products of 
$QH^*_{\be}(M,\mV)$ and $QH^*(N)$ respectively,  
$\tau\in H^*(M)$ and $\pair{\cdot}{\cdot}^N$ is the Poincar\'{e} 
pairing of $N$.
The notation $|_{H_2(N)\to H_2(M)}$ means to replace 
$Q^d$, $d\in H_2(N)$ with $Q^{i_*(d)}$. 
\end{theorem}

\begin{notation*}
As we have seen, we can define the quantum cohomology over 
various ground rings.  
In order to describe different versions in a unified way, 
we use the letter $K$ to denote the ground ring. 
In Section \ref{sec:QDM}, $K$ means 
\begin{equation}
\label{eq:choiceofK}
K =\begin{cases} 
\C & \text{untwisted}, \\ 
\C[\lambda] & \text{$\bc=\be$, convex $\mV$}, \\ 
\C(\!(\lambda^{-1})\!) & \text{$\bc=\be$, general $\mV$}, \\
\C[\![\bs]\!] & \text{general $\bc$}.  
\end{cases}  
\end{equation} 
Also we use the following shorthand:  
\[
K[\![Q,t]\!] := K[\![Q^1,\dots,Q^r]\!][\![t^0,\dots,t^s]\!], \quad 
K[\![Q]\!] := K[\![Q^1,\dots,Q^r]\!]. 
\]
\end{notation*}

\subsection{QDM} 
\label{subsec:QDM}
Let $K$ be a topological ring $K$ 
endowed with a valuation $v\colon K\rightarrow
\R\cup\{\infty\}$. 
The space 
$K\{\hbar,\hbar^{-1}\}\!\}$ of Laurent power series 
is defined to be
\[
K\{\hbar,\hbar^{-1}\}\!\} := 
\left\{
\sum_{n\in \Z} a_n \hbar^n \;\Big|\;
a_n\in K, \ \lim_{n\to \infty} v(a_n) =\infty, \ 
\inf_{n\in \Z} v(a_n)> -\infty \right\}. 
\]
We define $K\{\hbar\}$ (resp.\ $K\{\!\{\hbar^{-1}\}\!\}$) 
to be the subspace of $K\{\hbar,\hbar^{-1}\}\!\}$ consisting of the 
positive (resp.\ negative) power series in $\hbar$. 
These become rings when $K$ is complete. 
The valuation on $K=\C[\![\bs]\!]$ was defined in (\ref{eq:valofCs}) 
and that on $K=\C(\!(\lambda^{-1})\!)$ 
is given by $v(\sum_{n}a_n\lambda^{-n})=\min\{n|a_n\neq 0\}$. 
These are complete. 
The valuations on $K=\C[\lambda]$ and $\C$ 
are induced from that on $\C(\!(\lambda^{-1})\!)$ and 
give them the discrete topology.  
For example, we have 
\begin{gather*}
\C(\!(\lambda^{-1})\!)\{\hbar,\hbar^{-1}\}\!\} = 
\C(\!(\hbar^{-1})\!)(\!(\lambda^{-1})\!), \quad 
\C(\!(\lambda^{-1})\!)\{\hbar\} = \C[\hbar](\!(\lambda^{-1})\!), \\
\C[\lambda]\{\hbar,\hbar^{-1}\}\!\} 
=\C(\!(\hbar^{-1})\!)[\lambda], \quad 
\C[\lambda]\{\hbar\}=\C[\hbar,\lambda].  
\end{gather*} 

\begin{definition} 
The big and small QDMs are modules  
\begin{align*}
QDM_{\bc}(M,\mV) &=H^*(M)\otimes 
K\{\hbar\}[\![Q,t]\!], \quad \text{(big)} \\
SQDM_{\bc}(M,\mV) &= H^*(M)\otimes K\{\hbar\}[\![Q]\!],  
\quad\text{(small)}   
\end{align*}
endowed with the actions of the Dubrovin connection: 
\begin{align*}
\nabla^\hbar_i &= \hbar \parfrac{}{t^i} + p_i * \quad 
\text{($0\le i\le s$, defined only for $QDM_{\bc}(M,\mV)$)}, \\
\smnabla^\hbar_a &= \hbar Q^a \parfrac{}{Q^a} + p_a * \quad 
(1\le a\le r)  
\end{align*}
where $*$ denotes the $S^1$-equivariant 
big/small quantum product of $(M,\mV)$ defined in (\ref{eq:defofquantumprod}).   
The Dubrovin connection is known to be flat. 
\end{definition}

The big QDM is considered as a flat bundle over  
a formal neighborhood of the origin in $\C^r \times H^*(M,\C)$. 
Here, $\C^r$ is a partial compactification of $H^2(M,\C^*)$ given 
by a choice of co-ordinates $Q^1,\dots, Q^r$. The small QDM 
is the restriction of the big one to $\C^r\times \{0\}$. 
The flat connection has a logarithmic singularity along $Q^1\cdots Q^r =0$. 
By the divisor equation, 
the big QDM restricted on the locus $\C^r\times H^2(M,\C)$ 
can be recovered from the small one. 
We also consider the small QDM as a module over 
the Heisenberg algebra $\mD$: 
\[
\mD=K\{\hbar\}[\![Q^1,\dots,Q^r]\!] 
     \langle \Qp_1,\dots,\Qp_r\rangle
\] 
whose generators satisfy the following commutation relations: 
\[
[\Qp_a,Q^b]=\hbar\delta_a^b Q^b,\quad [\Qp_a,\Qp_b]=[Q^a,Q^b]=0.
\]
Here, $\mD$ acts on $SQDM_{\bc}(M,\mV)$ by 
\[
\Qp_a \longmapsto \smnabla^\hbar_a, \quad 
Q^a \longmapsto \text{multiplication by } Q^a.
\]

\subsection{Fundamental solution and $J$-function} 
A fundamental solution to the big QDM is 
a formal section $L$ of the endomorphism bundle satisfying 
\begin{align}
\nonumber
&L\in \End(H^*(M))\otimes 
K\{\hbar,\hbar^{-1}\}\!\}[\![Q,t]\!], 
\\  
\label{eq:fundsol_big}
& \nabla^\hbar_j\circ L = L \circ \hbar \parfrac{}{t^j}, \quad (0\le j\le s) \\ 
\label{eq:fundsol_small} 
& \smnabla^\hbar_a \circ L = L \circ (\hbar Q^a \parfrac{}{Q^a} +p_a\cup), 
\quad (1\le a\le r).   
\end{align}
An explicit form of a fundamental solution is given by 
the gravitational descendants 
(see \cite[Section 1.3]{pandharipande}. 
Note that our choice of the sign of $\hbar$ is opposite to \cite{pandharipande}.) 
\begin{equation}
\label{eq:fundsol_grav}
\pairV{Lp_i}{p_j}:=\pairV{p_i}{p_j}+
\sum_{d\in \Lambda, n\ge 0,(d,n)\neq (0,0)}
\frac{1}{n!}
\corrV{\frac{p_i}{-\hbar-\psi_1},p_j,
\overbrace{\tau,\dots,\tau}^{\text{$n$ times}}}_d Q^d. 
\end{equation} 
The fraction $\frac{1}{-\hbar-\psi_1}$ here should be expanded 
as a power series in $\hbar^{-1}$. 
This $L$ does not contain positive powers in $\hbar$ and belongs to 
$\End(H^*(X))\otimes K\{\!\{\hbar^{-1}\}\!\}[\![Q,t]\!]$.   
Using the TRR, one can show that $Lp_i$, 
$0\le i\le s$ form a basis of $\nabla^\hbar$-parallel sections 
(see \cite[Proposition 2]{pandharipande}). 
Thus the equation (\ref{eq:fundsol_big}) holds for this $L$. 
By using the divisor equation, 
we can rewrite the above $L$ as (see \cite{pandharipande})
\begin{align*}
\pairV{Lp_i}{p_j}&=\pairV{e^{-\tau/\hbar}p_i}{p_j} \\ 
& +\sum_{d\in \Lambda\setminus\{0\}}\sum_{n\ge 0}
 \frac{1}{n!}
\corrV{\frac{e^{-(t^0+\tau_2) /\hbar} p_i}{-\hbar-\psi_1},p_j,
 \overbrace{\tau_{\ge 4},\dots,\tau_{\ge 4}}^{\text{$n$ times}}}_d 
 \prod_{a=1}^r (e^{t^a}Q^a)^{\pair{p_a}{d}}, 
\end{align*}
where $\tau= \sum_{i=0}^s t^ip_i = t^0p_0 + \tau_2 + \tau_{\ge 4}$ with 
$\tau_2\in H^2(M)$ and $\tau_{\ge 4} \in H^{\ge 4}(M)$. 
In the expression like $e^{-\tau/\hbar}$, $\tau$ is considered 
to be an operator acting on cohomology by the cup product.   
Therefore, we can decompose $L$ in the form 
$L=S\circ e^{-(t^0+\tau_2)/\hbar}$, where 
$S$ is an element of 
$\End(H^*(M))\otimes K\{\!\{\hbar^{-1}\}\!\}[\![Q,t]\!]$.  
It follows from the above expression that $S$ satisfies 
\begin{equation}
\label{eq:divisoreq_for_S}
\parfrac{}{t^a} S = Q^a \parfrac{}{Q^a}S.  
\end{equation}
The equation (\ref{eq:fundsol_small}) 
follows from (\ref{eq:fundsol_big}) and (\ref{eq:divisoreq_for_S}). 
\begin{proposition} 
The above fundamental solution $L=L(Q,\tau,\hbar)$ 
given by gravitational descendants 
is characterized by the following condition: 

i) initial condition:  $L(0,0,\hbar)=\id$ and 

ii) differential equations: 
\begin{align*}
 &\hbar\parti{a}L-L\circ(p_a\cup)+(p_a*)\circ L=0,\quad (1\le a\le r)\\
 &\hbar\parfrac{}{t^j} L+ (p_j*)\circ L=0, \quad (0\le j\le s).
\end{align*}
Moreover this satisfies  


iii) unitarity: 
\[
\pairV{L(Q,\tau,-\hbar)p_i}{L(Q,\tau,\hbar)p_j}=
\pairV{p_i}{p_j},
\]

iv) divisor equation: 
\begin{equation}
\label{eq:divisor}
(\parfrac{}{t^a} - Q^a\parfrac{}{Q^a}) L + L \circ (\frac{p_a}{\hbar}\cup) =0.
\end{equation} 
\end{proposition} 
\begin{proof} 
The unitarity is stated in \cite{givental-elliptic}. 
Set $\ov{L}=L|_{\hbar \mapsto -\hbar}$. 
By the Frobenius property $\pairV{p_i*p_j}{p_k}=\pairV{p_i}{p_j*p_k}$ 
and the differential equation for $L$, 
we have  
\begin{align*} 
\hbar \parfrac{}{t^k} \pairV{\ov{L}p_i}{Lp_j} & =
\pairV{p_k*\ov{L}p_i}{Lp_j} - \pairV{\ov{L}p_i}{p_k*Lp_j} =0, \\  
\hbar Q^a\parfrac{}{Q^a} \pairV{\ov{L}p_i}{Lp_j} 
&= -\pairV{\ov{L}(p_a\cup) p_i}{Lp_j} + \pairV{\ov{L}p_i}{L(p_a\cup) p_j}.    
\end{align*}
Since the operation $p_a\cup$ is nilpotent, 
$(\hbar Q^a\parfrac{}{Q^a})^n \pairV{\ov{L}p_i}{Lp_j}$ is zero 
for a sufficiently big $n$. 
This shows that $\pairV{\ov{L}p_i}{Lp_j}$ is a constant 
and is equal to $\pairV{\ov{L}p_i}{Lp_j}|_{Q=\tau=0}=\pairV{p_i}{p_j}$. 
The divisor equation follows from (\ref{eq:divisoreq_for_S}). 
\end{proof} 

The unitarity means that $L(Q,\tau,\hbar)^{-1}$ is 
the adjoint of $L(Q,\tau,-\hbar)$. Put $g_{ij}=\int_{M} p_i\cup p_j$ and 
let $(g^{ij})_{0\le i,j\le s}$ be the matrix inverse to 
$(g_{ij})_{0\le i,j\le s}$. 
Using the unitarity, we calculate $L^{-1}$ as 
\begin{align}
\label{eq:J-function}
&L^{-1}p_i= p_i + \sum_{d\in \Lambda, n\ge 0, (d,n)\neq (0,0)} \sum_{j,k} 
\frac{1}{n!} \frac{p_k}{\bc(\mV)}g^{kj}
\corrV{\frac{p_j}{\hbar-\psi_1},p_i,
\overbrace{\tau,\dots,\tau}^{\text{$n$ times}} }_d Q^d \\
\nonumber
&=e^{\tau/\hbar}p_i+
  e^{(t^0+\tau_2)/\hbar}\sum_{d\in\Lambda\setminus \{ 0\} } 
  \sum_{n\ge 0}\sum_{j,k}\frac{1}{n!}\frac{p_k}{\bc(\mV)}g^{kj}
  \corrV{\frac{p_j}{\hbar-\psi_1},p_i,
  \overbrace{\tau_{\ge 4},\dots,\tau_{\ge 4}}^{\text{$n$ times}}}_d 
  e^{\pair{\tau}{d}}Q^d. 
\end{align}
\begin{definition}
\label{def:J-funct}
The {\it big $J$-function} of $(M,\mV)$ is a 
cohomology-valued formal function 
defined as $J:=L^{-1}p_0$, where $L$ is the fundamental solution 
in (\ref{eq:fundsol_grav}) and 
$p_0$ is the unit.  
This is an element of 
$H^*(M)\otimes K\{\!\{\hbar^{-1}\}\!\}[\![Q,t]\!]$. 
The {\it small $J$-function} is the restriction of the big one to 
$\tau=0 \in H^*(M)$.  
\end{definition}

\begin{remark}
\label{rem:differenceofJ}
The small $J$-function is slightly different from what the author used 
in the previous paper \cite{iritani-EFC}. 
The $J$-function $J_{\rm prev}(Q,\hbar)$ in \cite{iritani-EFC}
is related to the above small $J$-function $J_{\rm small}(Q,\hbar)$ by 
\[
J_{\rm prev}(Q,\hbar) = 
e^{\sum_{a=1}^r p_a\log Q_a/\hbar} J_{\rm small}(Q,\hbar),  
\]
where $J_{\rm prev}$ lives in 
$H^*(X)\otimes K\{\!\{\hbar^{-1}\}\!\}[\![Q]\!][\log Q^1,\dots,\log Q^r]$. 
It follows from \eqref{eq:J-function} that the big $J$-function 
restricted to $\tau_2\in H^2(X)$ is related to the small $J$-function by 
\[
J_{\rm big}(Q,\tau_2,\hbar) = e^{\tau_2/\hbar} 
J_{\rm small}(e^{t^1}Q^1,\dots,e^{t^r}Q^r,\hbar), \quad 
\tau_2=\sum_{a=1}^r t^a p_a\in H^2(X).  
\]
\end{remark} 
The big QDM is generated by the unit section $p_0$ 
as a $D$-module. In fact, $QDM_{\bc}(M,\mV)$ is generated by $p_0$ 
and its derivatives by the Dubrovin connection:  
\[
\nabla^\hbar_1 p_0 = p_1, \dots, \nabla^\hbar_s p_0 = p_s. 
\] 
as a $K\{\hbar\}[\![Q,t]\!]$-module. 
Because the big $J$-function is defined to be the inverse image of $p_0$ 
under $L$, it plays the role of a
generator of the big QDM. 
In fact, if a differential operator $P$ annihilates $J$: 
\[
P( Q,t, \hbar Q^a\parfrac{}{Q^a}+p_a, \hbar\parfrac{}{t^i},\hbar) 
J(Q,\tau,\hbar) =0, 
\]
then we have a relation in $QDM_{\bc}(M,\mV)$: 
\[
P(Q,t,\smnabla_a^\hbar, \nabla_i^\hbar,\hbar) p_0 =0    
\]
and \emph{vice versa}. 
On the other hand, the small QDM is generated by 
the unit section under the $H^2$-generation condition;
in this case, the $J$-function plays the role of a generator of 
the small QDM. 
\begin{proposition}[{\cite[Theorem 2.4]{iritani-EFC}}]
\label{prop:Jgenerates}
Let $M$ be a smooth projective variety and 
$\mV$ be a vector bundle on $M$. 
Assume that the total cohomology ring of $M$ is generated by the
second cohomology group. 
Then the small QDM is 
generated by the unit section $p_0$ as a $\mD$-module. 
Moreover, we have an isomorphism of $\mD$-modules 
$\mD/\frI\cong SQDM_{\bc}(M,\mV)$ which sends $1$ to $p_0$,   
where $\frI$ is the left ideal of $\mD$ 
consisting of elements $f(Q,\Qp,\hbar)\in \mD$ satisfying 
\[
f(Q^1,\dots,Q^r,\hbar Q^1\parfrac{}{Q^1}+p_1, \dots, 
\hbar Q^r\parfrac{}{Q^r}+p_r,\hbar) J_{\rm small}(Q,\hbar) = 0.  
\]
\end{proposition}
\begin{proof}
This theorem was shown in \cite[Theorem 2.4]{iritani-EFC} 
for convex $\mV$. 
The proof applies without change to this more general case. 
\end{proof} 

In this paper, we are mainly interested in the quantum cohomology 
of toric varieties and their twists. 
In this case, the $H^2$-generation of the total cohomology always holds. 
When the total cohomology ring is generated by $H^2(M)$, 
we have the following reconstruction theorem by 
Kontsevich and Manin \cite{kontsevich-manin}, 
whose generalization is the main theme of this paper.  
\begin{theorem}[Kontsevich and Manin]
\label{thm:KMreconstruction}
If the total cohomology ring $H^*(M)$ is generated by 
the second cohomology group $H^2(M)$ as a ring, 
the big quantum cohomology $QH^*(M)$ can be 
reconstructed from the small quantum cohomology $SQH^*(M)$. 
In other words, if we know all the three point functions 
of the form $\langle p_a,p_i, p_j \rangle_{0,3,d}$  
with $1\le a\le r$, $1\le i,j\le s$, 
we can reconstruct all genus $0$, $n$-point functions for all $n$.  
\end{theorem}

\section{Equivariant Floer Theory}
\label{sec:EFT}
In this section, we review the equivariant Floer theory 
for toric complete intersections 
\cite{givental-ICM,givental-homologicalgeom,iritani-EFC}. 
The exposition here slightly differs from 
what was given in \cite{iritani-EFC} 
so that it contains a little more general case where 
$\mV$ is not necessarily nef.  
See also \cite{audin} for background materials 
on toric varieties.

\subsection{Toric varieties} 
\label{subsec:toric}
A projective toric variety $X$ can be defined 
by the following data: 
\begin{itemize}
\item[(1)] algebraic torus $\T_\C\cong (\C^*)^r$ 
and its maximal compact subgroup $\T \cong (S^1)^r$;  
\item[(2)] $N$-tuple of integral vectors $u_1,\dots,u_N$ 
in the weight lattice $\Hom(\T_\C,\C^*)$; 
\item[(3)] a vector $\eta$ in $\Lie(\T)^\vee:=\Hom(\T_\C,\C^*)\otimes \R$. 
\end{itemize}
The integral vectors $u_1,\dots,u_N$ define a homomorphism 
$\T_\C\rightarrow (\C^*)^N$ ($\T\rightarrow (S^1)^N$) and 
an action of $\T_\C$ (resp. $\T$) on $\C^N$.   
Define 
\[
\mA_\eta:= \left\{ I \subset \{1,2,\dots,N\} \;\Big|\; 
\sum_{i\in I} \R_{>0} u_i \ni \eta
\right \}.  
\]
This is a subset of the power set $\mP(\{1,2,\dots,N\})$. 
Then $X$ is defined as a GIT quotient of $\C^N$ by $\T_\C$. 
\[
X:=\C^N/\!/_\eta \T_{\C}= \mU_\eta/\T_\C, \quad 
\mU_\eta:=\C^N \setminus \bigcup_{I\notin \mA_\eta} \C^I,  
\]
where $\C^I\subset \C^N$ denote the co-ordinate subspace 
$\{(z_1,\dots,z_N); z_i=0 \text{ for } i\notin I\}$. 
The toric variety $X$ is non-empty, smooth and compact 
under the following conditions. 
\begin{itemize}
\item[(A)] $\{1,\dots, N\}\setminus\{i\} \in \mA_\eta$ 
for all $1\le i\le N$.  
\item[(B)] For any $I\in \mA_\eta$, $\{u_i\}_{i\in I}$ generate $\Hom(\T_\C,\C^*)$ 
as a $\Z$-module.
\item[(C)] If $\sum_{i=1}^N c_i u_i =0$ with $c_i\ge 0$, then $c_j=0$ for all $j$. 
\end{itemize}
Henceforth we assume these. A representation $\rho\in \Hom(\T_\C,\C^*)$ of $\T_\C$ 
gives rise to a line bundle $\mL_\rho$: 
\begin{equation}
\label{eq:linebundleonX}
\mL_\rho := \C \times_{\rho} \mU_\eta \longrightarrow X,   
\end{equation} 
where $\T_\C$ acts on $\C\times \C^N$ by 
the representation $\rho$ on the $\C$ factor and in the given way 
on the $\C^N$ factor. 
By this, we have a natural identification:  
\[
\Hom(\T_\C, \C^*) \cong \Pic(X) \overset{c_1}{\cong} H^2(X,\Z).  
\]
Let $\mu\colon \C^N \rightarrow \Lie(\T)^\vee$ be the moment map for 
the $\T$ action on $\C^N$. This is given by 
$\mu(z_1,\dots, z_N) = \sum_{i=1}^N \frac{1}{2} |z_i|^2 u_i$. 
We can give $X$ also as a symplectic quotient: 
\[
X \cong \mu^{-1}(\eta)/\T. 
\]
The vector $\eta\in \Lie(\T)^\vee$ is identified 
with the class of the reduced symplectic form. 
The K\"{a}hler cone $C_X$ of $X$ 
\emph{i.e.} the cone of K\"{a}hler classes 
are given by 
\begin{equation}
\label{eq:Kaehlercone}
C_X = \bigcap_{I\in \mA_\eta} \sum_{i\in I} \R_{>0} u_i \subset 
\Hom(\T_\C, \C^*)\otimes \R \cong H^2(X,\R).  
\end{equation} 
We have $\eta\in C_X$. 
Take a nef integral basis $p_1,\dots, p_r$ of $H^2(X,\Z)\cong \Hom(\T_\C,\C^*)$ 
as we did in Section \ref{sec:QDM}. 
Here $p_a\in \ov{C}_X$ since $p_a$ is nef. 
We write 
\begin{equation}
\label{eq:classesu}
u_i = \sum_{a=1}^r m_i^a p_a 
\end{equation}
The vector $u_i\in \Hom(\T_\C,\C^*)$ represents the Poincar\'{e} dual 
of the torus invariant divisor $\{z_i=0\}$, where $z_i$ is the 
standard $i$-th co-ordinate on $\C^N$.  

We will also consider the equivariant Floer cohomology of $X$ twisted 
by a vector bundle $\mV$. We take $\mV$ as a sum 
$\mV=\mV_1\oplus\mV_2\oplus \dots \oplus \mV_l$ of line bundles $\mV_i$.  
We put 
\begin{equation}
\label{eq:classesv}
v_i := c_1(\mV_i) =\sum_{a=1}^r l_i^a p_a \in H^2(X,\Z).  
\end{equation} 

\subsection{Algebraic model for the universal cover of free loop spaces} 

Following Givental \cite{givental-homologicalgeom} and 
Vlassopoulos \cite{vlassopoulos}, 
we introduce the algebraic model $L_X$ for the universal cover $\widetilde{X^{S^1}}$ 
of the free loop space $X^{S^1}=\Map(S^1,X)$. 
The universal cover $\widetilde{X^{S^1}}$ is given by 
\[
\widetilde{X^{S^1}} = \Map(S^1,\mU_\eta)/ \Map_0(S^1,\T_\C),   
\]
where $\Map_0(S^1,\T_\C)$ denotes the set of contracting loops in $\T_\C$.  
We will replace $\Map(S^1,\mU_\eta)$ with 
the space of Laurent polynomial loops and $\Map_0(S^1,\T_\C)$ with $\T_\C$
to construct the algebraic model $L_X$. 
Define  
\[
L_X:=\C[\zeta,\zeta^{-1}]^N/\!/_\eta \T_{\C} = L\mU_\eta/\T_\C, \quad 
L\mU_\eta :=\C[\zeta,\zeta^{-1}]^N\setminus 
\bigcup_{I\notin\mA_\eta} \C[\zeta,\zeta^{-1}]^I. 
\]
Here, $\zeta$ is considered to be a parameter of loops. 
The $\T_\C$ action on $\C[\zeta,\zeta^{-1}]^N \cong \C^N\otimes \C[\zeta,\zeta^{-1}]$ 
is induced from the given action on $\C^N$ and the trivial action on 
$\C[\zeta,\zeta^{-1}]$. 
We write a general point on $L_X$ as 
\[
[z_1(\zeta),\dots, z_N(\zeta)] \in L_X, \quad \text{where }  
z_i(\zeta)=\sum_{\nu\in \Z} a_{i\nu}\zeta^\nu \in \C[\zeta,\zeta^{-1}].   
\]
The variables $\{a_{i\nu}\}$ play the role of homogeneous co-ordinates on $L_X$. 
We can write $L_X$ also as a symplectic quotient. 
The moment map 
$\mu_\infty\colon \C[\zeta,\zeta^{-1}]^N \rightarrow \Lie(\T)^\vee$ 
of $\T$ action is given by 
$\mu_{\infty}[z_1(\zeta),\dots,z_r(\zeta)]=\sum_{1\le i\le N,\nu\in \Z} 
\frac{1}{2} |a_{i\nu}|^2 u_i$ and we have 
\[
L_X \cong \mu_{\infty}^{-1}(\eta)/\T.  
\]
Note that we use the same $\eta$ as we used to define $X$. 
The space $L_X$ here is infinite dimensional, but   
by cutting off Laurent series at a finite exponent, 
we can write $L_X$ as an inductive limit of smooth projective toric 
varieties (or compact symplectic toric manifolds). 
We endow $L_X$ with the inductive limit topology.

The algebraic model $L_X$ has several properties analogous to 
the universal cover $\widetilde{X^{S^1}}$ of free loop space. 
First, this has an $S^1$-action of the loop rotation: 
\[
[z_1(\zeta),\dots, z_N(\zeta)] \mapsto [z_1(e^{-\sqrt{-1}\theta}\zeta), \dots, 
z_N(e^{-\sqrt{-1}\theta}\zeta)], \quad e^{\sqrt{-1}\theta} \in S^1.  
\]
Second, this $S^1$ action is Hamiltonian with respect to 
the reduced symplectic form on $L_X$. The Hamiltonian $H$ is given by 
\[
H[z_1(\zeta),\dots,z_N(\zeta)] = \frac{1}{2} \sum_{i=1}^N\sum_{\nu\in\Z} 
\nu |a_{i\nu}|^2  \quad \text{on the level set } \mu_{\infty}^{-1}(\eta).   
\]
This is an analogue of the classical action functional ``$\int p dq$" 
on $\widetilde{X^{S^1}}$. 
Third, the group $H_2(X,\Z) \cong \pi_1(X^{S^1})$ of covering transformations  
on $\widetilde{X^{S^1}}$ acts on the model $L_X$. 
For $d\in H_2(X,\Z)$, we define a covering transformation $Q^d$ as 
\[
Q^d\colon [z_1(\zeta),\dots, z_N(\zeta)] \mapsto 
[\zeta^{-\pair{u_1}{d}} z_1(\zeta),\dots,\zeta^{-\pair{u_N}{d}} z_N(\zeta)]
\] 
Take a basis of $H_2(X,\Z)$ dual to $p_1,\dots, p_r$ and let 
$Q^1,\dots,Q^r$ be the corresponding covering transformations. 
Then we have 
$Q^d = (Q^1)^{\pair{p_1}{d}} \cdots (Q^r)^{\pair{p_r}{d}}$. 

We introduce $S^1$-equivariant cohomology classes $P_1,\dots, P_r\in H_{S^1}^2(L_X)$ 
corresponding to the basis $p_1,\dots,p_r$ of $H^2(X,\Z)$. 
By the natural isomorphism $H^2(X,\Z)\cong \Hom(\T_\C, \C^*)$, 
$p_a$ gives a representation of $\T_\C$ and 
gives a line bundle $\mL_{p_a}$ on $L_X$:
\[
\mL_{p_a} = \C\times_{p_a} L\mU_\eta \longrightarrow L_X. 
\]
This $\mL_{p_a}$ admits an $S^1$-action such that $S^1$ acts on 
Laurent polynomials by the loop rotation $\zeta\mapsto e^{-\sqrt{-1}\theta}\zeta$ 
and trivially on the $\C$ factor. We define 
\[
P_a := c_1^{S^1}(\mL_{p_a}) \in H_{S^1}^2(L_X).  
\]
It is easy to see that the pull-back by the covering transformation 
$Q^b$ and 
the multiplication by $P_a$ satisfies the commutation relation 
\[
[P_a, (Q^b)^*] = \hbar \delta_a^b (Q^b)^*
\]
as operators acting on $H^*_{S^1}(L_X)$. Here, $\hbar$ is a generator 
of the $S^1$-equivariant cohomology of a point. 
This commutation relation, discovered by Givental 
\cite{givental-ICM,givental-homologicalgeom}, 
yields a $D$-module structure on the equivariant Floer cohomology. 

\subsection{Equivariant Floer cohomology} 
The equivariant Floer theory here is considered to be a
Morse-Bott theory on $L_X$ with respect to the Hamiltonian function $H$. 
The critical set of $H$ is equal to the $S^1$-fixed point set in $L_X$  
and is given by the union of copies $X_d$ of $X$ over $d\in H_2(X,\Z)$, where  
\[
X_d := \{[a_1\zeta^{\pair{u_1}{d}}, \dots, a_N \zeta^{\pair{u_N}{d}}] 
\in L_X \; |\; a_i \in \C\}.  
\]
The gradient vector field of $H$ generates a flow $\phi_t$ on $L_X$ given by 
\[
\phi_t[z_1(\zeta),\dots, z_N(\zeta)] = 
[z_1(e^{-t}\zeta),\dots, z_N(e^{-t}\zeta)], \quad t\in \R. 
\]
Let $L_d^\infty$ be the closure (with respect to the inductive limit 
topology) of the stable manifold of $X_d$:  
\begin{align*}
L_d^{\infty}:&=\ov{\{z(\zeta)\in L_X |\;
\lim_{t\to \infty}\phi_t(z)\in X_d\}} \\ 
&= \left\{[z_1(\zeta),\dots,z_N(\zeta)]\in L_X \; \Big|\; z_i(\zeta) = 
\sum_{\nu\ge \pair{u_i}{d}} a_{i\nu} \zeta^\nu \right\}. 
\end{align*}
Similarly, let $L_{-\infty}^d$ be the 
closure of the unstable manifold of $X_d$: 
\[
L_{-\infty}^d := \left\{[z_1(\zeta),\dots,z_N(\zeta)]\in L_X\; 
\Big | \; z_i(\zeta) = \sum_{\nu\le \pair{u_i}{d}} a_{i\nu}\zeta^\nu
\right\}. 
\] 
We set $L_{d_1}^{d_2}=L_{d_1}^\infty \cap L_{-\infty}^{d_2}$. 
This is a finite dimensional smooth toric variety and 
considered to be a compactification of 
the union of gradient flowlines connecting $X_{d_1}$ and $X_{d_2}$. 
Note that the infinite dimensional spaces $L_X$ and $L_d^\infty$ here 
are quotients of open subsets of $\C^\infty$ by the torus $\T_\C$;   
the open subset here is the complement of union of 
infinite codimensional subspaces in $\C^\infty$, so in particular 
contractible. Therefore, $L_X$ and $L_d^\infty$ are homotopy equivalent to 
the classifying space $B\T_\C$ and we have 
\[
H^*(L_X) \cong H^*(L_d^\infty) \cong 
\C[c_1(\mL_{p_1}), \dots, c_1(\mL_{p_r})].  
\]
Repeating the same argument for the Borel construction 
$L_X\times_{S^1} ES^1$, $L_{d}^\infty\times_{S^1} ES^1$, we obtain  
\[
H_{S^1}^*(L_X) \cong H_{S^1}^*(L_d^{\infty})\cong \C[P_1,\dots,P_r,\hbar]. 
\]

First we explain the construction of equivariant Floer cohomology in 
the case of toric variety itself ($\mV=0$).  
Introduce a partial order ($\le$) on $H_2(X,\Z)$ as 
\[
d_1 \le d_2 \Longleftrightarrow 
L_{d_1}^{\infty}\subset L_{d_2}^\infty, 
\quad d_1, d_2 \in  H_2(X,\Z). 
\] 
Then $H_2(X,\Z)$ is a directed set. 
When $d_1 \le d_2$, there exists a push-forward 
$H_{S^1}^*(L_{d_1}^\infty)\rightarrow H_{S^1}^*(L_{d_2}^\infty)$ 
of $S^1$-equivariant cohomology.  
This push-forward is defined by 
the multiplication by the Euler class of 
the (finite rank) normal bundle $\mN$ of $L_{d_1}^\infty$
in $L_{d_2}^\infty$: 
\[
\cup\Euler_{S^1}(\mN) \colon
H_{S^1}^*(L_{d_1}^\infty) \cong \C[P_1,\dots,P_r,\hbar] \rightarrow 
\C[P_1,\dots,P_r,\hbar] \cong H_{S^1}^*(L_{d_2}^\infty). 
\]
where the Euler class is given by 
\[\Euler_{S^1}(\mN) = \prod_{i=1}^N 
\prod_{\nu=\pair{u_i}{d_2}}^{\pair{u_i}{d_1}-1}
\left(\sum_{a=1}^r m_i^a P_a -\nu \hbar\right).
\] 
By this push-forward, $\{H_{S^1}^*(L_d^\infty)\}_{d\in H_2(X,\Z)}$ 
forms an inductive system. 
Define the semi-infinite $S^1$-equivariant 
cohomology $H_{S^1}^{\infty/2}(L_X)$ as the inductive limit:  
\[
H^{\infty/2}_{S^1}(L_X):=\injlim_{d} H^*_{S^1}(L_d^{\infty}).  
\]
The covering transformation $Q^{d'}$ acts on the inductive system 
by pull-backs 
\[
(Q^{d'})^* \colon 
H^{*}_{S^1}(L_d^\infty) \cong 
H^{*}_{S^1}(L_{d+d'}^\infty) 
\]
so induces a map $(Q^{d'})^*\colon H^{\infty/2}(L_X)\rightarrow 
H^{\infty/2}(L_X)$. The equivariant class $P_a$ also acts on 
the semi-infinite cohomology by the cup product. 
These actions satisfy the commutation relation 
\[
[P_a,(Q^b)^*]=\hbar\delta_a^b(Q^b)^* \quad \text{on } 
H^{\infty/2}_{S^1}(L_X). 
\]
Hence $H^{\infty/2}_{S^1}(L_X)$ becomes a module over 
the polynomial Heisenberg algebra. 
The semi-infinite cohomology defined here contains 
classes of ``semi-infinite" cycles. 
For example, the class  
\[
\Delta := \text{the image of } 1\in H_{S^1}^*(L_0^\infty)
\text{ in }H_{S^1}^{\infty/2}(L_X)
\]
can be viewed as the Poincar\'{e} dual 
of the fundamental class of the space $L_0^\infty$. 
Define a filtration of $H_{S^1}^{\infty/2}(L_X)$ as 
\begin{align*}
F^n(H_{S^1}^{\infty/2}(L_X)) &:= 
\sum_{d\in \Box_n} H_{S^1}^*(L_d^\infty) 
= \sum_{d\in \Box_n} Q^d \C[P_1,\dots,P_r,\hbar]\Delta, \\
\text{where }
\Box_n & := \{d\in H_2(X,\Z)\;|\; 
\sum_{a=1}^r \pair{p_a}{d}\ge n, \ 
\pair{p_a}{d}\ge 0, \forall a\}.  
\end{align*}
\begin{definition}
The $S^1$-equivariant Floer cohomology $FH^*_{S^1}(L_X)$ 
is defined to be the completion of $F^0(H^{\infty/2}_{S^1}(L_X))$ 
with respect to the above filtration. 
Then $FH^*_{S^1}(L_X)$ becomes a module over  
$\C[\hbar][\![Q^1,\dots,Q^r]\!]\langle P_1,\dots,P_r \rangle$. 
\end{definition} 

Second, we explain the case where 
we have a vector bundle $\mV$ on $X$. 
We introduce another fiberwise $S^1$
action on $\mV$ and work $T^2$-equivariantly. 
We use $\hbar$ and $\lambda$ for generators of $T^2$-equivariant 
cohomology of a point, where 
$\hbar$ corresponds to the $S^1$-action rotating loops  
and $\lambda$ corresponds to the additional fiberwise $S^1$-action. 
Recall that $\mV$ is a sum of line bundles 
$\mV_1\oplus\cdots \oplus\mV_l$ with the first Chern class 
$v_i=c_1(\mV_i)$. 
Let $\mV_{i,\nu}$ be a line bundle on 
$L_X$ given by $v_i\in H^2(X,\Z)\cong \Hom(\T_\C,\C^*)$: 
\[
\mV_{i,\nu} := \C \times_{v_i} L\mU_\eta 
\longrightarrow L_X
\] 
endowed with the $S^1$-action (corresponding to the loop rotation): 
\[
[v, (z_1(\zeta),\dots,z_N(\zeta))] \mapsto 
[e^{-\sqrt{-1}\nu\theta}v, 
(z_1(e^{-\sqrt{-1}\theta}\zeta),\dots, z_N(e^{-\sqrt{-1}\theta}\zeta))], 
\quad e^{\sqrt{-1}\theta} \in S^1.  
\]
The additional $S^1$ action acts on $\mV_{i,\nu}$ 
by scalar multiplication on each fiber and trivially on the base. 
For a function $f\colon \{1,\dots,l\} \rightarrow \Z$, 
let $\mV_f^\infty$ be the infinite dimensional $T^2$-equivariant 
vector bundle on $L_X$: 
\[
\mV_f^\infty := \bigoplus_{i=1}^l \bigoplus_{\nu\ge f(i)} \mV_{i,\nu}.
\] 
For $d\in H_2(X,\Z)$, we define $\mV_d^\infty$ by regarding 
$d$ as a function $i\mapsto \pair{v_i}{d}$ on $\{1,\dots,l\}$. 
Define a partial order on $H_2(X,\Z)\times \Z^l$ as 
\[
(d_1,f_1) \le (d_2,f_2) \Longleftrightarrow 
L_{d_1}^\infty \subset L_{d_2}^\infty \text{ and } 
f_1(i)\le f_2(i) \ \forall i\in \{1,\dots, l\}.
\]
We define for $(d,f)\in H_2(X,\Z)\times \Z^l$, 
\[
H_{T^2}^*(L_d^\infty/\mV_f^\infty) := H_{T^2}^*(L_d^\infty) 
\cong \C[P_1,\dots,P_r\lambda,\hbar].  
\]
We could imagine $H_{T^2}^*(L_d^\infty/\mV_f^\infty)$ 
as the cohomology of the zero-locus $s^{-1}(0)\subset L_d^\infty$ 
if we have a transverse section $s$ of 
$\mV_f^\infty$ over $L_d^\infty$. In this case, an element of 
$H_{T^2}^*(L_d^\infty/\mV_f^\infty)$ could be viewed as 
the restriction to $s^{-1}(0)$ of a 
cohomology class of the ambient $L_d^\infty$.  
When $(d_1,f_1)\le (d_2,f_2)$, we define the push-forward 
$H_{T^2}^*(L_{d_1}^\infty/\mV_{f_1}^\infty)\rightarrow 
H_{T^2}^*(L_{d_2}^\infty/\mV_{f_2}^\infty)$ by 
multiplying the Euler class of the normal bundle $\mN$ 
of $L_{d_1}^\infty$ in $L_{d_2}^\infty$ 
and the quotient $\mV_{f_1}^\infty/\mV_{f_2}^\infty$: 
\[
\cup \Euler_{S^1}(\mN) \Euler_{T^2}(\mV_{f_1}^\infty/\mV_{f_2}^\infty) 
\colon H_{T^2}^*(L_{d_1}^\infty/\mV_{f_1}^\infty)
\longrightarrow  H_{T^2}^*(L_{d_2}^\infty/\mV_{f_2}^\infty). 
\] 
Note that we have a natural inclusion 
$\mV_{f_2}^\infty \subset \mV_{f_1}^\infty$ with finite rank quotient.  
The $S^1$-equivariant characteristic class of the quotient is given by 
\[
\Euler_{T^2}(\mV_{f_2}^\infty/\mV_{f_1}^\infty) 
= \prod_{i=1}^l \prod_{\nu=f_1(i)}^{f_2(i)-1}
\left(\sum_{a=1}^r l_i^a P_a -\nu\hbar + \lambda\right). 
\]
As before, $\{H_{T^2}(L_d^\infty/\mV_f^\infty)\}_{(d,f)}$ forms an 
inductive system. 
Let $\hH_{T^2}^{\infty/2}(L_{X/\mV})$ denote its inductive limit. 
This time, we define $H^{\infty/2}_{T^2}(L_{X/\mV})$ as a {\it submodule} 
generated by the images of $H_{T^2}^*(L_d^\infty/\mV_d^\infty)$ 
in the inductive limit $\hH_{T^2}^{\infty/2}(L_{X/\mV})$: 
\[
H_{T^2}^{\infty/2}(L_{X/\mV}) := \sum_{d\in H_2(X,\Z)} 
H_{T^2}^*(L_d^\infty/\mV_d^\infty) 
\subset \hH_{T^2}^{\infty/2}(L_{X/\mV}).  
\]
This semi-infinite cohomology is considered to contain  
semi-infinite cycles like 
``$[L_d^\infty]\cap \Euler_{T^2}(\mV_d^\infty)$''. 
The covering transformation $Q^{d'}$ acts on the inductive system 
by the pull-backs: 
\[
(Q^{d'})^* \colon H_{T^2}^*(L_d^\infty/\mV_f^\infty) 
\cong H_{T^2}^*(L_{d+d'}^\infty/\mV_{f+d'}^\infty),  
\]
so acts on the inductive limit. 
This action together with the multiplication by $P_a$ 
preserves the submodule $H_{T^2}^{\infty/2}(L_{X/\mV})$ 
and makes it a module over the polynomial Heisenberg algebra. 
Let $\Delta$ be the image of $1\in H_{T^2}^0(L_0^\infty/\mV_0^\infty)$ 
in $H_{T^2}^{\infty/2}(L_{X/\mV})$. 
As in Section \ref{sec:QDM}, we set 
\[
K := \begin{cases}
\C[\lambda] & \text{if $\mV$ is convex 
\emph{i.e.} $v_i$ is nef for all $i$}, \\
\C(\!(\lambda^{-1})\!) & \text{otherwise}.
\end{cases}
\]
We will consider the semi-infinite cohomology over $K\{\hbar\}$ 
(see Section \ref{subsec:QDM} for this notation). 
Define a filtration on $H_{T^2}^{\infty/2}(L_{X/\mV})$ by:
\begin{equation}
\label{eq:filtr_semiinf}
F^n(H_{T^2}^{\infty/2}(L_{X/\mV})) := 
\sum_{d\in \Box_n} H_{T^2}^*(L_d^\infty/\mV_d^\infty)
= \sum_{d\in \Box_n}
Q^d \C[\lambda,\hbar,P_1,\dots,P_r] \Delta. 
\end{equation}
\begin{definition} 
The $T^2$-equivariant Floer cohomology $FH_{T^2}^*(L_{X/\mV})$ 
for a pair $(X,\mV)$ is defined to be the completion of 
$F^0(H_{T^2}^{\infty/2}(L_{X/\mV}))
\otimes_{H_{T^2}^*({\rm pt})}K\{\hbar\}$ 
with respect to the above filtration. 
This is a module over 
$\mD\cong K\{\hbar\}[\![Q^1,\dots,Q^r]\!]\langle P_1,\dots, P_r\rangle$.  
\end{definition} 

\begin{remark*}
In \cite{givental-homologicalgeom}, 
in case of hypersurfaces in the projective space, 
the semi-infinite class $\Delta$ was given in the form of 
the infinite product 
\[
\Delta = \prod_{i=1}^N \prod_{\nu<0} 
\left(\sum_{i=1}^r m_i^aP_a-\nu\hbar \right) 
\cup 
\prod_{i=1}^l \prod_{\nu\ge 0} 
\left(\sum_{i=1}^r l_i^a P_a -\nu\hbar + \lambda\right). 
\]
This is here interpreted as an element 
in the inductive limit. 
\end{remark*}

\subsection{Localization map, pairing and the freeness of Floer cohomology}
A solution to the $D$-module $FH_{T^2}^*(L_{X/\mV})$, 
\emph{i.e.} a $D$-module homomorphism from $FH_{T^2}^*(L_{X/\mV})$
to the trivial $D$-module $K\{\hbar,\hbar^{-1}\}\!\}[\![Q^1,\dots,Q^r]\!]$, 
is given by the localization of cohomology 
on the $S^1$-fixed set of the algebraic model $L_X$. 
For $\alpha\in \hH_{T^2}^{\infty/2}(L_{X/\mV})$ and $d\in H_2(X,\Z)$, 
we define a localization at $X_d\cong X$ by 
\[
\Loc_d(\alpha) := 
\frac{\alpha_{(d',f)}}
{\Euler_{S^1}(\mN)\Euler_{T^2}(\mV_d^\infty/\mV_f^\infty)}\Bigg|_{X_d}
\in H^*(X)\otimes \C[\hbar,\hbar^{-1}](\!(\lambda^{-1})\!), 
\]
where $\alpha_{(d',f)}\in H_{T^2}^*(L_{d'}^\infty/\mV_f^\infty)$ 
is a representative of $\alpha$ such that $(d',f) \ge (d,d)$ 
and $\mN$ is the normal bundle of $L_d^\infty$ in $L_{d'}^\infty$. 
Here, we expand the inverse Euler classes 
around $\lambda=\infty$ so that $\Loc_d(\alpha)$ is defined over 
$\C[\hbar,\hbar^{-1}](\!(\lambda^{-1})\!)$. 
It was shown in \cite[Lemma 4.4]{iritani-EFC} that 
on the submodule $H_{T^2}^{\infty/2}(L_{X/\mV})$, 
$\Loc_d$ takes values in $H^*(X)\otimes \C[\hbar,\hbar^{-1}][\lambda]$ 
if $\mV$ is convex. 
Thus, in any case, $\Loc_d(\alpha)\in K\{\hbar,\hbar^{-1}\}\!\}$ 
for $\alpha\in H_{T^2}^{\infty/2}(L_{X/\mV})$. 
Moreover, we can show that the support $\{d\in H_2(X,\Z)\;|\;
\Loc_d(\alpha)\neq 0\}$ is contained in a cone of the form 
$d_0+\Lambda$, where $\Lambda$ is the Mori cone of $X$ and 
$d_0\in H_2(X,\Z)$ depends on $\alpha$.  
Therefore, we can define the map $\Loc$ as 
\begin{align*}
\Loc\colon  H_{T^2}^{\infty/2}(L_{X/\mV}) 
& \rightarrow H^*(X)\otimes 
K\{\hbar,\hbar^{-1}\}\!\}[\![Q^1,\dots,Q^r]\!]
[(Q^1)^{-1},\dots,(Q^r)^{-1}], \\ 
\alpha & \mapsto \sum_{d\in H_2(X,\Z)} \Loc_d(\alpha) Q^d.  
\end{align*}
This was denoted by $\Xi$ in \cite{iritani-EFC}. 
The following proposition was shown for convex $\mV$ in 
\cite{iritani-EFC}, but the proof there applies to a general $\mV$. 
\begin{proposition}[{\cite[Proposition 4.5,4.6]{iritani-EFC}}]
\label{prop:LocSol}  
The map $\Loc$ is a homomorphism of 
$H_{T^2}^*({\rm pt})\otimes\C[(Q^1)^{\pm},\dots,(Q^r)^\pm]$-modules 
and injective. Moreover, it satisfies the differential equation:
\[
\Loc(P_a \alpha) = (\hbar Q^a\parfrac{}{Q^a} + p_a\cup) \Loc(\alpha). 
\]
Here, $Q^a$ acts on the $H_{T^2}^{\infty/2}(L_{X/\mV})$ 
by the pull-back $(Q^a)^*$. 
It follows from the injectivity of $\Loc$ 
that the filtration defined in (\ref{eq:filtr_semiinf}) is Hausdorff. 
This $\Loc$ induces an embedding of 
$K\{\hbar\}[\![Q^1,\dots,Q^r]\!]$-modules:
\[
\Loc\colon 
FH_{T^2}^{*}(L_{X/\mV}) \rightarrow H^*(X)\otimes 
K\{\hbar,\hbar^{-1}\}\!\}[\![Q^1,\dots,Q^r]\!].  
\]
\end{proposition} 

By the localization map, we can define the $I$-function 
\cite{givental-mirrorthm-toric}  
for $(X,\mV)$ as $I_{X,\mV}(Q,\hbar):=\Loc(\Delta)$. 
The $I$-function is given explicitly as 
\begin{equation}
\label{eq:I-funct}
I_{X,\mV}(Q,\hbar) := \sum_{d\in \Lambda} 
\prod_{i=1}^N 
\frac{\prod_{\nu=-\infty}^{0}(u_i+\nu \hbar)}
     {\prod_{\nu=-\infty}^{\pair{u_i}{d}}(u_i+\nu\hbar)}
\prod_{j=1}^l 
\frac{\prod_{\nu=-\infty}^{\pair{v_j}{d}}(v_j+\nu\hbar+\lambda)}
     {\prod_{\nu=-\infty}^{0}(v_j+\nu\hbar+\lambda)}
Q^d. 
\end{equation} 
In Section \ref{sec:main}, we will show 
that the $I$-function determines the big QDM of $(X,\mV)$ 
by a (generalized) mirror transformation. 

We can also introduce a pairing on our Floer cohomology. 
We have an anti-$S^1$-equivariant automorphism  
$\barop\colon L_X\to L_X$ defined by 
$z_i(\zeta) \mapsto z_i(\zeta^{-1})$.  
This induces a map $\barop$:  
\begin{align*}
\barop\colon H_{T^2}^*(L_d^\infty) & \longrightarrow 
H_{T^2}^*(L_{-\infty}^{-d}), \\
f(P_1,\dots,P_r,\lambda,\hbar)& \longmapsto f(P_1,\dots,P_r,\lambda,-\hbar).
\end{align*} 
For $\alpha$, $\beta$ in $\hH_{T^2}^{\infty/2}(L_{X/\mV})$, 
define $\int_{L_{X/\mV}}\ov{\alpha}\cup \beta
\in \C[\hbar](\!(\lambda^{-1})\!) $ by 
\[
\int_{L_{X/\mV}}\ov{\alpha} \cup \beta 
= \int_{L_{d}^{d'}} 
\ov{\alpha_{(-d',-f')}}\beta_{(d,f)} \prod_{i=1}^l 
\frac{\prod_{\nu=-\infty}^{f'(i)} \Euler_{T^2}(\mV_{i,\nu}) 
      \prod_{\nu=f(i)}^{\infty} \Euler_{T^2}(\mV_{i,\nu})}
     {\prod_{\nu=-\infty}^\infty \Euler_{T^2}(\mV_{i,\nu})}. 
\]
Here, $\alpha_{(-d',-f')}$ and $\beta_{(d,f)}$ are representatives 
of $\alpha$ and $\beta$ such that $L_{d}^\infty \subset L_{d'}^\infty$. 
The Euler class in the denominator is 
expanded around $\lambda=\infty$. 
It is easy to check that this is independent of 
the choice of representatives. 
When $\mV$ is convex, 
the above pairing restricted to the submodule 
$H_{T^2}^{\infty/2}(L_{X/\mV})$ is the same as 
what is given in \cite{iritani-EFC}, so takes values in 
$\C[\lambda,\hbar]$. 
For $\alpha,\beta$ in $H_{T^2}^{\infty/2}(L_{X/\mV})$, 
we define another pairing $(\ov{\alpha},\beta)$ as
\[
(\ov{\alpha},\beta) := \sum_{d\in H_2(X,\Z)} Q^{-d} 
\int_{L_{X/\mV}} \ov{\alpha} \cup (Q^d)^*(\beta). 
\]
\begin{proposition}[{\cite[Proposition 4.8]{iritani-EFC}}]
The pairing $(\ov{\cdot},\cdot)$ on $H_{T^2}^{\infty/2}(L_{X/\mV})$ 
takes values in 
$K\{\hbar\}[\![Q^1,\dots,Q^r]\!][(Q^1)^{-1},\dots,(Q^r)^{-1}]$. 
We have 
\[
(\ov{\alpha},\beta) = \int_X \ov{\Loc(\alpha)} \cup \Loc(\beta) 
\cup \Euler_{S^1}(\mV),   
\]
where $\ov{\Loc(\alpha)}=\Loc(\alpha)|_{\hbar\mapsto -\hbar}$ is 
defined by flipping the sign of $\hbar$.  
The pairing $(\ov{\cdot},\cdot)$ can be extended on 
$FH_{T^2}^*(L_{X/\mV})$ via the right-hand side 
and takes values in $K\{\hbar\}[\![Q^1,\dots,Q^r]\!]$ 
on $FH_{T^2}^*(L_{X/\mV})$.  
\end{proposition}
\begin{remark*}
In \cite{iritani-EFC}, the author introduced 
the Floer homology $FH_*$ and defined a pairing 
between $FH_*$ and $FH^*$. 
The pairing above is equivalent to 
that in \cite{iritani-EFC} when we identify Floer homology 
with cohomology via the Poincar\'{e} duality map 
$\barop\colon FH^* \cong FH_*$ given in \cite{iritani-EFC}. 
\end{remark*} 

\begin{theorem}[{\cite[Theorem 4.10]{iritani-EFC}}]
\label{thm:freenessofefc}
The $T^2$-equivariant Floer cohomology $FH_{T^2}^*(L_{X/\mV})$
is a free module of rank $\dim H^*(X)$ 
over $K\{\hbar\}[\![Q^1,\dots,Q^r]\!]$. 
Moreover, the localization $\Loc_0$ at $X_0\cong X$ gives 
a canonical isomorphism:
\[
\Loc_0 \colon FH_{T^2}^*(L_{X/\mV}) \Big /
\sum_{a=1}^r Q^a FH_{T^2}^*(L_{X/\mV})\cong
H^*(X) \otimes K\{\hbar\}.  
\]
satisfying $\Loc_0(P_a\alpha)=p_a \Loc_0(\alpha)$, $\Loc_0(\Delta)=1$. 
In particular, for a set of polynomials 
$T_0,\dots, T_s\in \C[x_1,\dots,x_r]$, 
if $\{T_i(p_1,\dots,p_r)\}_{i=0}^s$ forms a basis of $H^*(X)$, 
then $\{T_i(P_1,\dots,P_r)\Delta\}_{i=0}^s$ gives a free basis 
of $FH_{T^2}^*(L_{X/\mV})$ over $K\{\hbar\}[\![Q^1,\dots,Q^r]\!]$. 
\end{theorem} 

We can prove this theorem by 
using the localization map and the pairing defined above. 
The proof for convex $\mV$ in \cite{iritani-EFC} 
again applies to this general case. 
In Section \ref{subsec:def_AQDM}, we will see that   
$FH_{T^2}^{*}(L_{X/\mV})$ is the small AQDM from this theorem. 

\section{Reconstruction of Big AQDM}
\label{sec:AQDM} 

In this section, we formulate the big and small 
{\it abstract} quantum $D$-modules (AQDM) and prove 
a reconstruction theorem from the small AQDM to the big one. 
We start from their definitions and explain canonical frames, 
affine structures on the base and reconstruction theorem. 


\subsection{Definitions}
\label{subsec:def_AQDM}

Let $K$ be an integral domain with a valuation 
$v_K\colon K\rightarrow \R\cup\{\infty\}$. 
We assume $K$ is complete with respect to the topology defined by $v_K$. 
See Section \ref{subsec:QDM} for notation used here.  
We assume that $K$ contains the field $\Q$ of rational numbers 
and $v_K|_\Q=0$.  
Basic examples of $K$ include $\Q$, $\C$, $\C[\lambda]$, 
$\C(\!(\lambda^{-1})\!)$, $\C[\![\bs]\!]$, $H_T^*({\rm pt})$ and 
$H_T^*({\rm pt})(\!(\lambda^{-1})\!)$, 
where the last two rings appear when we consider 
the $T$-equivariant quantum cohomology and its twist. 
We set the ideal 
$\frm_K=\{x\in K\;|\; v_K(x)>0\}$ if $v_K(K)\subset \R_{\ge 0}$ 
and $\frm_K=\{0\}$ otherwise. 
Define
\[
\mO:= K[\![Q^1,\dots,Q^r]\!], \quad 
\mO^\hbar:= K\{\hbar\}[\![Q^1,\dots,Q^r]\!].  
\]
Let $B$ be the formal neighborhood of zero in the affine space 
$\A_\mO^{s+1}$ over $\mO$
endowed with the structure sheaf $\mO_B$ and the sheaf 
$\mO_B^\hbar$ of algebra: 
\[
\mO_B := \mO[\![t^0,\dots,t^s]\!], \quad 
\mO_B^\hbar := \mO^\hbar[\![t^0,\dots,t^s]\!]. 
\]
The rings $\mO,\mO_B$ have a natural topology. 
In general, the topology on the 
formal power series ring $K[\![x_1,\dots,x_k]\!]$ over $K$
is given by the following fundamental neighborhood 
system of zero: 
\[
U_n := \left\{\sum_{i_1,\dots,i_k\ge 0} a_{i_1,\dots,i_k} x_1^{i_1} 
\cdots x_k^{i_k}\;\Big|\; v_K(a_{i_1,\dots,i_k})\ge n \text{ if } 
i_1+\cdots+i_k\le n \right\}. 
\]
A set of elements $\hatt^0,\dots,\hatt^s$ of $\mO_B$ 
is called a {\it co-ordinate system} on $B$ if  
$\hatt^i|_{Q=t=0}$ belongs to $\frm_K$ 
and if $\mO_B$ is topologically generated by $\hatt^i$ as an $\mO$-algebra, 
\emph{i.e.} $\mO_B=\mO[\![\hatt^0,\dots,\hatt^s]\!]$. 
By an {\it $\mO$-valued point} on $B$, we mean 
a continuous $\mO$-algebra homomorphism 
$\phi\colon \mO_B\rightarrow \mO$. 
For example, any co-ordinate system $\hatt^0,\dots,\hatt^s$ on $B$ 
defines an $\mO$-valued point $\hatt^0=\cdots=\hatt^s=0$. 
Conversely, any $\mO$-valued point $\bb$ on $B$ is given by 
$\hatt^0=\cdots=\hatt^s=0$ for some co-ordinate system. 
In this case, $\hatt^0,\dots,\hatt^s$ are called 
{\it co-ordinates centered at $\bb$}. 
Let $\mD_B$ denote the Heisenberg algebra defined by  
\[
\mD_B := \mO_B^\hbar \langle \Qp_1,\dots,\Qp_r, 
\bwp_0,\dots, \bwp_s \rangle/\mI_B.
\]
Here, $\mI_B$ is the two-sided ideal generated by 
\begin{gather*}
[\Qp_a,f(Q,t,\hbar)]-\hbar Q^a \parfrac{}{Q^a}f(Q,t,\hbar), \quad 
[\bwp_i, f(Q,t,\hbar)] -\hbar\parfrac{}{t^i}f(Q,t,\hbar), \\ 
[\Qp_a,\Qp_b], \quad 
[\bwp_i,\bwp_j], \quad 
[\Qp_a, \bwp_i],  
\end{gather*}
where $f(Q,t,\hbar)\in \mO^\hbar_B$. 
When $B$ is zero-dimensional over $\mO$, \emph{i.e.} $\mO_B=\mO$,  
$\mD_B$ coincides with the Heisenberg algebra $\mD$ 
in Section \ref{subsec:QDM}. 
The generators $\Qp_a, \bwp_i$ of $\mD_B$ 
depend on the choice of co-ordinates on $B$. 
Another co-ordinate system $\hatt^0,\dots,\hatt^s$ 
on $B$ yields the following generators 
$\hQp_a, \hbwp_i$ of $\mD_B$: 
\begin{equation}
\label{eq:genchangeDB}
\hQp_a = \Qp_a + \sum_{j=0}^s Q^a \parfrac{t^j(Q,\hatt)}{Q^a} 
\bwp_j, \quad 
\hbwp_i 
= \sum_{j=0}^s \parfrac{t^j(Q,\hatt)}{\hatt^i} \bwp_j.  
\end{equation} 
For a $\mD_B$-module $\mEh$ and an $\mO$-valued point $\bb$ on $B$, 
we set 
\[
\mEh_{\bb}:= \mEh\Big/\sum_{i=0}^s t^i \mEh, \quad  
\mEh_{\bb,0}:= \mEh_{\bb}\Big/\sum_{a=1}^r Q^a \mEh_{\bb}, \quad
\mE_{\bb,0}:= \mEh_{\bb,0}/\hbar\mEh_{\bb,0}, 
\]
where $t^0,\dots,t^s$ are co-ordinates centered at $\bb$. 
When $B$ is zero-dimensional over $\mO$, 
we will use the notation $\mEh_{0}$, $\mE_{0}$ instead of 
$\mEh_{\bb,0}, \mEh_{\bb,0}$.  
Let $p_a$ be the operator acting on $\mEh_{\bb,0}$ and $\mE_{\bb,0}$  
induced from $\Qp_a$ and 
$\wp_i$ be the operator acting on $\mE_{\bb,0}$ 
induced from $\bwp_i$.
Note that $p_a$ does not depend on the choice of 
co-ordinates $t^i$ on $B$ whereas $\wp_i$ does. 
Then $\mEh_{\bb}$ has the structure of a $\mD$-module, 
$\mEh_{\bb,0}$ becomes a $K\{\hbar\}[p_1,\dots,p_r]$-module 
and $\mE_{\bb,0}$ is a 
$K[p_1,\dots,p_r, \wp_0,\dots,\wp_s]$-module. 
\begin{definition}
\label{def:AQDM} 
An {\it abstract quantum $D$-module} or {\it AQDM} on the base $B$ 
is a $\mD_B$-module $\mEh$ satisfying the following conditions 
at an $\mO$-valued point $\bb$ on $B$: 
\begin{itemize}
\item[(1)] $\mEh$ is a finitely generated free $\mO_B^\hbar$-module. 
In particular, $\mEh_{\bb,0}$ (resp.\ $\mE_{\bb,0}$) 
is a free $K\{\hbar\}$-module 
(resp.\ free $K$-module) of finite rank. 
\item[(2)] 
There exists a splitting $\Phi_0\colon 
\mE_{\bb,0}\rightarrow \mEh_{\bb,0}$ 
such that the induced map 
$\Phi_0\colon \mE_{\bb,0}\otimes_K K\{\hbar\} \rightarrow \mEh_{\bb,0}$ 
is an isomorphism of $K\{\hbar\}[p_1,\dots,p_r]$-modules.
\item[(3)] 
Matrix elements of $p_a^n\in \End_K(\mE_{\bb,0})$
(with respect to a certain $K$-basis) 
have valuations bounded from below uniformly in $n\ge 0$. 
\end{itemize}
A {\it small AQDM} is an AQDM over the zero-dimensional base $B$, 
\emph{i.e.} $\mO_B=\mO$. 
A {\it big AQDM} is an AQDM satisfying the additional condition 
\begin{itemize}
\item[(4)] There exists an element $e_0\in \mE_{\bb,0}$ 
such that 
$\{\wp_i e_0\}_{i=0}^s$ is a $K$-basis of $\mE_{\bb,0}$. 
\end{itemize}
\end{definition} 
The condition (2) above   
is equivalent to that $p_a\in \End(\mEh_{\bb,0})$ is represented by an 
$\hbar$-independent matrix in a suitable 
$K\{\hbar\}$-basis of $\mEh_{\bb,0}$. 
The condition (3) is technically necessary; we will 
use this to construct a fundamental solution below.  
This is always satisfied when $v_K(K)\subset \R_{\ge 0}$. 
The condition (4) implies that the rank of $\mEh$ is 
same as the dimension $s+1$ of the base space $B$ over $\mO$. 
\begin{remark*}
In \cite{iritani-EFC}, we assumed that $\{e_0,p_1e_0,\dots,p_re_0\}$ 
is part of a $K$-basis of $\mE_{0}$ for some $e_0\in \mE_{0}$
for a small AQDM $\mEh$. 
We will see this condition 
in Proposition \ref{prop:divisor_AQDM} below. 
We will see in Theorem \ref{thm:canonicalframe} below that  
the conditions (2), (3), (4) do not depend on 
the choice of an $\mO$-valued point $\bb$, \emph{i.e.} 
if they hold at one $\bb$, then they hold at any other point. 
\end{remark*}
\begin{example} 
(1) The small and big QDM for a pair $(M,\mV)$ in Section \ref{sec:QDM} 
are the small and big AQDM respectively. 
The actions of $\Qp_a$ and $\bwp_i$ are given respectively by 
the Dubrovin connections $\smnabla^\hbar_a$ and $\nabla^\hbar_i$. 
We have a natural identification $\mE_{\bb,0}\cong H^*(M)\otimes K$ 
and the actions of $p_a$ and $\wp_i$ correspond to the cup products 
by $p_a$ and $p_i$ respectively. 
We can set $e_0$ to be the unit class $1\in H^0(M)$. 

(2) Theorem \ref{thm:freenessofefc} shows that 
the equivariant Floer cohomology $FH_{T^2}^*(L_{X/\mV})$ 
becomes a small AQDM, where the action of $\Qp_a$ and $Q^a$ 
is given by the cup product by $P_a\in H^2_{S^1}(L_X)$ 
and the pull back by $Q^a$ respectively.   
We have $\mE_{\bb,0}\cong H^*(X)\otimes K$ and 
the action of $p_a$ corresponds to the cup product by $p_a$.  
Note that we can take $\Phi_0$ as the inverse of $\Loc_0$.  
\end{example} 

A free basis (trivialization) of $\mEh$ over $\mO_B^\hbar$ 
is not a priori given.  
By the conditions (1) and (2), we can choose a splitting 
$\Phi\colon \mE_{\bb,0}\rightarrow \mEh$ inducing an isomorphism 
of $\mO^\hbar_B$-modules $
\Phi\colon \mE_{\bb,0}\otimes_K \mO_B^\hbar \rightarrow \mEh $ 
such that 
\begin{align}
\label{eq:frame_zero}
\begin{split}
&\text{\it the map } \  
\Phi_0\colon \mE_{\bb,0}\otimes_K K\{\hbar\} \rightarrow \mEh_{\bb,0}
\ \ \text{\it induced from $\Phi$}  \\ 
&\text{\it is an isomorphism of $K\{\hbar\}[p_1,\dots,p_r]$-modules.}
\end{split}
\end{align} 
We call such a splitting $\Phi$ (resp. $\Phi_0$)  
a {\it frame} of $\mEh$ (resp. $\mEh_{\bb,0}$). 
Two frames $\Phi$ and $\hPhi$ of $\mEh$ are related by a 
{\it gauge transformation} $G$ in 
$\Aut_{\mO^\hbar_B}(\mE_{\bb,0}\otimes_K \mO^\hbar_B)$
as $\Phi \mapsto \hPhi:=\Phi\circ G$. 
For a gauge transformation $G$, we need to assume that 
$G_{\bb,0}:=G|_{Q=t=0}\in 
\Aut_{K\{\hbar\}}(\mE_{\bb,0}\otimes_K K\{\hbar\})$ 
is a homomorphism of $K\{\hbar\}[p_1,\dots,p_r]$-modules 
and $G|_{Q=t=\hbar=0}=\id_{\mE_{\bb,0}}$ 
so that $\hPhi=\Phi\circ G$ gives a new frame. 
Using a frame $\Phi$, 
we can write down the $\mD_B$-module structure of $\mEh$ as follows. 
Define operators $\smnabla_a^\hbar$, $1\le a\le r$ and 
$\nabla_i^\hbar$, $0\le i\le s$ acting on 
$\mE_{\bb,0}\otimes_K \mO_B^\hbar$ by 
\[
\Phi(\smnabla_a^\hbar v) = \Qp_a \Phi(v), \quad 
\Phi(\nabla_i^\hbar v) = \bwp_i \Phi(v), \quad  
 v\in \mE_{\bb,0}\otimes_K 
\mO_B^\hbar. 
\]
These operators commute and take of the form 
\begin{align*}
\smnabla_a^\hbar = \hbar Q^a\parfrac{}{Q^a} + A_a(Q,t,\hbar), \quad 
\nabla_i^\hbar =\hbar t^i \parfrac{}{t^i} + \Omega_i(Q,t,\hbar) 
\end{align*} 
where $A_a, \Omega_i\in \End(\mE_{\bb,0})\otimes_K \mO^\hbar_B$.  
They corresponds to the Dubrovin connection of QDMs 
and are considered to be a flat connection on the trivial bundle 
$\mEh_{\bb,0}\times B\rightarrow B$. 
Note that by the condition (2) and our choice of $\Phi$, 
$p_a=A_a(0,0,\hbar)$ is independent of $\hbar$. 
Under a gauge transformation $G$, 
$A_a$ and $\Omega_i$ are transformed into  
\begin{equation}
\label{eq:gaugetransf_conn}
G^{-1} A_a G + \hbar G^{-1} Q^a \parfrac{G}{Q^a}, \quad 
G^{-1} \Omega_i G + \hbar G^{-1} \parfrac{G}{t^i}. 
\end{equation} 
Under the generator change (\ref{eq:genchangeDB}) of $\mD_B$ 
induced from a co-ordinate change on $B$, 
$A_a$ and $\Omega_i$ are transformed into
\begin{equation}
\label{eq:connchange_coordchange}
\hA_a = A_a + \sum_{j=0}^s Q^a \parfrac{t^j(Q,\hatt)}{Q^a} \Omega_j, \quad 
\Omega_{\hat{i}} =  \sum_{j=0}^s \parfrac{t^j(Q,\hatt)}{\hatt^i} \Omega_j.
\end{equation} 


\subsection{Canonical frames}
A frame $\Phi$ of an AQDM $\mEh$ is said to be {\it canonical} 
if the associated connection operators 
$A_a(Q,t,\hbar)$, $1\le a\le r$ and  $\Omega_i(Q,t,\hbar)$, 
$0\le i\le s$ are independent of $\hbar$. 
In case of QDMs, $A_a$ and $\Omega_i$ are given as 
the quantum multiplications by $p_a$ and  $p_i$, 
so they are by definition $\hbar$-independent.  
Therefore the QDM is a priori endowed with a canonical frame. 
In general, if the connection operators $A_a$ and $\Omega_i$ 
are $\hbar$-independent, we have  
\begin{gather}
\label{eq:commutingconn}
[A_a, A_b] = [A_a,\Omega_i] = [\Omega_i,\Omega_j]=0, \\
\label{eq:integrableconn}
Q^b\parfrac{A_a}{Q^b} = Q^b\parfrac{A_a}{Q^b}, \quad
\parfrac{A_a}{t^i} = Q^a\parfrac{\Omega_i}{Q^a}, \quad 
\parfrac{\Omega_i}{t^j} = \parfrac{\Omega_j}{t^i}.   
\end{gather} 
These follow from the flatness: 
$[\smnabla_a^\hbar,\smnabla^\hbar_b]=[\smnabla^\hbar_a,\nabla^\hbar_i
]=[\nabla^\hbar_i,\nabla^\hbar_j]=0$.  
In particular, $\mE_{\bb,0}\otimes \mO_B$ has the structure of 
an $\mO_B[A_a, \Omega_i]$-module; 
if moreover it is generated by $e_0$ as an $\mO_B[A_a,\Omega_i]$-module 
(this holds when $\mEh$ is a big AQDM),  
$\mE_{\bb,0}\otimes \mO_B$ has the structure of 
an $\mO_B[A_a,\Omega_i]$-algebra such that $e_0$ is the unit. 
This algebra is an analogue of the quantum cohomology. 

Let $t^0,\dots,t^s$ be a co-ordinate system centered at 
the base point $\bb$. Take a frame $\Phi$ of $\mEh$. 
A {\it fundamental solution} is an element $L=L(Q,t,\hbar)$ of 
$\End(\mE_{\bb,0})\otimes K\{\hbar,\hbar^{-1}\}\!\}[\![Q,t]\!]$ satisfying 
\begin{align}
\label{eq:fundsol_Omega}
\nabla^\hbar_i L &= 
\hbar \parfrac{}{t^i} L + \Omega_i L =0, \\
\label{eq:fundsol_A}
\smnabla^\hbar_a L -L p_a &= 
\hbar Q^a \parfrac{}{Q^a} L + A_a L - L p_a = 0, 
\end{align}
where $0\le i\le s, 1\le a\le r$. 
Here, $p_a=A_a|_{Q=t=0}\in \End(\mE_{\bb,0})$. 
These equations correspond to 
(\ref{eq:fundsol_big}), (\ref{eq:fundsol_small}) in case of QDM. 
We normalize a fundamental solution $L$ by the condition: 
\begin{equation}
\label{eq:normalizeL}
L(0,0,\hbar) = \id_{\mE_{\bb,0}}  
\end{equation} 

\begin{proposition}
\label{prop:Lexists}
There exists a unique fundamental solution $L$ 
normalized by (\ref{eq:normalizeL}). 
This $L$ depends on the base point $\bb$ and the frame $\Phi$.  
If $\Phi$ is a canonical frame, then 
$L=\id_{\mE_{\bb,0}} + O(\hbar^{-1})$. 
\end{proposition} 
\begin{proof}
The corresponding statement in the small case is proven in 
Proposition 3.5 in \cite{iritani-EFC}. 
For a function $f=f(Q^1,\dots, Q^r,t^0,\dots, t^s,\hbar)$, we write 
\[
f^{k,l} = f(Q^1,\dots,Q^k,0,\dots,0, t^0,\dots,t^{l-1},0,\dots, 0,\hbar).  
\]
Suppose by induction that we have a solution 
$L^{k,l}$ satisfying (\ref{eq:fundsol_Omega}) with $0\le i \le l-1$ and 
(\ref{eq:fundsol_A}) with $1\le a\le k$ and 
\begin{equation}
\label{eq:biggerthank}
A_a^{k,l} L^{k,l}  = L^{k,l} p_a    
\end{equation} 
for $k+1\le a\le r$. 
When $k=l=0$, these are satisfied for $L^{0,0}=\id_{\mE_{\bb,0}}$. 
We expand $L^{k+1,l}$ and $A_{k+1}^{k+1,l}$ as 
\begin{align*} 
L^{k+1,l} = L^{k,l} \sum_{n=0}^\infty 
T_n (Q^{k+1})^n, \quad 
A_{k+1}^{k+1,l} = \sum_{n=0}^\infty B_n (Q^{k+1})^n, 
\end{align*} 
where $T_n$ and $B_n$ are $\End(\mE_{\bb,0})$-valued functions 
in $Q^1,\dots,Q^k,t^0,\dots,t^{l-1},\hbar$ and $T_0=\id_{\mE_{\bb,0}}$ 
and $B_0=A_{k+1}^{k,l}$.  
From (\ref{eq:fundsol_A}) with $a=k+1$, we have 
\begin{align*}
n \hbar L^{k,l}T_n +\sum_{i=0}^n B_{n-i} L^{k,l}T_{i} - 
L^{k,l} T_n p_{k+1} = 0.   
\end{align*} 
Since $B_0 L^{k,l}=A_{k+1}^{k,l} L^{k,l} =L^{k,l} p_{k+1}$ 
by (\ref{eq:biggerthank}), we have 
\[
\left (1+\frac{1}{n\hbar} \ad(p_a)\right ) T_n = \frac{1}{n\hbar} 
\sum_{i=0}^{n-1} (L^{k,l})^{-1}B_{n-i} L^{k,l} T_{i}.  
\]
This determines $T_n$ in terms of $T_0,\dots,T_{n-1}$ because 
the operator $1+\ad(p_a)/(n\hbar)$ is invertible in 
$\End(\mE_{\bb,0})\otimes_K K\{\hbar,\hbar^{-1}\}\!\}$ by 
the condition (4): In fact, 
\[
\left(1+\frac{\ad(p_a)}{n\hbar}\right)^{-1} B 
= \sum_{m=0}^\infty \frac{1}{(n\hbar)^m}
\sum_{i=0}^m (-1)^i {m \choose i} p_a^i B p_a^{m-i}, \quad 
B\in \End(\mE_{\bb,0}). 
\] 
Note that the right-hand side is well-defined in 
$\End(\mE_{\bb,0})\otimes_K K\{\hbar,\hbar^{-1}\}\!\}$ 
by the condition (4). Thus we obtain $L^{k+1,l}$. 
Now we need to check that (\ref{eq:fundsol_A}), $1\le a\le k$   
holds for $L^{k+1,l}$. We have for $1\le a\le k$ 
\begin{align*}
\smnabla^\hbar_{k+1}(\smnabla^\hbar_a L^{k+1,l} - L^{k+1,l}p_a) 
&= \smnabla^\hbar_a \smnabla^\hbar_{k+1} L^{k+1,l} - L^{k+1,l} p_{k+1}p_a 
\\
&= (\smnabla^\hbar_a L^{k+1,l} -L^{k+1,l} p_a) p_{k+1}. 
\end{align*} 
and $(\smnabla^\hbar_a L^{k+1,l} - L^{k+1,l}p_a)|_{Q^{k+1}=0}=
\smnabla^\hbar_a L^{k,l}-L^{k,l}p_a=0$. 
Thus $T:=\smnabla_a^\hbar L^{k+1,l}-L^{k+1,l}p_a$ satisfies 
$\smnabla_{k+1}^\hbar T- Tp_{k+1}=0$ 
and $T|_{Q^{k+1}=0}=0$. This differential equation for $T$ 
turns out to have a unique solution, thus we have $T=0$.  
We also need to check 
\[
A_a^{k+1,l} L^{k+1,l} - L^{k+1,l} p_a =0 
\]
for $k+2\le a\le r$. But the left-hand side equals 
$\smnabla^\hbar_a L^{k+1,l}- L^{k+1,l} p_a$ restricted to 
$Q^{k+2}=\cdots=Q^r=t^l=\cdots=t^s=0$, 
thus this follows from the same argument. 
In the same way, we can also check that 
(\ref{eq:fundsol_Omega}), $0\le i\le l-1$ hold for $L^{k+1,l}$. 

Next we solve for $L^{k,l+1}$ 
under the same induction hypothesis. 
We expand $L^{k,l+1}$ and $\Omega_l^{k,l+1}$ as 
\[
L^{k,l+1} = \sum_{n=0}^\infty U_n (t^{l})^n, \quad 
\Omega_l^{k,l+1} = \sum_{n=0}^\infty C_n (t^l)^n,   
\]
where $U_n$ and $C_n$ are $\End(\mE_{\bb,0})$-valued functions 
in $Q^1,\dots,Q^k, t^0,\dots,t^{l-1}$ and 
$U_0=L^{k,l}$ and $C_0=\Omega_{l}^{k,l}$. 
By (\ref{eq:fundsol_Omega}) with $i=l$, we have 
\[
(n+1) \hbar U_{n+1} + \sum_{i=0}^{n} C_{n-i} U_i =0.   
\] 
Therefore, $U_{n+1}$ is recursively determined by $U_0,\dots,U_n$
and we obtain $L^{k,l+1}$. By a routine argument, we can check that 
(\ref{eq:fundsol_Omega}), $0\le i\le l$ and (\ref{eq:fundsol_A}), 
$1\le a\le k$ hold for $L^{k,l+1}$ and that 
$A_a^{k,l+1}L^{k,l+1}=L^{k,l+1}p_a$ holds for $k+1\le a \le r$.   
This completes the induction step. 
The last statement follows from the above procedure. 
\end{proof}

The fundamental solution $L$ above 
depends on a choice of frame. 
For another frame $\hPhi=\Phi\circ G$, 
the corresponding fundamental solution $\hL$ is given by 
\begin{equation}
\label{eq:gaugetransf_fundsol}
\hL=G^{-1} \circ L \circ G_{\bb,0}
\end{equation} 
where $G_{\bb,0}=G|_{Q=t=0}$. 
This easily follows from that $\hL|_{Q=t=0}=\id_{\mE_{\bb,0}}$ 
and that $G_{\bb,0}$ commutes with the action of $p_a$. 

\begin{definition}
Take a frame $\Phi$ and the fundamental solution $L$ associated to $\Phi$. 
Also choose an element $e_0\in \mE_{\bb,0}$. 
We assume $e_0$ satisfies the condition (4) in case of the big AQDM. 
The {\it $J$-function} of the AQDM $\mEh$ is defined to be 
$J(Q,t,\hbar):= L^{-1}e_0 \in \mE_{\bb,0}\otimes_K
K\{\hbar,\hbar^{-1}\}\!\}[\![Q,t]\!]$. 
This depends on the choice of $\bb$, $\Phi$ and $e_0$. 
\end{definition}


In case of QDM, this definition coincides with the original $J$-function 
(Definition \ref{def:J-funct}). 
In equivariant Floer theory, we choose a frame $\Phi$ of 
$FH_{T^2}^*(L_{X/\mV})$ such that $\Phi_0$ coincides with $\Loc_0^{-1}$.   
Then by Proposition \ref{prop:LocSol}, the inverse $L^{-1}$ of 
the corresponding fundamental solution is given by  
\[
\begin{CD}
L^{-1}\colon H^*(X) \otimes K\{\hbar\}[\![Q]\!] 
@>{\Phi}>> FH_{T^2}^*(L_{X/\mV}) 
@>{\Loc}>>  H^*(X)\otimes 
K\{\hbar,\hbar^{-1}\}\!\}[\![Q]\!].  
\end{CD}
\]
We take $e_0$ to be the image of $\Delta$. 
and choose $\Phi$ so that it satisfies $\Phi(e_0)=\Delta$. 
Then the $J$-function of $FH_{T^2}^*(L_{X/\mV})$ as an AQDM   
is given by the $I$-function (\ref{eq:I-funct}) because 
$I_{X,\mV}=\Loc(\Delta)=L^{-1}(e_0)$. 

Because the $J$-function is the inverse image of $e_0$ under $L$, 
we have a relation in $\mEh$ 
\[
f(Q,t,\Qp_a, \bwp_i, \hbar) \Phi(e_0) = 0, \quad f\in \mD_B,  
\] 
if and only if the $J$-function satisfies the differential equation 
\[
f(Q,t,\hbar Q^a\parfrac{}{Q^a}+p_a, \hbar\parfrac{}{t^i},\hbar) 
J(Q,t,\hbar) =0. 
\]
Because the big AQDM is generated by $e_0$ as a $\mD_B$-module 
by the condition (4), 
it is reconstructed by the $J$-function as a $\mD_B$-module. 
When we view $\mEh$ as a module over 
$\mD_B':=\mO_B^\hbar\langle \Qp_1,\dots,\Qp_r \rangle 
\subset \mD_B$, we have the following analogue of 
Proposition \ref{prop:Jgenerates}. 

\begin{proposition}
\label{prop:Jgenerates_AQDM}
Let $\mEh$ be an AQDM and $\Phi$ be a frame of $\mEh$. 
Let $J(Q,t,\hbar)$ be the $J$-function associated with $\Phi$ 
and $e_0\in \mE_{\bb,0}$. 
Assume that $\mE_{\bb,0}$ is generated by 
$e_0$ as a $K[p_1,\dots,p_r]$-module. 
Then $\mEh$ is generated by $\Phi(e_0)$ as a 
$\mD_B':=\mO_B^\hbar\langle \Qp_1,\dots,\Qp_r\rangle$-module. 
In other words, $\mEh$ is generated by 
$J(Q,t,\hbar)$ as a $\mD_B'$-module, i.e. 
we have an isomorphism of $\mD'_B$-module 
$\mD_B'/\frI\cong \mEh$ which sends $1$ to $\Phi(e_0)$, 
where $\frI$ is the left ideal
of $\mD_B'$ consisting of elements $f(Q,t,\Qp,\hbar)\in \mD_B'$ 
satisfying 
\[
f(Q, t, \hbar Q^a\parfrac{}{Q^a}+p_a, \hbar) J(Q,t,\hbar) =0  
\] 
In particular, the equivariant Floer cohomology in Section 
\ref{sec:EFT} is generated by 
the $I$-function (\ref{eq:I-funct}) as a $\mD$-module. 
\end{proposition}
\begin{proof}
This was stated in Theorem 3.17 in \cite{iritani-EFC} 
for small AQDMs. 
We omit the proof since it is elementary and 
is similar to that of Theorem 2.4 in \cite{iritani-EFC}. 
\end{proof} 

\begin{theorem}
\label{thm:canonicalframe}
For a given frame $\Phi_0$ of $\mEh_{\bb,0}$ satisfying 
(\ref{eq:frame_zero}), there exists  
a unique canonical frame $\Phi_{\rm can}$ which induces $\Phi_0$. 
The existence of a canonical frame implies that 
the conditions (2), (3) (and (4) in case of big AQDM) 
in Definition \ref{def:AQDM} hold at any $\mO$-valued point $\bb'$. 
\end{theorem} 
\begin{proof}
The proof here is similar to that in \cite[Theorem 3.9]{iritani-EFC}. 
We owe the idea of the Birkhoff factorization here 
to Guest \cite{guest}. This method also appeared in \cite{coates-givental}. 
We begin with an arbitrary frame $\Phi$ which induces $\Phi_0$. 
Let $L$ be the fundamental solution associated with $\Phi$. 
We claim that there exist 
$L_+ \in \End(\mE_{\bb,0})\otimes K\{\hbar\}[\![Q,t]\!]$ and 
$L_- \in \End(\mE_{\bb,0})\otimes K\{\!\{\hbar^{-1}\}\!\}[\![Q,t]\!]$ 
such that 
\begin{equation}
\label{eq:Birkhoff}
L= L_+ L_-, \quad L_-|_{\hbar=\infty} = \id_{\mE_{\bb,0}}, \quad 
L_+|_{Q=t=0} = \id_{\mE_{\bb,0}}. 
\end{equation} 
By expanding $L$, $L_+$ and $L_-$ in power series in $Q$ and $t$, 
we can solve for each coefficient of $L_+$ and $L_-$ recursively. 
We refer the reader to \cite[Theorem 3.9]{iritani-EFC} for details. 
In Equation (\ref{eq:calofLpos}) below, we will also 
give a formula of $L_+$ in terms of $L^{-1}$. 
By substituting $L$ with $L_+L_-$ in the differential equations 
(\ref{eq:fundsol_Omega}), (\ref{eq:fundsol_A}), we have 
\begin{align*}
L_+^{-1} \hbar \parfrac{L_+}{t^i} + L_+^{-1} \Omega_i L_+ 
&= - \hbar \parfrac{L_-}{t^i}L_-^{-1},   \\
L_+^{-1} \hbar Q^a \parfrac{L_+}{Q^a}  + L_+^{-1} A_a L_+ 
&= - \hbar Q^a \parfrac{L_-}{Q^a} L_-^{-1} + L_- p_a L_-^{-1}. 
\end{align*} 
The left-hand sides belong to 
$\End(\mE_{\bb,0})\otimes K\{\hbar\}[\![Q,t]\!]$ and 
the right-hand sides belong to 
$\End(\mE_{\bb,0})\otimes K\{\!\{\hbar^{-1}\}\!\}[\![Q,t]\!]$. 
Therefore, the both hand sides are $\hbar$-independent. 
Hence, the gauge transformation by $G:=L_+$ makes 
the connection operators $\hbar$-independent and   
$\Phi_{\rm can}=\Phi\circ G$ is a canonical frame. 
Note that $G$ does not change $\Phi_0$ 
since $G|_{Q=t=0} = \id_{\mE_{\bb,0}}$. 

Next we show the uniqueness of a canonical frame. 
Let $\Phi$ and $\Phi'$ be two canonical frames inducing the 
same frame $\Phi_0$ of $\mEh_{\bb,0}$. 
Then there exists a gauge transformation $G$ such that 
$\Phi' = \Phi \circ G$ and $G|_{Q=t=0} = \id_{\mE_{\bb,0}}$.  
Let $L$ and $L'$ be the fundamental solutions associated to $\Phi$ 
and $\Phi'$ respectively.  
Then by (\ref{eq:gaugetransf_fundsol}), we have $G= LL'^{-1}$.   
Because $\Phi$ and $\Phi'$ are canonical, 
we have $LL'^{-1} =\id_{\mE_{\bb,0}}+O(\hbar^{-1})$ by Proposition 
\ref{prop:Lexists}. On the other hand, since 
$G$ does not contain negative powers of $\hbar$, we have 
$G=LL'^{-1}|_{\hbar=\infty} = \id_{\mE_{\bb,0}}$, so $\Phi=\Phi'$. 

Finally, we see that the conditions (2), (3), (4) 
hold at any $\mO$-valued point $\bb'$. 
Let $A_a(Q,t)$ and $\Omega_i(Q,t)$ be 
the connection operators associated to a canonical frame. 
The frame $\Phi$ induces 
an identification $\mE_{\bb,0}\otimes K\{\hbar\}\cong \mEh_{\bb',0}$. 
Under this isomorphism, $p_a$ action on $\mEh_{\bb',0}$ corresponds 
to an $\hbar$-independent operator 
$A_a(0,t(\bb')|_{Q=0})\in \End(\mE_{\bb,0})\otimes K\{\hbar\}$. Thus 
the condition (2) holds at $\bb'$. 
The condition (4) implies that $\{\Omega_i(Q,t)e_0\}_{i=0}^s$ 
is an $\mO_B$-basis of $\mE_{\bb,0} \otimes \mO_B$. 
By specialization, we know (4) holds also at $\bb'$. 
The condition (3) is trivial when $v_K(K)\subset \R_{\ge 0}$, 
so we can assume $v_K(K)\not\subset \R_{\ge 0}$. 
In this case, $t(\bb')|_{Q=0}\in \frm_K=\{0\}$, so we have
$A_a(0,t(\bb')|_{Q=0}) = p_a$. 
From this the condition (3) at $\bb'$ is clear. 
\end{proof} 
%

\begin{remark*}
The positive part $L_+$ of the Birkhoff factorization, 
which serves as a gauge transformation to a canonical frame, 
may be calculated by the following formula: 
\begin{equation}
\label{eq:calofLpos}
L_+ v=\sum_{k=0}^\infty (\id - \pi_+\circ L^{-1})^k v 
\ \text{ for }v\in \mE_{\bb,0}. 
\end{equation} 
Here, $\pi_+$ is the map 
$\pi_+ \colon \mE_{\bb,0}\otimes_K K\{\hbar,\hbar^{-1}\}\!\}[\![Q,t]\!]
\rightarrow \mE_{\bb,0} \otimes_K K\{\hbar\}[\![Q,t]\!]$ 
induced from the projection 
$K\{\hbar,\hbar^{-1}\}\!\} = 
K\{\hbar\} \oplus \hbar^{-1} K\{\!\{\hbar^{-1}\}\!\} \to 
K\{\hbar\}$. 
To show this, first note that the right-hand side 
converges in $(Q,t)$-adic topology since 
$\pi_+\circ L^{-1} = \id + O(Q,t)$ on 
$\mE_{\bb,0}\otimes K\{\hbar\}$. 
From the Birkhoff factorization, we have $v=\pi_+ (L_-^{-1} v) 
= \pi_+ L^{-1} (L_+v)$ for $v\in \mE_{\bb,0}$. 
Therefore 
$\sum_{k=0}^\infty (\id -\pi_+ L^{-1})^k v =  
\sum_{k=0}^\infty (\id - \pi_+ L^{-1})^k \pi_+ L^{-1} (L_+v)
= \sum_{k=0}^\infty ((\id -\pi_+ L^{-1})^k - 
(\id -\pi_+ L^{-1})^{k+1})(L_+v) =L_+v$. 
Note also that the right-hand side of (\ref{eq:calofLpos}) 
is not $K\{\hbar\}$-linear but only $K$-linear in $v$. 
\end{remark*} 


\subsection{Affine structures} 
\label{subsec:affine}
We describe an affine structure on the base space of the big AQDM. 
The base space of the big QDM 
has a natural linear structure and the 
covariant derivative along a linear vector field
$\partial/\partial t^i$ gives   
\[
\nabla^\hbar_i 1 = p_i* 1 = p_i\cup 1. 
\]
where $\nabla^\hbar_i$ is the Dubrovin connection 
in Section \ref{sec:QDM}. 
This motivates the following definition:  
\begin{theorem-definition}
\label{thm-def:affine}
Let $\mEh$ be a big AQDM endowed with a canonical frame $\Phi$. 
Let $e_0$ be an element of $\mE_{\bb,0}$ satisfying the condition 
(4) in Definition \ref{def:AQDM}. 
There exists a co-ordinate system $\hatt^0,\dots,\hatt^s$ on $B$ such that 
\begin{equation}
\label{eq:affinecoord}
\Phi(\hwp_i e_0) = \hbwp_i \Phi(e_0), \quad 
\Phi(p_a e_0) = \hQp_a \Phi(e_0) 
\end{equation} 
for $0\le i\le s$ and $1\le a\le r$. Here, 
$\hQp_a$ and $\hbwp_i$ are generators of $\mD_B$ corresponding 
to the co-ordinates $\hatt^i$ and $\hwp_i\in \End(\mE_{\bb,0})$ 
is induced from $\hbwp_i$. 
The co-ordinate system having this property 
is unique up to affine transformations over $K$, 
thus defines an \emph{affine structure} on $B$ over $K$. 
This affine structure depends on the choice of $\Phi$ 
and $e_0$. 
We call such $\hatt^0,\dots,\hatt^s$ \emph{a flat co-ordinate system}
associated with $\Phi$ and $e_0$. 

(ii) The $J$-function associated with the 
canonical frame $\Phi$ and $e_0$ has the $\hbar^{-1}$-expansion 
of the form
\begin{equation}
\label{eq:asymptoticsofJ}
J = e_0 + \frac{1}{\hbar}
\sum_{i=0}^{s} \hatt^i (\hwp_i e_0) + O(\hbar^{-2}).  
\end{equation} 
\end{theorem-definition} 
\begin{proof}
Let $t^0,\dots,t^s$ be an arbitrary co-ordinate system 
centered at $\bb$. Let $\Qp_a$, $\bwp_i$ be the corresponding 
generators of $\mD_B$. 
Let $A_a(Q,t)$ and $\Omega_i(Q,t)$ be the connection operators 
associated with the canonical frame $\Phi$.  
By Proposition \ref{prop:Lexists} and $J=L^{-1} e_0$, 
the $J$-function associated with the canonical frame 
is of the form $J=e_0 + O(\hbar^{-1})$. Therefore, 
we can define an element $\hatt^i\in \mO_B$ by the expansion 
(note that $\{\wp_ie_0\}$ is a $K$-basis of $\mE_{\bb,0}$): 
\[
J = e_0 + \frac{1}{\hbar}
\sum_{i=0}^s \hatt^i (\wp_ie_0) + O(\hbar^{-2}). 
\]
First we check that $\hatt^0,\dots,\hatt^s$ 
form a co-ordinate system on $B$. We have 
\begin{align*}
\hbar\parfrac{}{t^j} J &= L^{-1} (\Omega_i e_0) = 
\Omega_j e_0 + O(\hbar^{-1}), \\
(\hbar Q^a\parfrac{}{Q^a} +p_a) 
J &= L^{-1} (A_a e_0) = A_ae_0 + O(\hbar^{-1}).  
\end{align*}
Hence, by comparing the leading terms in the $\hbar^{-1}$-expansions, 
\[
\sum_{i=0}^s \parfrac{\hatt^i}{t^j} (\wp_i e_0) = \Omega_j e_0 
= \wp_j e_0 + O(Q,t), \quad 
p_ae_0 + \sum_{i=0}^s Q^a \parfrac{\hatt^i}{Q^a} (\wp_i e_0) 
= A_a e_0. 
\]
From the first equation, 
the Jacobi matrix $(\partial \hatt^i/\partial t^j)=\delta^i_j+O(Q,t)$ 
turns out to be invertible. 
Therefore $\hatt^0,\dots,\hatt^s$ gives a co-ordinate system. 
On the other hand, it follows from the above formulas that 
\begin{equation}
\label{eq:affinestr_1} 
\wp_i e_0 = \sum_{j=0}^s \parfrac{t^j(Q,\hatt)}{\hatt^i}\Omega_j e_0, 
\quad 
p_a e_0 = A_a e_0 - \sum_{a=1}^r Q^a \parfrac{\hatt^i(Q,t)}{Q^a} 
(\wp_i e_0).  
\end{equation}
By the chain rule  
\[
0 = \sum_{k=0}^s \parfrac{\hatt^i(Q,t)}{t^k} 
Q^a\parfrac{t^k(Q,\hatt)}{Q^a} + Q^a\parfrac{\hatt^i(Q,t)}{Q^a} 
\]
and (\ref{eq:affinestr_1}), we find 
\[
p_a e_0 = A_a e_0 + \sum_{a=1}^r \sum_{k=0}^s 
\parfrac{\hatt^i(Q,t)}{t^k} Q^a \parfrac{t^k(Q,\hatt)}{Q^a} (\wp_i e_0)
= (A_a + \sum_{a=1}^r Q^a \parfrac{t^j(Q,\hatt)}{Q^a} \Omega_j) e_0. 
\]
This and the first equation of (\ref{eq:affinestr_1}) 
show that 
\[
\wp_i e_0 = \Omega_{\hat{i}} e_0, \quad p_a e_0 = \hat{A}_a e_0 
\]
where $\Omega_{\hat{i}}$ and $\hat{A}_a$ are connection operators 
corresponding to the new co-ordinates $\hatt^i$ 
(see (\ref{eq:connchange_coordchange})). 
Because $\wp_i = \hwp_i$ in this case, 
this is equivalent to (\ref{eq:affinecoord}). 
Thus we proved the existence of flat co-ordinates and (ii). 

Let ${t^0}',\dots, {t^s}'$ be another flat co-ordinate system. 
Let $\bwp_i',\Qp_a'\in \mD_B$ and $\wp_i'\in \End(\mE_{\bb,0})$ 
be the corresponding operators. 
Since two bases $\{\hwp_ie_0\}_{i=0}^s$ and $\{\wp_i'e_0\}_{i=0}^s$ 
of $\mE_{\bb,0}$ are related by an element in $GL(s+1,K)$,  
it follows from (\ref{eq:connchange_coordchange}) and 
(\ref{eq:affinecoord}) that the matrix 
$\partial {t^i}'(Q,\hatt)/\partial \hatt^j$ is in $GL(s+1,K)$. 
Similarly, we know that $Q^a (\partial {t^i}'(Q,\hatt)/\partial Q^a)=0$. 
Therefore, ${t^i}' = c^i+ \sum_{j=0}^s b^i_j \hatt^j$ 
for some $c^i, b^i_j\in K$. 
\end{proof} 

In case of the big QDM, 
linear co-ordinates on cohomology give a flat co-ordinate system 
in the above sense.  
We also have an analogue of the string and divisor equations 
in the context of big AQDM. 
 
\begin{proposition}
\label{prop:divisor_AQDM} 
Let $\mEh$ be a big AQDM. 
Let $\Phi$ be a canonical frame of $\mEh$ and 
$e_0$ be an element of $\mE_{\bb,0}$ 
satisfying the condition (4) in Definition \ref{def:AQDM}. 
Assume that $\{e_0, p_1e_0,\dots,p_re_0\}$ is part of a $K$-basis 
of $\mE_{\bb,0}$. There exists a flat co-ordinate system 
$t^0,\dots,t^s$ on $B$ associated with $\Phi$ and $e_0$ 
such that 
\begin{equation}
\label{eq:pequalswp}
\wp_0 e_0 = e_0, \quad \wp_a e_0 = p_a e_0, \quad 1\le a\le r. 
\end{equation}
Then the fundamental solution $L$ associated with $\Phi$ satisfies 
\[
(\parfrac{}{t^a} - Q^a \parfrac{}{Q^a}) L + 
L\circ \frac{p_a}{\hbar}= 0, \quad 
\hbar \parfrac{}{t^0} L = L 
\]
and the corresponding connection operators $A_a$ and $\Omega_i$ 
are independent of $t^0$ and satisfy 
\[
A_a = \Omega_a, \quad 
(\parfrac{}{t^a}- Q^a \parfrac{}{Q^a})\Omega_i=0, \quad \Omega_0=\id,  
\]
where $1\le a \le r$ and $0\le i\le s$. 
\end{proposition} 
\begin{proof}
By assumption, we can make flat co-ordinates satisfy (\ref{eq:pequalswp}) 
by a $K$-linear co-ordinate change. 
Then by the differential equations for $L$, we have 
\begin{align*}
\left(\hbar\parfrac{}{t^a} - (\hbar Q^a\parfrac{}{Q^a} + p_a)
\right) L^{-1} \wp_j e_0 & = 
L^{-1}(\nabla^\hbar_a -\smnabla^\hbar_a) \wp_j e_0 \\
&= L^{-1} (\Omega_a - A_a) \Omega_j e_0  \\
&= L^{-1} \Omega_j (\Omega_a -A_a) e_0 = 0
\end{align*} 
where we used $\wp_je_0=\Omega_j e_0$ in the second line 
and used (\ref{eq:commutingconn}) and 
$A_a e_0 =p_a e_0 =\wp_a e_0= \Omega_a e_0$ in the third line. 
Thus we have 
$(\hbar\parfrac{}{t^a} - (\hbar Q^a\parfrac{}{Q^a} + p_a))
L^{-1}=0$. 
Similarly, 
\[
\hbar \parfrac{}{t^0} L^{-1} \wp_j e_0 = L^{-1} \Omega_0 \wp_j e_0
= L^{-1}\Omega_0 \Omega_j e_0 = L^{-1} \Omega_j \wp_0 e_0 
= L^{-1} \Omega_j e_0 = L^{-1} \wp_j e_0. 
\] 
Thus $\hbar\parfrac{}{t^0} L^{-1} = L^{-1}$. 
The differential equations for $L$ follow from these. 
The statement on the connection operators easily follows from the 
differential equations for $L$. The details are left to the reader. 
\end{proof}


\begin{remark*}
This proposition shows that the connection $A_a$, $\Omega_i$ 
are functions in $Q^1 e^{t^1},\dots, Q^r e^{t^r}, t^{r+1},\dots, t^s$
in suitable flat co-ordinates. 
The variable $q^a$ used in the previous paper \cite{iritani-EFC} 
corresponds to the combination $q^a = Q^a e^{t^a}$.  
There, the small AQDM was a module over 
the Heisenberg algebra generated by 
$\boldsymbol{p}_a$, $q^b$ with commutation relation 
$[\boldsymbol{p}_a,q^b]=\delta_a^b \hbar q^b$. 
In \cite{iritani-EFC}, we allow a co-ordinate change in $q$-variables 
to solve for flat co-ordinates $\hq^a$. 
This method works for a small AQDM satisfying a certain ``nef" assumption 
(Assumption 3.15 in \cite{iritani-EFC}) but 
does not work beyond nef case. 
In this paper, introducing redundant co-ordinates $Q^a$ and $t^i$,  
we solve for flat co-ordinates $\hatt^i$ 
for the big AQDM (which will be reconstructed 
from the small one in the next section). 
Note that we keep a co-ordinate $Q^a$ fixed and 
allow a co-ordinate change only in $t$-variables. 
\end{remark*}

\subsection{Reconstruction}
For an AQDM $\mEh$ over $B$ and an $\mO$-valued point $\bb$ on $B$, 
a fiber $\mEh_{\bb}$ at $\bb$ becomes a small AQDM. 
Here we consider the reconstruction of a big AQDM from its fiber. 
The reconstruction theorem below is considered 
as a generalization of Kontsevich and Manin's reconstruction theorem 
(Theorem \ref{thm:KMreconstruction}). 

Let $\mEh$ be an AQDM over $B$. 
Let $t^0,\dots,t^s$ be a co-ordinate system on $B$. 
Let $B'$ be the formal neighborhood of zero in 
the affine space $\A_{\mO}^{l+1}$ over $\mO$ 
endowed with the structure sheaf 
$\mO_{B'}=\mO[\![\epsilon^0,\dots,\epsilon^l]\!]$ 
and the sheaf $\mO^\hbar_{B'} = \mO^\hbar
[\![\epsilon^0,\dots,\epsilon^l]\!]$ of algebras. 
By {\it a map} $f\colon B'\rightarrow B$, 
we mean a continuous $\mO$-algebra 
homomorphism $f^* \colon \mO_B \rightarrow \mO_{B'}$.  
More concretely, this is given by $s+1$ elements  
$t^0(Q,\epsilon),\dots, t^s(Q,\epsilon)$ in $\mO_{B'}$ 
such that $t^i(0,0) \in \frm_K$, where $t^i(Q,\epsilon)=f^*t^i$. 
We define a pulled back $\mD_{B'}$-module $f^*\mEh$ as 
\[
f^*\mEh = \mO^\hbar_{B'} \otimes_{\mO^\hbar_{B}} \mEh 
\]
where the actions of $\Qp'_a, \bwp'_i\in \mD_{B'}$ 
(corresponding to co-ordinates $\epsilon^i$) are given by 
\begin{align*}
\Qp_a' &= \hbar Q^a \parfrac{}{Q^a} \otimes 1 + 
1\otimes (\Qp_a+\sum_{j=0}^s Q^a\parfrac{t^j(Q,\epsilon)}{Q^a} \bwp_j), \\
\bwp_i' &= \hbar \parfrac{}{\epsilon^i} \otimes 1+ 
1\otimes \sum_{j=0}^s \parfrac{t^j(Q,\epsilon)}{\epsilon^i} \bwp_j.   
\end{align*} 
It is easy to check that $f^*\mEh$ becomes an AQDM over $B'$.

\begin{theorem}
\label{thm:reconstruction} 
Let $\mEh$ be a small AQDM. 
Assume that $\mE_{0}$ is generated by a single element $e_0$ 
as a $K[p_1,\dots,p_r]$-module. 
Then there exist a big AQDM $\hmEh$ over a base $B$, an $\mO$-valued 
point $\bb$ on $B$ and an 
isomorphism $\varphi\colon \mEh \cong \hmEh_{\bb}$ of $\mD$-modules  
satisfying the following universal property: 
For any AQDM $\mFh$ over $B'$ with an isomorphism 
$\varphi'\colon \mEh \cong \mFh_{\bb'}$ of $\mD$-modules 
at $\bb'$ on $B'$, 
there exist a unique map $f\colon B'\rightarrow B$ and 
a unique isomorphism $\theta\colon \mFh \cong f^*\hmEh$ 
of $\mD_{B'}$-modules such that 
$f(\bb')= \bb$ and $\theta_{\bb'} \circ \varphi' = \varphi$. 
Here, $\theta_{\bb'}\colon \mFh_{\bb'} \cong f^*\hmEh_{\bb}$ 
is an isomorphism of $\mD$-modules induced from $\theta$.  
In particular, $\hmEh$ is unique up to isomorphisms.
\end{theorem} 

\begin{proof} 
{\bf Construction of $\hmEh$:}     
Take a canonical frame $\Phi$ of $\mEh$. 
Let $A_a(Q)$, $1\le a\le r$ be the associated 
($\hbar$-independent) connection operators. 
Let $\{\be_0,\dots,\be_s\}$ be a basis of $\mE_{0}$ over $K$. 
Let $B$ be a space with co-ordinates $t^0,\dots,t^s$ 
and $\bb$ be the point defined by $t^0=\cdots=t^s=0$. 
We set $\mO_B=\mO[\![t^0,\dots,t^s]\!]$ and 
$\mO^\hbar_B=\mO^\hbar[\![t^0,\dots,t^s]\!]$ as usual. 
We claim that there exist unique flat connections 
$\smnabla^\hbar_a$, $1\le a\le r$ and 
$\nabla^\hbar_i$, $0\le i\le s$ acting on 
$\mE_0 \otimes_K \mO^\hbar_B$ of the form:
\[
\smnabla^\hbar_a = \hbar Q^a\parfrac{}{Q^a}+A_a(Q,t), 
\quad 
\nabla^\hbar_i = \hbar \parfrac{}{t^i} + \Omega_i(Q,t) 
\]
such that $A_a(Q,0)=A_a(Q)$ and $\Omega_i(Q,t) e_0 = \be_i$. 
Then $\hmEh :=\mE_0\otimes_K \mO^\hbar_B$ has the structure of 
a big AQDM endowed with a canonical frame  
such that $\hmEh_{\bb}=\mE_0\otimes\mO^\hbar$ 
is isomorphic to $\mEh$ via $\Phi$ (\emph{i.e.} $\varphi=\Phi^{-1}$). 
%
We expand $A_a(Q,t)$ and $\Omega_i(Q,t)$ as 
\[
A_a(Q,t) = \sum_{n\ge 0} A_a^{(n)}(Q,t), \quad 
\Omega_i(Q,t) = \sum_{n \ge 0} \Omega_i^{(n)}(Q,t),  
\]
where $A_a^{(n)}(Q,t)$ (resp. $\Omega_i^{(n)}(Q,t)$) is 
the degree $n$ part of $A_a(Q,t)$ (resp. $\Omega_i(q,t)$)  
with respect to the variables $t^0,\dots,t^s$. 
We also write $A_a^{\le n} :=\sum_{k=0}^n A_a^{(k)}$ and 
$\Omega_i^{\le n} = \sum_{k=0}^n \Omega_i^{(k)}$. 
Then we must have 
\begin{equation}
\label{eq:cond1_basis}
\Omega_i^{(n)}(Q,t)e_0 = \delta_{n,0} \be_i
\end{equation} 
The flatness of the connection implies 
(see (\ref{eq:commutingconn}), (\ref{eq:integrableconn}))
\begin{gather}
\label{eq:cond2_comm}
\sum_{k+l=n}[A_a^{(k)},A_b^{(l)}]=0, \quad 
\sum_{k+l=n}[A_a^{(k)},\Omega_j^{(l)}]=0, \quad 
\sum_{k+l=n}[A_a^{(k)},\Omega_j^{(l)}]=0, \\
\label{eq:cond3_integ}
\bpartial_b A_a^{(n)}=\bpartial_a A_b^{(n)},\quad 
\partial_j \Omega_i^{(n)}=\partial_i \Omega_j^{(n)}, \quad 
\partial_j A_a^{(n)}=\bpartial_a \Omega_j^{(n-1)}, 
\end{gather}
where $\bpartial_a = Q^a \partial/\partial Q^a$ and 
$\partial_i = \partial/\partial t^i$. 
Since $\mE_0$ is generated by $e_0$ as a $K[p_1,\dots,p_r]$-module 
and $p_a = A_a^{(0)}|_{Q=0}$, we know that $\mE_0\otimes_K \mO$ is 
generated by $e_0$ as an $\mO[A_1^{(0)},\dots,A_r^{(0)}]$-module. 
Thus, $\mE_0\otimes_K \mO$ 
has a unique $\mO[A_1^{(0)},\dots,A_r^{(0)}]$-algebra 
structure such that $e_0$ is a unit. 
Since $\Omega^{(0)}_i$ commutes with $A_a^{(0)}$ and 
$\Omega^{(0)}_i e_0 = \be_i$, 
$\Omega^{(0)}_i$ is identified with the multiplication 
by $\be_i$ in this algebra $\mE_0\otimes_K \mO$. 

Suppose by induction that 
we have $A_a^{(k)}$ for $0\le k \le m-1$ and $\Omega_i^{(k)}$ 
for $0\le k\le m-1$ satisfying all the conditions (\ref{eq:cond1_basis}), 
(\ref{eq:cond2_comm}), (\ref{eq:cond3_integ}) up to $n=m-1$. 
First, we can solve for $A_a^{(m)}$ uniquely 
using the third equation of (\ref{eq:cond3_integ}) for $n=m$.  
This is possible since the integrability condition is satisfied:   
$\partial_i(\bpartial_a\Omega_j^{(m-1)})=
\bpartial_a\partial_i\Omega_j^{(m-1)}=
\bpartial_a\partial_j\Omega_i^{(m-1)}=
\partial_j(\bpartial_a\Omega_i^{(m-1)})$. 
We need to check the first equation of (\ref{eq:cond3_integ}) 
with $n=m$. 
This follows from 
$\partial_i(\bpartial_b A_a^{(m)}-\bpartial_a A_b^{(m)})
=\bpartial_b\bpartial_a \Omega_i^{(m-1)}-
\bpartial_a\bpartial_b\Omega_i^{(m-1)}=0$. 
Also we need to check the first equation of (\ref{eq:cond2_comm})
with $n=m$. This follows from 
\begin{align*}
\partial_i &\sum_{k+l=m}[A_a^{(k)},A_b^{(l)}]
=\sum_{k+l=m}[\bpartial_a\Omega_{i}^{(k-1)},A_b^{(l)}]
  +\sum_{k+l=m}[A_{a}^{(k)},\bpartial_b\Omega_i^{(l-1)}] \\
&=\sum_{k+l=m-1}\left\{
  \bpartial_a[\Omega_i^{(k)},A_b^{(l)}]-
  [\Omega_i^{(k)},\bpartial_a A_b^{(l)}]
  +\bpartial_b[A_a^{(k)},\Omega_i^{(l)}]
  -[\bpartial_b A_a^{(k)},\Omega_i^{(l)}]
  \right\} \\
&=\sum_{k+l=m-1}[\Omega_i^{(k)},
  \bpartial_b A_a^{(l)}-\bpartial_a A_b^{(l)}]=0.
\end{align*}
Secondly, we solve for $\Omega_j^{(m)}$. 
Let $\frm\subset \mO_B$ be the ideal generated by $t^0,\dots,t^s$. 
The equation (\ref{eq:cond2_comm}) with $n\le m$ imply that 
$\mE_0\otimes (\mO_B/\frm^{m+1})$ has the structure 
of an $(\mO_B/\frm^{m+1})[A_a^{\le m}]$-module. 
Since this is again generated by $e_0$, this has a unique 
$(\mO_B/\frm^{m+1})[A_a^{\le m}]$-algebra structure 
such that $e_0$ is a unit. 
By (\ref{eq:cond1_basis}) and (\ref{eq:cond2_comm}), 
$\Omega_i^{\le m}$ commutes with $A_a^{\le m}$ on 
$\mE_0\otimes (\mO_B/\frm^{m+1})$ and $\Omega_i^{\le m}e_0=\be_i$,  
$\Omega_i^{\le m}$ is uniquely determined as the 
multiplication by $\be_i$ in this algebra. 
We need to check 
$\partial_j \Omega^{\le m}_i = \partial_i \Omega^{\le m}_j$. 
We have 
\begin{align*}
[\partial_j \Omega^{\le m}_i-\partial_j \Omega^{\le m}_i,
A_a^{\le m}]  
= & \partial_j [\Omega_i^{\le m}, A_a^{\le m}] 
- \partial_i [\Omega_j^{\le m}, A_a^{\le m}] \\
& + [\Omega_i^{\le m},\bpartial_a \Omega_j^{\le m-1}] 
- [\Omega_j^{\le m},\bpartial_a \Omega_i^{\le m-1}] \\ 
\equiv   
\bpartial_a[\Omega_i^{\le m},\Omega_j^{\le m-1}] &  
- [\bpartial_a \Omega_i^{\le m}, \Omega_j^{\le m-1}] 
- [\Omega_j^{\le m}, \bpartial_a \Omega_i^{\le m-1}] 
\equiv 0 \mod \frm^m  
\end{align*} 
Thus in particular, the action of 
$\partial_j \Omega^{\le m} - \partial_j \Omega^{\le m}_i$
on $\mE_0\otimes (\mO_B/\frm^m)$ commutes with $A_a^{\le m-1}$.  
On the other hand, 
$(\partial_j \Omega^{\le m}_i - \partial_i \Omega^{\le m}_j)e_0
= \partial_j \be_i - \partial_i \be_j =0$. Since 
$\mE_0\otimes (\mO_B/\frm^m)$ is generated by $e_0$ 
as an $(\mO_B/\frm^m)[A^{\le m-1}_a]$-module, we must have 
$\partial_j \Omega^{\le m} - \partial_j \Omega^{\le m}_i=0$. 
This completes the induction step.

\noindent
{\bf Universal property:}   
Let $\mFh$ be an AQDM over $B'$ with an 
isomorphism $\varphi' \colon \mEh \cong \mFh_{\bb'}$ 
of $\mD$-modules. 
Let $\epsilon^0,\dots,\epsilon^l$ be co-ordinates on $B'$ 
centered at $\bb'$. 
First we construct a map $f\colon B' \rightarrow B$ and 
an isomorphism $\theta\colon \mFh \cong f^*\hmEh$
such that $f(\bb')=\bb$ and $\theta_{\bb'} \circ \varphi' = \varphi$. 
The canonical frame $\Phi$ of $\mEh$ used in the first half of 
the proof induces a frame $\Phi'_{\bb'}$ of $\mFh_{\bb'}$ 
via the isomorphism $\varphi'$. 
\[
\begin{CD}
\mE_0 \otimes \mO^\hbar @>{\Phi}>> \mEh \\
@V{\varphi'_0\otimes \id}VV    @V{\varphi'}VV \\
\mF_{\bb',0} \otimes \mO^\hbar @>{\Phi'_{\bb'}}>> \mFh_{\bb'} 
\end{CD} 
\quad\quad \text{$\varphi_0'\colon \mE_0\rightarrow \mF_{\bb',0}$ 
is the map induced from $\varphi'$.} 
\]
By Theorem \ref{thm:canonicalframe}, 
$\Phi'_{\bb'}$ uniquely extends to a canonical frame 
$\Phi'\colon \mF_{\bb',0}\otimes \mO^\hbar_{B'} \rightarrow \mFh$. 
Let $J'(Q,\epsilon,\hbar)$ be the $J$-function of $\mFh$ 
associated with $\Phi'$ and $\varphi'_0(e_0)\in \mF_{\bb',0}$. 
Define $\hatt^i(Q,\epsilon)\in \mO_{B'}$ by the $\hbar^{-1}$-expansion
of $J'$: 
\[
J' = \varphi'_0(e_0) + \frac{1}{\hbar} 
\sum_{i=0}^s \hatt^i(Q,\epsilon) \varphi'_0(\be_i) + O(\hbar^{-2}). 
\]
Note that $\{\varphi'_0(\be_i)\}_{i=0}^s$ 
is a basis of $\mF_{\bb,0}$. 
The functions $t^i(Q,\epsilon):=\hatt^i(Q,\epsilon)-\hatt^i(Q,0)$
determine a map $f\colon B' \rightarrow B$ such that $f(\bb')=\bb$. 
We define an $\mO_{B'}$-module 
isomorphism $\theta\colon \mFh\rightarrow f^*\hmEh$ by  
\[
\begin{CD}
\theta \colon \mFh @>{\Phi'^{-1}}>> \mF_{\bb',0}\otimes_K \mO^\hbar_{B'}
@>{\varphi'^{-1}_0\otimes \id}>> \mE_0 \otimes_K \mO^\hbar_{B'} 
= f^* \hmEh.  
\end{CD} 
\]
Then $\theta_{\bb'} \circ \varphi'$ coincides with $\Phi^{-1}=\varphi$ 
by the above commutative diagram. 
We show that $\theta$ is a homomorphism of $\mD_{B'}$-modules. 
Let $A_a'(Q,\epsilon)$, $1\le a\le r$ and $\Omega_i'(Q,\epsilon)$, 
$0\le i\le l$ be the connection operators of $\mFh$ associated 
with the canonical frame $\Phi'$. Since $\varphi'$ gives an 
isomorphism of $\mD$-modules, 
\begin{equation}
\label{eq:varphiconjA}
A_a'(Q,0) = \varphi_0'A_a(Q) \varphi_0'^{-1}.
\end{equation} 
Let $L'$ be the fundamental solution of $\mFh$ associated with $\Phi$. 
Then one has 
\[
\hbar\parfrac{}{\epsilon^i} J' = 
\hbar\parfrac{}{\epsilon^i} L'^{-1} \varphi'(e_0) = 
L'^{-1} (\Omega'_i(Q,\epsilon) \varphi'(e_0)) = 
\Omega'_i(Q,\epsilon) \varphi'(e_0) + O(\hbar^{-1}).
\]
Hence, by comparing the leading term of the $\hbar^{-1}$-expansion 
\begin{equation}
\label{eq:Omega_ie_0}
\Omega'_i(Q,\epsilon)\varphi'_0(e_0)=
\sum_{j=0}^s \parfrac{\hatt^j(Q,\epsilon)}{\epsilon^i}
\varphi'_0(\be_j) 
=\varphi'_0
\left(\sum_{j=0}^s \parfrac{t^j(Q,\epsilon)}{\epsilon^i} 
\be_j\right).  
\end{equation} 
On the other hand, the connection operators $f^*A_a(Q,\epsilon)$, 
$1\le a\le r$ and $f^*\Omega_i(Q,\epsilon)$, $0\le i\le l$  
of $f^*\hmEh$ satisfy 
\[
f^*A_a(Q,0) = A_a(Q), \quad 
f^*\Omega_i(Q,\epsilon) e_0 = 
\sum_{j=0}^s \parfrac{t^j(Q,\epsilon)}{\epsilon^i} \be_j. 
\]
But the flat connections $\hbar Q^a\parfrac{}{Q^a}+ A'_a$ and 
$\hbar\parfrac{}{\epsilon^i}+\Omega'_i$ satisfying 
(\ref{eq:varphiconjA}) and (\ref{eq:Omega_ie_0}) are unique; 
this follows from the same argument as the first half. 
Therefore, we conclude
\[
A_a'(Q,\epsilon) = \varphi'_0 A_a(Q,\epsilon) \varphi'^{-1}_0, 
\quad 
\Omega_i'(Q,\epsilon) = \varphi'_0 \Omega_i(Q,\epsilon)\varphi'^{-1}_0.
\]
This shows that $\theta$ is a homomorphism of $\mD_{B'}$-modules. 

Finally, we show that such $f$ and $\theta$ are unique. 
Let $f\colon B'\rightarrow B$ and $\theta\colon \mFh\cong f^*\hmEh$ 
be an arbitrary map and an isomorphism of $\mD_{B'}$-modules 
satisfying $f(\bb')=\bb$ and 
$\theta_{\bb'} \circ \varphi' = \varphi$. 
Let $\hJ(Q,t,\hbar)$ be the $J$-function 
of $\hmEh$ associated with the given trivialization and 
$e_0\in \mE_0=\hmE_{\bb,0}$. 
Let $\hatt^i(Q,t) \in \mO_B$ be the flat co-ordinate 
given by the expansion (see Theorem-Definition \ref{thm-def:affine})
\[
\hJ = e_0 + \frac{1}{\hbar} \sum_{i=0}^s \hatt^i(Q,t) 
\be_i + O(\hbar^{-2}). 
\]
By applying the discussion preceding (\ref{eq:Omega_ie_0}) to 
this $\hJ$, we have 
$\be_i = \Omega_i(Q,t) e_0 = \sum_{j=0}^s 
(\partial \hatt^j(Q,t)/\partial t^i) \be_j$. 
Therefore, $\partial \hatt^j(Q,t)/\partial t^i=\delta^j_i$ and 
\begin{equation}
\label{eq:relation_t_hatt}
\hatt^i(Q,t) = t^i + g^i(Q),  \quad \text{ or } \quad 
t^i = \hatt^i(Q,t)-\hatt^i(Q,0). 
\end{equation}
for some $g^i\in \mO$. 
The given canonical frame of $\hmEh$ induces a canonical frame 
$\Phi'$ of $\mFh$ via $\theta$. 
From $\theta_{\bb'}\circ \varphi' = \varphi$ and 
Theorem \ref{thm:canonicalframe}, it follows that 
this canonical frame $\Phi'$ coincides with the $\Phi'$
in the previous paragraph.  
The $J$-function $J'$ 
of $\mFh$ with respect to $\Phi'$ and 
$\varphi'_0(e_0)\in \mF_{\bb',0}$ is given by 
\[
J' = \varphi'_0(f^*\hJ) 
= \varphi'_0(e_0) + \frac{1}{\hbar}\sum_{i=0}^s 
(f^*\hatt^i) \varphi'_0(\be_i) + O(\hbar^{-2}).    
\]
From this, (\ref{eq:relation_t_hatt}) and 
$f(\bb')=\bb$, it follows that the map $f$ coincides with 
what was given above. 
The uniqueness of $\theta$ follows from that 
$\theta$ intertwines the canonical frame $\Phi'$ of $\mFh$ with  
the given canonical frame of $f^*\hmEh$. 
\end{proof}

\begin{remark}
The {\it uniqueness} of the reconstruction of a big AQDM 
holds in a more general situation. 
For a small AQDM $\mE^\hbar$, $\mE:=\mEh/\hbar \mEh$ 
has the structure of an $\mO[\Qp_1,\dots,\Qp_r]$-module. 
Through a canonical frame, $\mE$ is identified with 
the $\mO[A_1(Q),\dots,A_r(Q)]$-module $\mE_0\otimes \mO$, 
where $A_a(Q)$ is the connection operator associated with 
the canonical frame. 
In the above proof, we used the fact that 
$\mE_0\otimes \mO$ is generated by $e_0$ as an  
$\mO[A_1(Q),\dots,A_r(Q)]$-module. 
Under the following weaker condition, 
\begin{equation}
\label{eq:condition_for_uniqueness}
\text{\it $\mE\otimes_{\mO}\hat{\mO}$ is generated by a single element 
as an $\hat{\mO}[\Qp_1,\cdots,\Qp_r]$-module}  
\end{equation} 
where $\hat{\mO}$ is the quotient field of $\mO$, 
the same argument proves that 
if there exists a big AQDM $\hmEh$ such that $\hmEh_{\bb} \cong \mEh$, 
then $\hmEh$ satisfies the universal property stated in the above theorem. 
In particular, $\hmEh$ is unique if it exists. 
For small quantum cohomology $SQH^*$, 
this condition (\ref{eq:condition_for_uniqueness}) corresponds to 
that $SQH^*$ is generated by $H^2$ at a 
generic value of the quantum parameter $Q$. 
For example, the small quantum cohomology having a tame semisimple point 
satisfies this condition since the Euler vector field (in $H^2$) 
generates the quantum cohomology algebra in a neighborhood of
a tame semisimple point. 
In this case, the uniqueness of the reconstruction agrees with 
the Dubrovin's reconstruction theorem \cite{dubrovin}. 
\end{remark}

\begin{remark} 
Hertling-Manin \cite[Theorem 2.5]{hertling-manin} 
proved a similar reconstruction theorem for (TE)-structures. 
The (TE) structure is endowed with a family of flat connections 
$d+ \frac{1}{\hbar} \Omega(t,\hbar)$
parametrized by $\hbar$ which has poles at $\hbar=0$. 
It also has a flat connection in the $\hbar$-direction 
which corresponds to the grading in the quantum cohomology. 
Our reconstruction theorem differs from theirs 
in that (i) we do not consider differential equations 
in the $\hbar$-direction, 
that (ii) the AQDM here is regular singular  
along $Q^1\cdots Q^r=0$ 
and that (iii) they reconstructed (TE) structures 
in analytic category. 
\end{remark}

We can now formulate {\it generalized mirror transformations}. 
\begin{definition}[Generalized mirror transformation] 
\label{def:GMT}
(i) We begin with a small AQDM $\mEh$  
such that $\mE_0$ is generated by a single element $e_0$ 
as a $K[p_1,\dots,p_r]$-module. 
The equivariant Floer cohomology $FH_{T^2}^*(L_{X/\mV})$
in Section \ref{sec:EFT} gives an example of such small AQDMs. 

(ii) Choose a frame $\Phi_0$ of $\mEh_0$ satisfying 
(\ref{eq:frame_zero}) and 
calculate the canonical frame $\Phi$ of $\mEh$ inducing $\Phi_0$ 
by the Birkhoff factorization of the fundamental solution $L$
(Theorem \ref{thm:canonicalframe}). 
Let $A_a(Q)$, $1\le a\le r$ be the connection operators 
associated with the canonical frame $\Phi$. 

(iii) Reconstruct a big AQDM $\hmEh$ on the base $B$ 
such that $\hmEh_{\bb} \cong \mEh$ as $\mD$-modules 
at an $\mO$-valued point $\bb$ on $B$ (Theorem \ref{thm:reconstruction}). 
More concretely, we take $\hmEh$ to be $\mE_0\otimes \mO_B^\hbar$ 
endowed with flat connections 
$\hbar Q^a \parfrac{}{Q^a}+ A_a(Q,t)$, $1\le a\le r$ and 
$\hbar \parfrac{}{t^i}+ \Omega_i(Q,t)$, $0\le i\le s$. 
Here, $t^0,\dots,t^s$ are co-ordinates on $B$ centered at $\bb$ 
and $A_a$, $\Omega_i$ are $\hbar$-independent. 
Given a $K$-basis $\{\be_i\}_{i=0}^s$ of $\mE_0$,  
$A_a$, $\Omega_i$ are uniquely determined by  
(\ref{eq:commutingconn}), (\ref{eq:integrableconn}), 
$A_a(Q,0)=A_a(Q)$ and $\Omega_i(Q,t)e_0=\be_i$. 

(iv) 
Take a flat co-ordinate system $\hatt^0,\dots,\hatt^s$ on $B$
(Theorem-Definition \ref{thm-def:affine}) 
of the form $\hatt^i =\hatt^i(Q,t)= t^i + g^i(Q)$, $g^i(Q)\in\mO$, 
$g^i(0)=0$ (see (\ref{eq:relation_t_hatt})). 
This new co-ordinate system gives the new connection 
operators $\hA_a(Q,\hatt)$, $\Omega_{\hat{i}}(Q,\hatt)$ 
(by (\ref{eq:connchange_coordchange})) which 
satisfy $\hA_a e_0 = p_a e_0$, $\Omega_{\hat{i}} e_0 = \hwp_i e_0$, 
where $p_a=\hA_a|_{Q=\hatt=0}$ 
and $\hwp_i= \Omega_{\hat{i}}|_{Q=\hatt=0}$. 
\end{definition}

By a generalized mirror transformation, 
the original small AQDM $\mEh$ becomes isomorphic to the 
restriction of $\hmEh$ to the locus $\hatt^i = \hatt^i(Q,0)$, 
where $\hatt^i(Q,t)$ are flat co-ordinates taken in (iv). 
We call this subspace $\{\hatt^i=\hatt^i(Q,0)\}$ of the 
$(Q,\hatt)$-space a {\it locus of $\mEh$}. 
When $\{p_a e_0\}_{a=1}^r$ is part of a $K$-basis of $\mE_0$,  
by Proposition \ref{prop:divisor_AQDM}, we can assume that 
the connection operators $\hA_a$, $\Omega_{\hat{i}}$ in (iv) 
moreover satisfy $\hA_a=\Omega_{\hat{a}}$, $\Omega_{\hat{0}}=\id$ 
and depend only on $Q^1 e^{\hatt^1},\dots,Q^re^{\hatt^r},
\hatt^{r+1},\dots, \hatt^s$. 
Then we can take 
\[
t^0,\quad \hq^a := Q^a e^{\hatt^a} \quad 1\le a\le r, \quad 
t^j, \quad r+1\le j\le s  
\]
to be {\it effective parameters} of $\hmEh$. 
In these effective parameters, the locus of the original small AQDM 
$\mEh$ is given by (see Figure \ref{fig:embedding})
\[
\hatt^0 = \hatt^0(Q,0), \quad 
\hq^a = Q^a e^{\hatt^a(Q,0)} \quad 1\le a\le r, \quad 
\hatt^j = \hatt^j(Q,0) \quad r+1\le j\le s. 
\] 
In ordinary mirror transformations, the locus of $\mEh$ is contained 
in the small quantum cohomology locus defined by 
$\hatt^{r+1}=\cdots=\hatt^{s}=0$ and $\hq^a=Q^a e^{\hatt^a(Q,0)}$
can be interpreted as a co-ordinate change on the $q$-space. 
The specific characteristic in generalized case 
is that the locus of $\mEh$ is not necessarily contained in 
the small quantum cohomology locus in the effective parameter 
space.  
\begin{figure}[htbp]
\begin{center}
\begin{picture}(400,100)
\put(80,30){\vector(1,0){200}}
\put(310,30){\makebox(0,0){$\hq^1,\dots,\hq^r$}}
\put(125,0){\vector(0,1){100}}
\put(165,95){\makebox(0,0){$\hatt^0,\hatt^{r+1},\dots,\hatt^s$}}
\put(295,71){\makebox(0,0){locus of $\mEh$}}
\thicklines 
\qbezier(100,0)(150,75)(210,45)
\qbezier(210,45)(250,25)(280,60)
\end{picture} 
\end{center} 
\caption{The effective parameter space of $\hmEh$}
\label{fig:embedding}
\end{figure}

We now discuss a relationship among the $J$-functions 
of $\mEh$ and $\hmEh$ and the locus of $\mEh$ in $B$. 
Let $\Phi_{\rm ini}$ be an arbitrary frame of $\mEh$. 
Let $L_{\rm ini}$ be the fundamental solution of $\mEh$ associated with 
$\Phi_{\rm ini}$ and $I(Q,\hbar):=L^{-1}_{\rm ini} e_0$ 
be the $J$-function of $\mEh$ associated with $e_0\in \mE_0$. 
The canonical frame $\Phi$ calculated in step (ii) 
is given by $\Phi=\Phi_{\rm ini}\circ G$,  
where $G=L_{{\rm ini},+}$ is the positive part  
of the Birkhoff factorization (\ref{eq:Birkhoff}) of $L_{\rm ini}$.  
(Here, $\Phi$ and $\Phi_{\rm ini}$ induce the same frame 
of $\mEh_0$.) 
By Proposition \ref{prop:Jgenerates_AQDM}, 
$\mEh$ is generated by $\Phi_{\rm ini}(e_0)$ as a $\mD$-module.  
Hence there exists $V(Q,\Qp,\hbar)\in \mD$ such that 
$V(Q,\Qp,\hbar)\Phi_{\rm ini}(e_0)= \Phi(e_0)=\Phi_{\rm ini}(Ge_0)$, 
\emph{i.e.} $V(Q,\smnabla^\hbar_a,\hbar) e_0=G e_0$.  
Note that we can assume $V(0,\Qp,\hbar)=1$ since 
$G|_{Q=0} =L_{{\rm ini},+}|_{Q=0}=\id$.  
Then by (\ref{eq:gaugetransf_fundsol}), 
the $J$-function $J$ of $\mEh$ associated with 
the new frame $\Phi$ and $e_0$ is given by
\begin{align}
\label{eq:fromItoJ}
J(Q,\hbar) &= L^{-1}_{\rm ini}(Ge_0) \\
\label{eq:fromItoJ_diff}
&= V(Q,\hbar Q^a\parfrac{}{Q^a}+p_a,\hbar) I(Q,\hbar).  
\end{align}
The big AQDM $\hmEh$ reconstructed in step (iii) is naturally 
equipped with a canonical frame (compatible with $\Phi$).  
Let $\hJ(Q,\hatt,\hbar)$ be the $J$-function of $\hmEh$ with 
respect to the given canonical frame, $e_0\in \mE_0=\hmEh_{\bb,0}$
and flat co-ordinates $\hatt^i=\hatt^i(Q,t)$ in (iv). 
(Recall that $\bb$ is given by $t=0$). 
Then it is easy to see that 
\begin{equation}
\label{eq:fromJtohJ}
\hJ(Q,\hatt(Q,0),\hbar) = J(Q,\hbar). 
\end{equation} 
The equations (\ref{eq:fromItoJ}), (\ref{eq:fromItoJ_diff}) and 
(\ref{eq:fromJtohJ}) 
give a relationship between the two $J$-functions $I$ and $\hJ$. 
Moreover, since $\hatt^i$ can be read off from the 
$\hbar^{-1}$-expansion of $\hJ$ (Theorem-Definition \ref{thm-def:affine}), 
we have by (\ref{eq:fromJtohJ}), (\ref{eq:fromItoJ}), 
$G=L_{{\rm ini},+}$ and (\ref{eq:calofLpos}),    
\begin{align*}
\sum_{j=0}^s \hatt^j(Q,0) (\hwp_je_0) 
= \Res_{\hbar=0} d\hbar \left \{
\sum_{k=0}^\infty (\id- L_{\rm ini}^{-1}\circ\pi_+)^k I(Q,\hbar)\right\}.  
\end{align*}
This formula calculates the locus of $\mEh$ perturbatively. 
For a nef complete intersection in toric varieties, 
the $I$-function takes of the form 
$I(Q,\hbar) = f(Q)e_0+O(\hbar^{-1})$. 
In this case, the left-hand side is simplified to the form 
$\Res_{\hbar=0}d\hbar (I(Q,\hbar)/f(Q))$. 
This recovers the original mirror transformation in 
\cite{givental-mirrorthm-projective, givental-mirrorthm-toric}. 

\subsection{Reconstruction via generators}
\label{subsec:reconst_gen}
In order to obtain the reconstruction, in view of its uniqueness, 
we only need to find at least one big AQDM 
which is an extension of the given small AQDM. 
Here, we consider the extension of a $\mD$-module 
via generators or $J$-functions. 

Let $\qmO^\hbar$ be a subring of $\mO^\hbar$ and 
set $\qmD:=\qmO^\hbar\langle \Qp_1,\dots,\Qp_r\rangle\subset \mD$.  
Let $\mJ(Q,\hbar)$ be an element of a suitable function space 
$\frF$ of functions in $Q^1,\dots,Q^r$ and $\hbar$. 
We assume that $\qmD$ acts on $\frF$ and 
that $\mJ(Q,\hbar)$ generates a $\qmD$-module 
$\qmEh=\qmD \mJ(Q,\hbar)\subset \frF$ 
which is finitely generated as an $\qmO^\hbar$-module. 
Here, $\qmEh$ is isomorphic to $\qmD/\mathsf{I}$ where $\mathsf{I}$ 
is the left ideal consisting of an element $f(Q,\Qp,\hbar)\in \qmD$ 
which annihilates $\mJ(Q,\hbar)$ 
(\emph{c.f.} Proposition \ref{prop:Jgenerates_AQDM}). 
We also set $\qmO^\hbar_B:=\qmO^\hbar[\![t^0,\dots,t^s]\!]$ and 
$\qmD_B:= \qmO_B^\hbar\langle
\Qp_1,\dots,\Qp_r,\bwp_0,\dots,\bwp_s\rangle \subset \mD_B$. 
Let $\qmD_B$ act on $\frF[\![t]\!]:=\frF[\![t^0,\dots,t^s]\!]$ 
by $\bwp_i \mapsto \hbar \partial /\partial t^i$. 
\begin{proposition}
\label{prop:reconst_gen} 
For a $\qmD_B$-module $\hqmEh\subset \frF[\![t]\!]$ 
generated by $\hmJ\in \frF[\![t]\!]$ such that 
$\hmJ(Q,t,\hbar) = \mJ(Q,\hbar)+ O(t)$, 
there exists an isomorphism $\hqmEh/\sum_{i=0}^s t^i \hqmEh
\cong \qmEh$ of $\qmD$-modules which sends $[\hmJ]$ to $\mJ$ 
if and only if 
$\hmJ(Q,t,\hbar)$ satisfies 
\begin{equation}
\label{eq:diff_t}
\hbar\parfrac{}{t^i} \hmJ(Q,t,\hbar) \in 
\qmO_B^\hbar\langle \Qp_1,\dots,\Qp_r, t^0\bwp_0,\dots,t^s\bwp_s 
\rangle \hmJ(Q,t,\hbar).  
\end{equation} 
In this case, $\hqmEh$ is a finitely generated $\qmO_B^\hbar$-module. 
If moreover $\qmEh$ is free over $\qmO^\hbar$, 
then $\hqmEh$ is free over $\qmO_B^\hbar$. 
\end{proposition} 
\begin{proof}
First we prove the ``if part".  
Set $\qmD'_B:=\qmO^\hbar_B
\langle \Qp_1,\dots,\Qp_r, t^0\bwp_0,\dots,t^s\bwp_s 
\rangle\subset \qmD_B$. 
It follows from (\ref{eq:diff_t}) that  
$\hqmEh$ is generated by $\hmJ$ as a $\qmD'_B$-module. 
Therefore one can define a surjective homomorphism of $\qmD$-modules:
\[
\hqmEh \ni f(Q,t,\Qp,t\bwp, \hbar)\hmJ \longmapsto 
f(Q,0,\Qp,0,\hbar)\mJ 
\in \qmEh, \quad f(Q,t,\Qp,t\bwp, \hbar) \in \qmD'_B.  
\]
We need to show that the kernel equals $\sum_{i=0}^s t^i \hqmEh$. 
To see this, it suffices to show that for $f \in \qmD$ 
satisfying $f \mJ=0$, there exists $g\in \qmD'_B$ 
such that $g|_{t=0}= f$ and $g\hmJ=0$. 
We will prove this by induction on $s$. 
As an inductive hypothesis, 
we can assume that there exists $g\in \qmD_B'$ such that 
$g|_{t=0} =f$ and $(g\hmJ)|_{t^s=0} = 0$. 
Let $T_i \mJ\in \qmEh$, $1\le i\le N$ be generators 
of $\qmEh$ as an $\qmO^\hbar$-module, where $T_i\in \qmD$. 
Let $n_0$ be the maximum of the degrees 
of $T_i$, $1\le i\le N$ with respect to $\Qp_1,\dots,\Qp_r$. 
Let $\qmD''_{B}$ be 
the subspace of $\qmO_B^\hbar\langle \Qp_1,\dots,\Qp_r\rangle
\subset\qmD'_B$ consisting of all elements which are 
of degree less than or equal to $n_0$ in $\Qp_1,\dots,\Qp_r$. 
We show by induction on $m$ that if $g \hmJ = O((t^s)^m)$ for 
$g\in \qmD'_{B}$, there exists $g_m\in (t^s)^m \qmD''_{B}$ 
such that $(g+g_{m})\hmJ = O((t^s)^{m+1})$.
Let $m$ be zero. There exists $c_i\in \qmO^\hbar$ 
such that $(g \hmJ)|_{t=0} = \sum_{i=1}^N c_i T_i \mJ$. 
Thus we can take $g_0$ to be $-\sum_{i=1}^Nc_i T_i\in \qmD_B''$. 
Let $m$ be positive. 
Set $\partial_s := \partial/\partial t^s$. 
First note that we have for a function $c(t^s)$ in $t^s$, 
\[
[\partial_s, c(t^s) (\hbar t^s \partial_s)^n ] 
= c'(t^s) (\hbar t^s\partial_s)^n 
+ c(t^s) \sum_{k=0}^{n-1} \binom{n}{k} 
(\hbar t^s \partial_s )^k \hbar^{n-k-1} (\hbar \partial_s). 
\]
This means $[\partial_s, \qmD_B']\subset \qmD_B' + 
\qmD_B' \hbar \partial_s$. 
Therefore, we can write $[\partial_s, g]= g' + g'' \hbar \partial_s$
for some $g',g''\in \qmD_B'$. 
By (\ref{eq:diff_t}), there exists 
$F \in \qmD'_B$ such that 
$(\hbar\partial_s) \hmJ = F \hmJ$. 
We calculate 
\begin{align*}
\hbar \partial_s (g \hmJ) & = \hbar(g'+g''\hbar\partial_s) \hmJ 
+ g \hbar \partial_s \hmJ = (\hbar (g'+g''F) + gF) \hmJ \\ 
& = (\hbar (g'+g''F) + [g,F]) \hmJ 
+ F g \hmJ 
= \hbar \phi \hmJ + O((t^s)^m) 
\end{align*} 
where we set $\phi = g' + g''F + [g,F]/\hbar$. 
Note that $\phi \in \qmD'_B$ 
because the commutator $[g,F]$ is divisible by $\hbar$. 
Since $\phi \hmJ = O((t^s)^{m-1})$, by the induction hypothesis, 
there exists $\phi_{m-1}\in (t^s)^{m-1}\qmD''_{B}$ 
such that $(\phi+\phi_{m-1}) \hmJ = O((t^s)^m)$. 
We can take $g_{m}\in (t^s)^m\qmD''_{B}$ 
such that $\phi_{m-1}= \partial g_{m}/\partial t^{s}$. 
Then we have 
\[
\hbar \partial_s ((g+g_{m})\hmJ) = \hbar 
(\phi + \phi_{m-1})  \hmJ  + g_{m} \hbar\partial_s \hmJ + 
O((t^s)^m) = O((t^s)^m).  
\] 
Thus $(g+g_m)\hmJ = O((t^s)^{m+1})$. 
This completes the proof of the ``if part". 

Next we show that 
$\hqmEh$ is finitely generated over $\qmO_B^\hbar$ 
when $\hqmEh/\frm \hqmEh \cong \qmEh$, where 
$\frm=\sum_{i=0}^s t^i \qmO_B^\hbar$.  
It follows by an elementary argument that 
$\{T_i\hmJ\}_{i=1}^N$ generates $\hqmEh$ over $\qmO_B^\hbar$ 
when $\{T_i\mJ\}_{i=1}^N$ generates $\qmEh$ over $\qmO^\hbar$. 
Here, we need the fact that $\qmO_B^\hbar$ is complete and 
$\hqmEh$ is Hausdorff for their $\frm$-adic topology
($\bigcap_{n\ge 0} \frm^n \hqmEh \subset 
\bigcap_{n\ge 0}\frm^n \frF[\![t]\!]=\{0\}$). 
See for instance 
\cite[Corollary 2, \textsection 3, Ch.VIII]{zariski-samuel}. 

Next we show that $\hqmEh$ is free over $\qmO_B^\hbar$ 
if $\qmEh$ is free over $\qmO^\hbar$. 
We take $\{T_i\mJ\}$ to be a free basis over $\qmO^\hbar$. 
We claim that $\{T_i\hmJ\}$ is a free basis over $\qmO_B^\hbar$. 
Since we know that $T_i\hmJ$ generate $\hqmEh$ over $\qmO_B^\hbar$, 
we have $(\hbar\partial /\partial t^j) (T_i\hmJ) = 
\sum_{k=1}^N \Omega_{ji}^k (T_k\hmJ)$ for some 
$\Omega_{ji}^k \in \qmO_B^\hbar$. 
For $c\in \qmO_B^\hbar$, let $\ord(c):=\sup\{n; c\in \frm^n\}$ 
denote the order of zero at $t=0$. 
Suppose we have $\sum_{i=1}^N c_i(T_i\hmJ)=0$ and set 
\[
\nu := \min_{1\le i\le N} ( \ord(c_i)) = \ord (c_{i_0}).  
\]
Because $T_i\mJ$ is a basis of $\qmEh$ 
and $\sum_{i=1}^N c_i|_{t=0} T_i \mJ =0$, 
we have $c_i|_{t=0}=0$. Thus $\nu>0$. 
If $0<\nu <\infty$, we calculate 
\begin{align*}
0= \hbar \parfrac{}{t^j} \sum_{i=1}^N c_i (T_i\hmJ) 
= \sum_{k=1}^N (\hbar \parfrac{c_k}{t^j} +\sum_{i=1}^N 
 c_i \Omega_{ji}^k) (T_k\hmJ).  
\end{align*} 
Here, $\ord(\hbar (\partial/\partial t^j) c_{i_0}
+\sum_{i=1}^N c_i \Omega_{ji}^{i_0})<\ord c_{i_0}=\nu$ 
for some $0\le j\le s$. Thus we get a relation with smaller $\nu$.  
Repeating this, we have contradiction 
since $\nu$ must be positive. 
Therefore $\nu=\infty$, \emph{i.e.} $c_i=0$. 

Finally we prove the ``only if" part. 
As we have shown, $\hqmEh/\frm \hqmEh \cong \qmEh$ implies that 
$\hqmEh$ is generated by $T_i\hmJ$, $1\le i\le N$ over $\qmO_B^\hbar$. 
Thus, we can find elements $c_i\in \qmO_B^\hbar$ such that 
$\hbar(\partial/\partial t^i) \hmJ= \sum_{i=1}^N c_i T_i \hmJ 
\in \qmD'_B \hmJ$.
\end{proof} 

\begin{example}
Let $\mEh$ be the small AQDM such that 
$\mE_0$ is generated by $e_0\in \mE_0$ as a $K[p_1,\dots,p_r]$-module.  
Let $J(Q,\hbar)$ be the $J$-function associated with 
some frame and $e_0$. 
We take $\frF$ to be $\mE_0^\hbar\otimes 
K\{\hbar,\hbar^{-1}\}\!\}[\![Q]\!]$ 
and let $\mD$ act on $\frF$ by 
$\Qp_a \mapsto \hbar Q^a(\partial/\partial Q^a)+p_a$. 
By Proposition \ref{prop:Jgenerates_AQDM}, 
$J(Q,\hbar)\in \frF$ generates $\mEh$.  
Let $\{T_i(p_1,\dots,p_r) e_0\}_{i=0}^s$ 
be a $K$-basis of $\mE_0$. 
Then, the following function $\hJ$ satisfies (\ref{eq:diff_t}) 
and formally reconstructs a big AQDM: 
\[
\hJ(Q,t,\hbar)  = 
\exp\left(\frac{1}{\hbar}
\sum_{i=0}^s t^i T_i(\hbar Q^a\parfrac{}{Q^a}+p_a)
\right) J(Q,\hbar)  \in \frF[\![t]\!]
\]
Note that $t^i$ here is not necessarily a flat co-ordinate. 
\end{example}

\begin{example}
We construct a formal mirror 
corresponding to the big quantum cohomology of $\Proj^n$. 
Following Givental \cite{givental-ICM,givental-mirrorthm-toric}, 
we take a mirror of the small quantum cohomology of $\Proj^n$ 
to be the oscillatory integral  
\begin{equation}
\label{eq:osci_int_P2}
\mJ(Q,\hbar) = \int_{\Gamma\subset (\C^*)^n} 
\exp(\frac{1}{\hbar} (x_1+\cdots+x_n+\frac{Q}{x_1x_2\cdots x_n})) 
\frac{dx}{x}  
\end{equation}
where $\Gamma$ is a non-compact Morse cycle of 
the Morse function $(x_1,\dots,x_n) 
\mapsto \Re(x_1+\cdots+x_n+\frac{Q}{x_1\cdots x_n})$ 
and $dx/x = \prod_{i=1}^n (dx_i/x_i)$. 
This $\mJ$ is a multi-valued function in $Q$ and $\hbar$. 
We take $\frF$ to be the space of multi-valued 
analytic functions on $(Q,\hbar)\in \C^*\times \C^*$, \emph{i.e.} 
single-valued functions on the universal cover of $\C^*\times\C^*$. 
Then $\mJ\in \frF$. 
We also take $\qmO^\hbar$ to be the 
subspace of $\mO^\hbar=\C[\hbar][\![Q]\!]$ (here $K=\C$) 
consisting of convergent functions for all $Q,\hbar\in\C$. 
Let $\qmD:=\qmO^\hbar\langle \Qp\rangle$ 
act on $\frF$ by $\Qp \mapsto \hbar Q(\partial/\partial Q)$. 
The oscillatory integral (\ref{eq:osci_int_P2}) 
satisfies the differential equation 
\begin{equation}
\label{eq:rel_smallQDM_P2}
((\hbar \bpartial)^{n+1} - Q ) \mJ(Q,\hbar) = 0, \quad 
\bpartial = Q (\partial/\partial Q) 
\end{equation} 
and generates the small QDM of $\Proj^n$: $SQDM(\Proj^n) \cong 
\mO^\hbar \otimes_{\qmO^\hbar}
(\qmD/\langle(\hbar\bpartial)^{n+1}-Q \rangle)$. 
Let $t^0,t^1,\dots,t^n$ be co-ordinates on $B$ 
and consider $\hmJ$ in $\frF[\![t^0,\dots,t^n]\!]$: 
\[
\hmJ(Q,t,\hbar) = 
\int_{\Gamma\in (\C^*)^2} \exp (\frac{1}{\hbar}(t^0 + 
\sum_{i=1}^n x_i + \frac{Q}{x_1\cdots x_n}
+\sum_{i=1}^n t^i x^1\cdots x^i))
\frac{dx}{x} 
\] 
Here, the integrand should be expanded as a power series in $t$. 
We can easily show that $\hmJ$ satisfies 
\begin{gather}
\label{eq:deform_bigQDMP2}
(\mR_1 \mR_2 \cdots \mR_n \hbar \bpartial - Q) 
\hmJ = 0,  \\ 
\label{eq:diff_t_check}
\hbar \partial_0 \hmJ = \hmJ, \quad 
\hbar \partial_i \hmJ = \mR_1 \mR_2 \cdots \mR_{i} \hmJ \quad 
(i\ge 1). 
\end{gather} 
where 
$\mR_i := \hbar \bpartial - \hbar t^{i}\partial_{i} 
- \cdots -\hbar t^n \partial_n$ and $\partial_i = \partial/\partial t^i$.  
The equation (\ref{eq:deform_bigQDMP2}) 
is a deformation of the relation (\ref{eq:rel_smallQDM_P2})
and (\ref{eq:diff_t_check}) shows that $\hmJ$ 
satisfies the condition (\ref{eq:diff_t}). 
Therefore, the $\qmD_B$-module $\hqmEh$ generated by $\hmJ$ 
is a free $\qmO^\hbar_B$-module of rank $n+1$ and an 
extension of the $\qmD$-module generated by $\mJ$. 
By the uniqueness of the reconstruction, 
$\mO_B^\hbar\otimes_{\qmO^\hbar_B} \hqmEh$ 
is canonically isomorphic to the big QDM of $\Proj^n$. 

On the other hand, 
the Jacobi ring of the phase function $t^0 + 
\sum_{i=1}^n x_i + Q/(x_1\cdots x_n)
+\sum_{i=1}^n t^i x^1\cdots x^i$ is of dimension 
$2n>\dim H^*(\Proj^n)$ for non-zero generic $t^i$. 
Thus, it seems to be important to treat $\hmJ$ 
as a formal power series in $t^i$. 
This formal mirror construction for big quantum cohomology 
will be generalized to toric varieties in general.   
We hope to discuss this problem elsewhere. 
In the literature, starting from a Laurent polynomial, 
Barannikov \cite{barannikov-proj} (for $\Proj^n$) 
and Douai and Sabbah \cite{douai-sabbah-I, douai-sabbah-II} 
(for weighted projective spaces or for a general tame function)  
constructed a Frobenius manifold of big quantum cohomology.  
It would be interesting to study a relationship between their 
methods and ours. 
\end{example} 

\section{The Proof of Generalized Mirror Transformations}
\label{sec:main}  
In this section, we prove that 
the big AQDM reconstructed from the equivariant Floer 
cohomology $FH_{T^2}(L_{X/\mV})$ by the generalized 
mirror transformation (Definition \ref{def:GMT})  
is isomorphic to the big QDM for a pair $(X,\mV)$.

\subsection{Review of quantum Lefschetz theorem}
\label{subsec:review_QL}
We review the quantum Lefschetz theorem by 
Coates and Givental \cite{coates-givental}. 
We use the notation in Section \ref{sec:QDM}. 
Let $M$ be a smooth projective variety and 
$\mV$ be a vector bundle on $M$. 
Let $\bc$ be a general multiplicative characteristic 
class (\ref{eq:general_c}). 
We take $K=\C[\![\bs]\!]$ to be the ground ring 
and consider the Gromov-Witten theory twisted by $\mV$ and $\bc$. 
The genus zero twisted Gromov-Witten potential $\mF_{\bs}$ 
is defined to be
\[
  \mF_{\bs}(\bt_0,\bt_1,\bt_2,\dots)=
  \sum_{n=0}^\infty \sum_{d\in \Lambda}
  \frac{Q^d}{n!}
  \corrV{\bt(\psi_1),\cdots,\bt(\psi_n)}_{d}, 
 \quad \bt(\psi_i) = \sum_{k=0}^\infty \bt_k \psi_i^k,  
\]
where $\bt_i$ is a $H^*(M)$-valued co-ordinate. 
Introduce an infinite dimensional space $\mH$ over $K[\![Q]\!]$ 
with a symplectic form $\Omega_\bs$:    
\begin{align*}
\mH & := H^*(M)\otimes K\{\hbar,\hbar^{-1}\}\!\}[\![Q]\!], \\
\Omega_{\bs}({\bof}(\hbar),\bg(\hbar))
&:=\Res_{\hbar=0}\pairV{\bof(-\hbar)}{\bg(\hbar)}d\hbar, 
\quad \bof,\bg\in \mH, 
\end{align*}
where $\pairV{\cdot}{\cdot}$ is the twisted 
Poincar\'{e} pairing (\ref{eq:twisted_Poincare}). 
The symplectic space $\mH$ is decomposed 
into two isotropic subspaces $\mH_+$, $\mH_-$ 
as $\mH=\mH_+\oplus \mH_-$: 
\[
\mH_{+}=H^*(M)\otimes K\{\hbar\}[\![Q]\!], \quad 
\mH_{-}=\hbar^{-1}H^*(M)\otimes K\{\!\{\hbar^{-1}\}\!\}[\![Q]\!].
\] 
We call this a {\it polarization} of $\mH$. 
By the polarization, we can regard $\mH$ as 
the cotangent bundle $T^*\mH_+$ of $\mH_+$ 
with $\mH_-$ being identified with a fiber of $T^*\mH_+$.  
We write a general element in $\mH_+$ as 
\[
\bq(\hbar):=\bq_0+\bq_1\hbar+\bq_2\hbar^2+\cdots \in \mH_+, 
\]
where $\bq_i$ is a $H^*(M)$-valued co-ordinate. 
We relate the co-ordinates $\bt_n$ to the co-ordinate $\bq_n$ 
by the following {\it dilaton shift}: 
\[ 
\bq(\hbar)=\bt(\hbar)-1 \hbar,  \quad \text{\emph{i.e.}} \quad 
\bq_n = \bt_n - \delta_{n,1} \hbar,
\]
where $1\in H^0(M)$ is the ring identity 
and $\bt(\hbar)=\sum_{n=0}^\infty \bt_n \hbar^n$.
We regard the potential $\mF_\bs(\bt_0,\bt_1,\dots)$ 
as a function on the formal neighborhood of 
$\bq(\hbar)=-\hbar$ in $\mH_+$. 
The differential $d\mF_\bs$ of the potential $\mF_\bs$ 
gives a section of $T^*\mH_+$ and its graph 
defines a formal germ $\mL_\bs \subset \mH\cong T^*\mH_+$ 
of a Lagrangian submanifold. 
More explicitly,  $\mL_\bs$ consists of 
points in $\mH$ of the form: 
\begin{equation}
\label{eq:pointoncone}
d_{\bq}\mF_{\bs} = 
\bt(\hbar) - \hbar + \sum_{n=0}^\infty \frac{1}{(-\hbar)^{n+1}}
\sum_{0\le j,k\le s}  
\parfrac{\mF_\bs(\bt_0,\bt_1,\dots)}{\bt_n^j} g^{jk} \frac{p_k}{\bc(\mV)}  
\end{equation}
where $\{p_j\}_{j=0}^s$ is a basis of $H^*(M)$ and  
we write $\bt_n = \sum_{j=0}^s \bt_n^j p_j$, 
$g_{ij}=\int_M p_i\cup p_j$ and $(g^{ij})_{0\le i,j\le s}$ 
is the matrix inverse to $(g_{ij})_{0\le i,j\le s}$. 
The $J$-function in Definition \ref{def:J-funct} 
is considered to be a family of elements on the cone. 
In fact, by (\ref{eq:J-function}) and (\ref{eq:pointoncone}), 
the $J$-function $J_\bs$ of the twisted theory satisfies  
\[
-\hbar J_\bs(Q,\tau,-\hbar)
= d_{-\hbar+\tau} \mF_{\bs} 
= \mL_\bs \cap \{-\hbar +\tau + \mH_-\}. 
\]
Thus, the derivatives of the $J$-function 
(recall that $\tau=\sum_{i=0}^s t^i p_i$) 
\begin{equation}
\label{eq:derivatives_J}
-\hbar\parfrac{}{t^j} J_\bs(Q,\tau,-\hbar) = L(Q,\tau,-\hbar)^{-1} p_j  
\end{equation} 
are tangent vectors to $\mL_\bs$.   
Coates and Givental \cite{coates-givental} observed that 
\begin{itemize}
\item The tangent space $\bL_\tau$ of $\mL_\bs$ at 
$d_{-\hbar+\tau}\mF_\bs$ a free module over $K\{\hbar\}[\![Q]\!]$ 
generated by the derivatives (\ref{eq:derivatives_J}) of the $J$-function:
\[ 
\bL_\tau= L(Q,\tau,-\hbar)^{-1}(H^*(M)\otimes K\{\hbar\}[\![Q]\!]) 
\subset \mH 
\]

\item We have $\hbar \bL_\tau \subset \mL_\bs$ and the tangent 
space to $\mL_\bs$ at any point on $\hbar \bL_\tau$ equals $\bL_\tau$. 
The Lagrangian $\mL_\bs$ is 
ruled by these tangent spaces: 
\begin{equation}
\label{eq:ruledcone}
\mL_\bs = \bigcup_{\tau\in H^*(M)} \hbar \bL_\tau
\end{equation} 
\end{itemize} 
In particular, $\mL_\bs$ is a {\it Lagrangian cone}.  
The family $\{\bL_{\tau};\tau\in H^*(M)\}$ of subspaces in 
$\mH$ forms a {\it semi-infinite variation of Hodge structures}  
in the sense of Barannikov \cite{barannikov-qpI}. 
In \cite{givental-symplecticgeom}, Givental explained that 
the above geometric properties of $\mL_\bs$ correspond to 
the string equation, the dilaton equation and 
the topological recursion relations for 
the genus zero Gromov-Witten invariants. 

\begin{theorem}[Coates-Givental \cite{coates-givental}]
\label{thm:coates-givental-cone}
The Lagrangian cone $\mL_\bs\subset (\mH,\Omega_\bs)$ 
of the twisted Gromov-Witten theory and 
that $\mL_0\subset (\mH,\Omega_0)$ of the untwisted theory 
are related by a symplectic transformation 
$\triangle \colon (\mH,\Omega_0) \rightarrow (\mH,\Omega_\bs)$: 
\[
\mL_\bs = \triangle \mL_0, \quad 
\triangle = \bc(\mV)^{-1/2} 
\exp\left(\sum_{m,l\ge 0} s_{2m-1+l}
\frac{B_{2m}}{(2m)!}\ch_l(\mV) \hbar^{2m-1}\right),
\]
where we set $s_{-1}=0$ and 
Bernoulli numbers $B_{2m}$ are defined by 
$x/(1-e^{-x})=x/2+\sum_{m=0}^\infty B_{2m}x^{2m}/(2m)!$. 
\end{theorem}


The $J$-function is obtained from the cone $\mL_\bs$ as a 
slice $\mL_\bs \cap \{-\hbar+\tau+\mH_-\}$ and 
its derivatives (\ref{eq:derivatives_J}) recover the fundamental 
solution $L$ of the big QDM. 
The Dubrovin connection and the quantum product 
can be obtained from the differential 
equations (\ref{eq:fundsol_big}), (\ref{eq:fundsol_small}) for $L$, 
so the big QDM with its canonical frame  
is determined by the cone $\mL_\bs$. 

Let $\mV$ be the sum $\mV_1\oplus \cdots \oplus \mV_l$ 
of line bundles and $\bc$ 
be the $S^1$-equivariant Euler class $\be$ 
(\ref{eq:equivariant_Eulerclass}). 
In this case, Coates and Givental found another family of 
elements lying on the cone $\mL_\bs$. 
Recall that the ground ring $K$ is $\C(\!(\lambda^{-1})\!)$
when $\bc=\be$.  
Set $v_i:=c_1(\mV_i)$ and 
let $J(Q,\tau,\hbar)=\sum_{d\in \Lambda} J_d(\tau,\hbar) Q^d$ 
be the untwisted $J$-function of $M$. 
The {\it hypergeometric modification}
$\Itw(Q,\tau,\hbar)$ of the $J$-function is defined to be
\begin{equation}
\label{eq:hypergeometric_modif}
\Itw(Q,\tau,\hbar):=\sum_{d\in \Lambda}  
\prod_{i=1}^l
\frac{\prod_{\nu=-\infty}^{\pair{v_i}{d}}(v_i+\nu\hbar+\lambda)}
     {\prod_{\nu=-\infty}^{0}(v_i+\nu\hbar+\lambda)}\cup 
J_d(\tau,\hbar) Q^d.
\end{equation} 

\begin{theorem}[Coates-Givental \cite{coates-givental}]  
\label{thm:coates-givental-hyp}
The family $\tau\mapsto 
-\hbar \Itw(Q,\tau,-\hbar)$ 
is lying on the Lagrangian cone of the 
$(\mV,\be)$-twisted Gromov-Witten theory of $M$. 
\end{theorem}

As we see below, 
the hypergeometric modification $\Itw$ recovers 
the cone $\mL_\bs$ as $J_\bs$ does.   
The above theorem and (\ref{eq:ruledcone}) imply that 
\[
-\hbar \Itw(Q,\tau,-\hbar) = -\hbar 
L(Q,\htau,-\hbar)^{-1} v(Q,\tau,-\hbar) \in \hbar\bL_{\htau}
\]
for some $v\in H^*(M)\otimes K\{\hbar\}[\![Q,t]\!]$ 
and $\htau = \htau(Q,\tau) \in H^*(M)\otimes K[\![Q,t]\!]$. 
Because $L(0,\htau,\hbar) = e^{\htau/\hbar}$ and 
$\Itw(0,\tau,\hbar)= e^{\tau/\hbar}$, we have 
$v(0,\tau,\hbar) = e^{(\tau-\htau(0,\tau))/\hbar}$. 
But $v$ does not contain negative powers in $\hbar$, so 
it follows that $\htau(0,\tau)=\tau$ and $v(0,\tau,\hbar)=1$. 
Hence, $\tau\mapsto \htau=\htau(Q,\tau)$ 
is an invertible co-ordinate change on $H^*(M)$.  
The derivatives of $\Itw$  
\begin{align*}
-\hbar \parfrac{}{t^j} \Itw(Q,\tau,-\hbar) 
= L(Q,\htau,-\hbar)^{-1} \nabla^{-\hbar}_j  
v(Q,\tau,-\hbar),
\end{align*} 
where $\nabla^{-\hbar}_j = 
-\hbar(\partial/\partial t^j) + 
(\partial \htau/\partial t^j)*$ is the Dubrovin connection, 
generate $\bL_{\htau}$ over $K\{\hbar\}[\![Q]\!]$ 
since $\nabla^{-\hbar}_j v(Q,\tau,-\hbar)=p_j + O(Q)$.  
Therefore the derivatives of $\Itw$ 
determines the cone $\mL_\bs$ by (\ref{eq:ruledcone}). 
These imply the following corollary 
(\emph{c.f.} Proposition \ref{prop:Jgenerates}, 
Section \ref{subsec:reconst_gen}):

\begin{corollary} 
\label{cor:modif_generates}
(1) Let $B$ and $\mD_B$ be as in Section \ref{sec:AQDM}. 
Let $\mD_B$ act on $H^*(M)\otimes K\{\hbar,\hbar^{-1}\}\!\}[\![Q,t]\!]$ 
by $\Qp_a\mapsto \hbar Q^a \partial/\partial Q^a+p_a$ and 
$\bwp_i \mapsto \hbar \partial/\partial t^i$. 
The $\mD_B$-module generated by $\Itw(Q,\tau,\hbar)$
(where $\tau=\sum_{i=0}^s  t^i p_i$) 
is isomorphic to $QDM_\be(M,\mV)$ 
under some co-ordinate change $\htau=\htau(Q,\tau)$. 
Here, $\htau$ is the natural $H^*(M)$-valued co-ordinate 
on the base of the QDM. 
We have $\htau(0,\tau)=\tau$ and the isomorphism 
$\mD_B \Itw \cong QDM_{\be}(M,\mV)$ 
sends $\Itw $ to a section $v(Q,\htau,\hbar)$ 
such that $v(0,\htau,\hbar)=1$. 

(2) Via the isomorphism in (1), $\mEh:=\mD_B \Itw$ becomes 
a big AQDM having a canonical frame $\Phi$ 
induced from the standard trivialization of $QDM_\be(M,\mV)$. 
The frame $\Phi_0\colon \mE_{\bb,0}\otimes K\{\hbar\} \to \mEh_{\bb,0}$ 
induced by $\Phi$ satisfies $\Phi_0[\Itw]=[\Itw]$, 
where $\bb=\{\tau=0\}\in B$.  
The co-ordinate $\htau=\htau(Q,\tau)$ are obtained as  
flat co-ordinates associated with $\Phi$ and 
$e_0=[\Itw] \in \mE_{0,\bb}$. 
\end{corollary} 
\begin{proof}
(1) is clear from the preceding discussion. 
(2) follows from that the classes $[\Itw]$ of $\Itw$ 
in both $\mEh_{0,\bb}$ and $\mE_{0,\bb}$ correspond 
to the same element 
$v(0,\htau,\hbar)=1=v(0,\htau,0)$ on the QDM side.  
\end{proof} 

\subsection{Cones of big AQDMs}
Here we study a relationship between cones and general big AQDMs.   
We use the notation in Section \ref{sec:AQDM}. 
Let $\mEh$ be a big AQDM on the base $B$ with a base point $\bb$. 
Let $t^0,\dots,t^s$ be a co-ordinate system on $B$. 
Here we take the infinite dimensional space $\mH$ to be 
\[
\mH = \mE_{\bb,0} \otimes K\{\hbar,\hbar^{-1}\}\!\}[\![Q]\!]
\]
\begin{definition}
Let $\Phi$ be a frame of $\mEh$ and 
$L=L(Q,t,\hbar)$ be the fundamental solution 
(Proposition \ref{prop:Lexists}) associated 
with an $\mO$-valued point $\bb$ on $B$ 
and $\Phi$. 
A {\it cone $\mL\subset \mH$ of the big AQDM $\mEh$} is defined to be 
the following subset of $\mH$:
\[
\mL := \bigcup_{\bb'\in B} 
\hbar \bL_{t(\bb')}  \quad 
\text{where} \quad 
\bL_{t(\bb')} := L(Q,t(\bb'),-\hbar)^{-1} 
(\mE_{\bb,0}\otimes K\{\hbar\}[\![Q]\!])
\subset \mH.
\]
Here, $\bb'$ ranges over the set of $\mO$-valued points in $B$. 
This depends on the choice of 
a base point $\bb$ and a frame $\Phi$. 
\end{definition} 

\begin{proposition}
\label{prop:symplectictransf_as_gaugetransf} 
Let $\mEh$ be a big AQDM 
endowed with a frame 
$\Phi\colon \mE_{\bb,0}\otimes\mO_B^\hbar\rightarrow \mEh$. 
Let $t^i$, $0\le i\le s$ be co-ordinates centered at $\bb$. 
Assume that the connection operators 
$A_a(Q,t,\hbar), 
\Omega_i(Q,t,\hbar)\in \End(\mE_{\bb,0})\otimes \mO_B^\hbar$ 
associated with $\Phi$ satisfy 
\[
A_a(0,t,\hbar) = p_a, \quad 
\Omega_i(0,t,\hbar) = \wp_i. 
\] 
Let  $\bb'$ be an $\mO$-valued point 
and $\Phi'\colon \mE_{\bb',0}\otimes \mO_B^\hbar\rightarrow \mEh$
be a frame such that for $G:=\Phi^{-1}\circ \Phi'\in 
\Hom(\mE_{\bb',0},\mE_{\bb,0})\otimes \mO_B^\hbar$,  
$G_0:=G|_{Q=0}$ is independent of $t$. 
Then the cone $\mL'$ associated with $\bb'$ and $\Phi'$ 
is related to the cone $\mL$ associated with $\bb$ and $\Phi$  
as 
\[
\mL'=G_{0}(-\hbar)^{-1}
\exp\left(-\frac{1}{\hbar}\sum_{i=0}^s c^i \wp_i\right)\mL, 
\] 
where $c^i=t^i(\bb')|_{Q=0}$.   
\end{proposition} 
\begin{proof}
Let $L(Q,t,\hbar)$ be the fundamental solution 
associated with $\bb$ and $\Phi$. 
By the assumption, we have 
$L(0,t,\hbar)=\exp(\sum_{i=0}^s t^i \wp_i/\hbar)$. 
We claim that the fundamental solution $L'$ 
associated with $\Phi'$ and $\bb'$ is given by 
(\emph{c.f.} (\ref{eq:gaugetransf_fundsol}))  
\[
L'(Q,t,\hbar) = G(Q,t,\hbar)^{-1} L(Q,t,\hbar) 
\exp(-\sum_{i=0}^s c^i \wp_i/\hbar)G_0(\hbar)
\]
This $L'$ satisfies the initial condition 
$L'(Q,t(\bb'),\hbar)|_{Q=0} = \id$. 
Let $A_a'$, $\Omega'_i$ be the connection operators 
associated with $\Phi'$. 
Then by (\ref{eq:gaugetransf_conn}), 
$A_a'|_{Q=0} =G_0^{-1} A_a|_{Q=0} G_0 = G_0^{-1} p_a G_0$ 
is independent of $t$. This implies 
$p_a = A'_a|_{Q=t=0} = G_0^{-1} p_a G_0$, \emph{i.e.} 
$G_0$ commutes with $p_a$. 
By this and (\ref{eq:gaugetransf_conn}), 
it easily follows that 
$(\hbar \partial/\partial t^i + \Omega_i')L'=0$ 
and $(\hbar (Q^a \partial/\partial Q^a) + A_a') L' + L'p_a=0$. 
Now the conclusion follows from the definition of $\mL$ and $\mL'$. 
\end{proof}

The symplectic transformation $\triangle$ 
in Theorem \ref{thm:coates-givental-cone} 
can be interpreted as the combination of 
a gauge transformation $G_0$ at $Q=0$ and 
a shift of the $K$-valued origin $\bb|_{Q=0}$. 
Note that the symplectic transformation  
does not change the $\mD_B$-module itself. 
By Theorem \ref{thm:coates-givental-cone} and Proposition 
\ref{prop:symplectictransf_as_gaugetransf}, 
we see that the twist of the quantum cohomology 
corresponds to the shift of the origin by 
$-\sum_{l\ge 1} s_{l-1}\ch_l(\mV)$ and 
the gauge transformation at $Q=0$: 
\[
G_0(\hbar)=\bc(\mV)^{1/2}\exp
\left(\sum_{m\ge 1,l\ge 0} s_{2m-1+l}
\frac{B_{2m}}{(2m)!}\ch_l(\mV)\hbar^{2m-1}\right).
\]


\subsection{Main Theorem} 
\label{subsec:main} 

We consider the following condition on 
a smooth projective toric variety $X$. 
\begin{condition}
\label{cond:ambienttoric}
There exists a smooth projective toric variety $Y$ 
such that $c_1(Y)$ is nef and 
$X$ is a complete intersection of {\it nef} 
toric divisors $D_1,\dots, D_k$ in $Y$. 
\end{condition}

We give a sufficient condition 
for Condition \ref{cond:ambienttoric} to hold. 

\begin{proposition}
\label{prop:embedding}
Assume that there exist nef toric divisors $D'_1,\dots, D'_l$ 
in $X$ and integers $n_1,\dots,n_l\ge 0$  
such that $c_1(X)+\sum_{i=1}^l n_i [D'_i]$ is nef. 
Then $X$ satisfies Condition \ref{cond:ambienttoric}.  
\end{proposition}
\begin{proof}
Recall the construction of a toric variety $X$ in 
Section \ref{subsec:toric}. 
We use the same notation for $X$ as there. 
The first Chern class of $X$ is given 
by the sum $c_1(X)=\sum_{i=1}^N u_i$ of toric divisor classes. 
Without loss of generality, 
we can assume $u_i=[D'_i]$ for $1\le i\le l$.  
By the assumption, $c_1(X)+\sum_{i=1}^l n_iu_i$ 
are in the closure $\ov{C}_X$ of the K\"{a}hler cone 
$C_X$ (\ref{eq:Kaehlercone}) 
and so are $u_i$, $1\le i\le l$. 
We define a toric variety $Y$ by the following data 
(see Section \ref{subsec:toric}):  
\begin{itemize}
\item[(1)] the same algebraic torus $\T_\C\cong (\C^*)^r$. 
\item[(2)] $(N+\sum_{i=1}^l n_i)$-tuple of integral vectors 
$(u_1,\dots,u_1,\dots,u_l,\dots,u_l,u_{l+1},\dots,u_N)$ 
in $\Hom(\T_\C,\C^*)$ 
where each $u_i$ is repeated $n_i+1$ times for $1\le i\le l$.  
\item[(3)] the same vector $\eta\in \Hom(\T_\C,\C^*)\otimes \R$. 
\end{itemize}
It is easy to check that these data satisfy the conditions 
(A)--(C) in Section \ref{subsec:toric}, 
so $Y$ is a smooth projective toric variety. 
We have $H^2(Y,\Z)\cong H^2(X,\Z)$ and $C_Y \cong C_X$.  
By the construction, $c_1(Y)$ is nef  
and $X$ is a complete intersections of $(\sum_{i=1}^l n_i)$ 
nef toric divisors in $Y$. 
\end{proof} 

By using the same method as above, we can show that 
every smooth projective toric variety $X$ 
can be realized as a complete intersection 
of {\it not necessarily nef} toric divisors in a smooth 
projective Fano toric variety $Y$.

\begin{theorem}[Main theorem] 
\label{thm:main}
Let $X$ be a smooth projective toric variety satisfying Condition 
\ref{cond:ambienttoric} 
and $\mV$ be a sum of line bundles on $X$.  
By the generalized mirror transformation 
in Definition \ref{def:GMT}, 
the equivariant Floer cohomology $\mEh = FH^*_{T^2}(L_{X/\mV})$ 
together with the choice of a frame $\Phi_0=\Loc_0^{-1}$ of $\mEh_0$ 
and $e_0=[\Delta]\in \mE_0$ reconstructs 
the big QDM $QDM_\be(X,\mV)$ twisted by $\mV$ and the 
$S^1$-equivariant Euler class $\be$.  
\end{theorem}

\begin{remark*}
In the previous preprint (version 4) \cite{iritani-genmir-v4}, 
the main theorem was stated for an arbitrary smooth projective 
toric variety. The proof there contained technical mistakes
when $X$ is a complete intersection 
of {\it not nef} toric divisors in a Fano toric variety $Y$. 
Here, we impose Condition \ref{cond:ambienttoric} on $X$ instead.  
The author, however, hopes that this condition 
will be removed in the future. 
\end{remark*} 
The rest of Section \ref{subsec:main} is 
devoted to the proof of Theorem \ref{thm:main}. 
Let $\mV$ be a sum $\mV_1\oplus \cdots \oplus \mV_l$ 
of line bundles over $X$ 
and $v_i=c_1(\mV_i)\in H^2(X)$. 
See Section \ref{subsec:toric}
for the notation on toric varieties and 
Section \ref{sec:QDM} for that on QDMs. 

\subsubsection{When $X$ is Fano} 
\label{subsubsec:FanoX}
By Givental's mirror theorem \cite{givental-mirrorthm-toric},
we know that for a Fano toric variety $X$, 
the small $J$-function $J_X(Q,\hbar)$ coincides with 
the $I$-function $I_X(Q,\hbar)$ of $X$ itself, 
where $I_X$ is (\ref{eq:I-funct}) with $\mV=0$. 
Let $\Itw_{\mV}(Q,\tau,\hbar)$ denote the 
hypergeometric modification (\ref{eq:hypergeometric_modif}) 
of the big $J$-function of $X$ with respect to $\mV$. 
Then it is immediate to see that 
\[
\Itw_{\mV} (Q,0,\hbar) = I_{X,\mV}(Q,\hbar)  
\]
where the right-hand side is the $I$-function in (\ref{eq:I-funct}). 
Let $B$ be the formal neighborhood of zero in $H(X)\otimes \mO$ 
and take the ground ring $K$ to be $\C(\!(\lambda^{-1})\!)$. 
The right-hand side generates $\mEh:=FH_{T^2}(L_{X/\mV})$ 
over the Heisenberg algebra $\mD$ 
by Proposition \ref{prop:Jgenerates_AQDM}. 
On the other hand, the $\mD_B$-module $\mFh$ generated 
by $\Itw_{\mV}$ is isomorphic to $QDM_\be(X,\mV)$ 
by Corollary \ref{cor:modif_generates}. 
By the same corollary,  
the class $[\Itw_{\mV}]$ in $\mF_{0,\bb}$ corresponds to $1\in H^*(X)$, 
so generates $\mF_{0,\bb}$ as a $K[p_1,\dots,p_r]$-module. 
Then by Proposition \ref{prop:Jgenerates_AQDM}, $\mFh$ 
is generated by $\Itw_\mV$ as a  
$\mO^\hbar_B\langle \Qp_1,\dots,\Qp_r\rangle$-module 
and the restriction of $\mFh$ to $\tau=0$ is
generated by $\Itw_\mV(Q,0)$ as a $\mD$-module. 
Therefore, $\mFh|_{\tau=0} \cong \mEh$ as $\mD$-modules. 
By the uniqueness of the reconstruction 
in Theorem \ref{thm:reconstruction}, $\mFh$ must 
be isomorphic to the big AQDM $\hmEh$ reconstructed from $\mEh$. 
Finally we check that the frame $\Loc_0^{-1}$ of $\mEh_0$ and 
$e_0=[\Delta]\in \mE_0$ correspond to the standard trivialization of 
the QDM at $Q=0$ and the unit $1\in H^*(X)$. 
The class $[\Delta]\in \mE_0$ corresponds to 
$[\Itw_{\mV}]\in \mF_{\bb,0}$, so corresponds  
to $1\in H^*(X)$ by Corollary \ref{cor:modif_generates}. 
By the same corollary, the frame $\Phi_0$ of $\mEh$ 
induced from the standard trivialization of $QDM_\be(X,\mV)$ 
satisfies $\Phi_0[\Delta]=[\Delta]$. 
Since $\Phi_0$ is a homomorphism of $K\{\hbar\}[p_1,\dots,p_r]$-modules 
(\ref{eq:frame_zero}), this completely determines $\Phi_0$. 
By Theorem \ref{thm:freenessofefc}, $\Loc_0^{-1}$ satisfies 
the same properties as $\Phi_0$, thus $\Phi_0=\Loc_0^{-1}$.

\subsubsection{From $\mV=0$ to $\mV\neq 0$} 
We will show that Theorem \ref{thm:main} holds for a pair $(X,\mV)$ 
when it holds for $X$ itself. 
We will assume that $\mV$ is a line bundle and set $v=c_1(\mV)$. 
The proof in the case where $\mV$ is a sum of line bundles 
requires only notational changes. 
By assumption, 
$FH_{S^1}^*(L_X)$ is isomorphic to the restriction of 
$QDM(X)$ to some $c(Q)\in H^*(X)\otimes \C[\![Q]\!]$ 
such that $c(0)=0$. 
By an argument similar to that showing (\ref{eq:fromItoJ_diff}), 
we have 
\begin{equation}
\label{eq:IX_JX}
I_X(Q,\hbar) = f(Q,\hbar Q\parfrac{}{Q}+p, \hbar) J_X(Q,c(Q),\hbar) 
\end{equation} 
for some $f(Q,\Qp,\hbar)\in \C[\hbar][\![Q]\!]\langle \Qp \rangle$ 
such that $f(0,\Qp,\hbar)=1$. 
Set 
\[
M_d(x) = 
\frac{\prod_{\nu=-\infty}^{\pair{v}{d}}(x+\lambda+\nu\hbar)}
{\prod_{\nu=-\infty}^0 (x+\lambda+\nu\hbar)} 
\in \lambda^{\pair{v}{d}} \C[\hbar,x][\![\lambda^{-1}]\!].  
\]
Put $v=\sum_{a=1}^r v^a p_a$ and 
$\Qla^a:=Q^a \lambda^{v^a}$.  
We define the modification of a general cohomology-valued 
$Q$-series $\sum_d Q^d H_d(t,\hbar) \in H^*(X)\otimes 
\C(\!(\hbar^{-1})\!)[\![Q,t]\!]$ to be 
\[
\modif(H):= \sum_{d} Q^d  M_d(v) H_d(t,\hbar)  
\in H^*(X) \otimes \C(\!(\hbar^{-1})\!)
[\![\lambda^{-1},\Qla,t]\!].  
\]
Note that $\C(\!(\hbar^{-1})\!)[\![\lambda^{-1},\Qla,t]\!]$ 
is a subring of $\C(\!(\hbar^{-1})\!)(\!(\lambda^{-1})\!)[\![Q,t]\!]$. 
For example, we have 
$I_{X,\mV}=\modif(I_X)$ and $\Itw_{\mV}=\modif(J_X)$. 
Set $\bpartial_a = Q^a \partial/\partial Q^a$ and 
$\partial_i=\partial/\partial t^i$. 
The modification of a differential operator 
$g=\sum_{d}Q^d g_d(t,\hbar\bpartial+p,\hbar \partial,\hbar)$ 
in $\C[\hbar,\hbar\bpartial+p][\![t]\!][\hbar\partial][\![Q]\!]$ 
is defined to be  
\[
\modif(g):= 
\sum_{d} Q^d M_d(\hbar \bpartial_v+v) 
g_d(t,\hbar\bpartial, \hbar \partial,\hbar)  
\in \C[\hbar,\hbar\bpartial+p]
[\![\lambda^{-1},t]\!][\hbar\partial][\![Q_\lambda]\!] 
\]
where $\bpartial_v= \sum_{a=1}^r v^a \bpartial_a$.  
\begin{lemma}
\label{lem:modif_hom} 
Let $g_1,g_2$ be differential operators in 
$\C[\hbar,\hbar\bpartial+p][\![t]\!][\hbar\partial][\![Q]\!]$ 
and $H$ be a cohomology-valued power series in 
$H^*(X)\otimes \C(\!(\hbar^{-1})\!)[\![Q,t]\!]$. 
We have  
\[
\modif(g_1g_2)=\modif(g_1)\modif(g_2), \quad 
\modif(g_1H)= \modif(g_1)\modif(H). 
\]
\end{lemma} 
\begin{proof}

We expand $g_1=\sum_{d} Q^d g_{1,d}$ and 
$g_2=\sum_d Q^d g_{2,d}$, where 
$g_{i,d}\in \C[\hbar,\hbar\bpartial+p]
[\![t]\!][\hbar\partial][\![Q]\!]$. 
By using $M_{d_1+d_2}(x) = M_{d_1}(\hbar \pair{v}{d_2} + x)M_{d_2}(x)$, 
\begin{align*}
\modif(g_1g_2) 
& = \sum_{d} Q^d M_d(\hbar \bpartial_v+v) 
\sum_{d=d_1+d_2}  g_{1,d_1} g_{2,d_2} \\
& = \sum_{d_1,d_2} Q^{d_1} Q^{d_2} 
M_{d_1}(\hbar \pair{d_2}{v}+ \hbar\bpartial_v+v) 
g_{1,d_1} M_{d_2}(\hbar\bpartial_v+v) g_{2,d_2} \\
& = \modif(g_1) \modif(g_2). 
\end{align*}
This proves the first equation. 
The proof of the second equation is similar. 
\end{proof} 

Put $c(Q)=\sum_{i=0}^s c^i(Q)p_i$. 
By using (\ref{eq:IX_JX})
and Lemma \ref{lem:modif_hom}, we have 
\begin{align}
\nonumber 
I_{X,\mV}(Q,\hbar) &= \modif(I_X(Q,\hbar)) \\ 
\nonumber
&= \modif(f) \modif(J_X(Q,c(Q),\hbar)) \\
\label{eq:IXV_modif_IV}
&= \modif(f) e^{\sum_{i=0}^s \modif(c^i(Q)\partial_i)} 
\Itw_\mV(Q,\tau,\hbar)\Bigr|_{\tau=0}  
\end{align} 
where we used 
$J_X(Q,c(Q),\hbar)= \exp(\sum_{i=0}^s
c^i(Q)\partial_i)J_X(Q,\tau,\hbar)|_{\tau=0}$
in the third line. 
By Corollary \ref{cor:modif_generates}, it follows that 
the $\mD_B$-module generated by $\Itw_\mV$ 
has the following free basis over 
$\mO^\hbar_B=\C[\hbar](\!(\lambda^{-1})\!)[\![Q,t]\!]$: 
\begin{equation}
\label{eq:basis_DItw} 
\hbar \partial_0 \Itw_\mV=\Itw_\mV, 
\quad \hbar \partial_1 \Itw_\mV, \quad 
\dots, \quad \hbar\partial_s\Itw_\mV. 
\end{equation} 
We consider a submodule $\mFh$ generated by these elements 
over $\qmO^\hbar_B:=\C[\hbar][\![\lambda^{-1},\Qla,t]\!]$. 
This is invariant under the action of $\hbar\bpartial_a+p_a$ 
and $\hbar\partial_i$, so has the structure of a 
$\qmD_B:=\qmO^\hbar_B\langle \hbar\bpartial+p,
\hbar\partial_i\rangle$-module.  
In fact, the connection matrices $A_a=(A^j_{ak})$, 
$\Omega_i=(\Omega^j_{ik})$ 
corresponding to $\hbar\bpartial_a+p_a$ 
and $\hbar\partial_i$ are given by 
\[
(\hbar \bpartial_a + p_a) \hbar \partial_k \Itw_\mV = 
\sum_{j=0}^s (\hbar \partial_j \Itw_\mV) A_{ak}^j, \quad 
\hbar \partial_i \hbar (\partial_k \Itw_\mV) 
= \sum_{j=0}^s (\hbar \partial_j \Itw_\mV) \Omega_{ik}^j  
\] 
It is obvious that the matrix entries $A^j_{ak}$, $\Omega^j_{ik}$ 
belong to $\mO_B^\hbar$, but from $\Itw_\mV\in 
H^*(X)\otimes \C(\!(\hbar^{-1})\!)[\![\lambda^{-1},\Qla,t]\!]$ 
and $\Itw_\mV|_{\Qla=0} = \exp(\tau/\hbar)$, 
it follows that they belong also to $\qmO_B^\hbar$. 
Put $\smnabla^\hbar_a = \hbar \bpartial_a+A_a$ and 
$\nabla^\hbar_i =\hbar \partial_i+\Omega_i$. In 
the basis (\ref{eq:basis_DItw}),  
the element $e^{\sum_{i=0}^s \modif(c^i(Q)\partial_i)} 
\Itw_\mV(Q,\tau,\hbar)$ is represented by 
\begin{equation}
\label{eq:exp_diff_1}
\exp(g/\hbar)   
[1,0,\dots,0]^{\rm T},\quad  
\text{where } 
g=\sum_{i=0}^s \sum_d Q^d M_d(\smnabla^\hbar_v) 
c^i_d \nabla^\hbar_i.  
\end{equation} 
Here $c^i(Q)=\sum_{d} Q^d c_d^i$ and 
$[1,0,\dots,0]^{\rm T}$ represents $\Itw_\mV$.  
\begin{lemma}
\label{lem:exp_diff_fac} 
Let $g$ be a differential operator in 
$\C[\hbar,\smnabla^\hbar][\![\lambda^{-1}, t]\!]
[\nabla^\hbar][\![\Qla]\!]$ 
such that $g|_{\Qla=0}=0$. 
Then there exist a differential operator 
$h$ in the same ring and 
$\frc^i(Q,\tau) \in \qmO_B :=\C[\![\lambda^{-1},\Qla,t]\!]$ 
such that $h|_{\Qla=0}=0$, $\frc^i(0,\tau)=0$, 
\[
\exp(g/\hbar) = \exp\left(\frac{1}{\hbar} \sum_{i=0}^s \frc^i(Q,\tau) \nabla^\hbar_i\right) \exp(h/\hbar) 
\]
and that the action of $h/\hbar$ preserves the module 
$\mFh\cong \C^{s+1}\otimes \qmO_B^\hbar$. 
\end{lemma} 
\begin{proof} 
We assume by induction that 
there exist $h_n \in \C[\hbar,\smnabla^\hbar][\![\lambda^{-1},t]\!]
[\nabla^\hbar][\![\Qla]\!]$ and 
$\frc^i_{n}(Q,\tau)\in \qmO_B$ such that 
\begin{equation}
\label{eq:factor_n}
\exp(g/\hbar) = 
\exp\left (\frac{1}{\hbar} \sum_{i=0}^s \frc^i_n(Q,\tau) 
\nabla^\hbar_i \right)   \exp(h_n/\hbar)  
\end{equation} 
and that 
\begin{equation} 
\label{eq:vanish_upton}
\frac{h_n}{\hbar} \mFh 
\subset \mFh + \frac{1}{\hbar} \sum_{|d|>n} \Qla^d \mFh. 
\end{equation} 
Here, we identify $\mFh$ with 
$\C^{s+1}\otimes \qmO_B^\hbar$ 
and $|d|=\sum_{a=1}^r \pair{p_a}{d}$. 
When $n=0$, we can take $\frc^i_0=0$. 
We have for $h_n=h_n(Q,t,\smnabla^\hbar,\nabla^\hbar,\hbar)$, 
\[
\Res_{\hbar=0} d\hbar ((h_n/\hbar) v) 
= h_n(Q,t,A_a|_{\hbar=0}, \Omega_i|_{\hbar=0},0) v(Q,t,0), 
\]
where any $v=v(Q,t,\hbar) \in 
\C^{s+1}\otimes \qmO_B^\hbar$. 
By the induction hypothesis, 
\[
[w_1(Q,t),\dots, w_s(Q,t)]^{\rm T} := 
h_n(Q,t,A_a|_{\hbar=0},\Omega_i|_{\hbar=0},0) 
[1,0,\dots,0]^{\rm T}
\]
is in $\sum_{|d|>n} \Qla^d 
\C^{s+1}\otimes \qmO_B$. 
Since $\Omega_i[1,0,\dots,0]^{\rm T}
=[0,\dots, 0, 1, 0,\dots,0]^{\rm T}$,  
where $1$ is in the $i$-th position, 
we know that the application of 
\[
h_n(Q,t,A_a|_{\hbar=0},\Omega_i|_{\hbar=0},0) 
- \sum_{i=0}^s w^i(Q,t) \Omega_i|_{\hbar=0}
\]
to the vector $[1,0,\dots,0]^{\rm T}$ is zero. 
Because $A_a|_{\hbar=0}$ and $\Omega_i|_{\hbar=0}$ commute 
with each other 
(by the flatness of $\smnabla^\hbar$ and $\nabla^\hbar$) and 
because $\C^{s+1}\otimes \qmO_B$ is generated 
by $[1,0,\dots,0]^{\rm T}$ as a $\C[A_a|_{\hbar=0}][\![\lambda^{-1},t]\!]
[\Omega_i|_{\hbar=0}][\![\Qla]\!]$-module, we know that 
the above endomorphism is identically zero. 
Put $\frc_{n+1}^i(Q,t) = \frc^i_n(Q,t) + w^i(Q,t)$.  
We define $\nabla_{\frc_n}:= 
(1/\hbar) \nabla^\hbar_{\frc_n} := 
(1/\hbar)\sum_{i=0}^s  \frc_n^i \nabla^\hbar_i$. 
We define $\nabla^\hbar_w$, $\nabla_w$ and 
$\nabla_{\frc_{n+1}}$ similarly. 
By Baker-Campbell-Hausdorff formula, we have 
\begin{align*}
\exp(-\nabla_{\frc_{n+1}})\exp(\frac{g}{\hbar})  
&= \exp(-\nabla_{\frc_{n}+w} + \frac{g}{\hbar}    
-\frac{1}{2}[\nabla_{\frc_n+w}, \frac{g}{\hbar}]  
+\cdots)    \\
& = \exp((-\nabla_{\frc_n} + \frac{g}{\hbar} 
-\frac{1}{2}[\nabla_{\frc_n},\frac{g}{\hbar}] 
+ \cdots) -  
(\nabla_w +\frac{1}{2}[\nabla_w,\frac{g}{\hbar}]  
+\cdots))  
\\
&=\exp((h_{n} - \nabla^\hbar_w - R_n)/\hbar). 
\end{align*} 
In the third line, we used the Baker-Campbell-Hausdorff formula 
for $e^{h_n/\hbar} = e^{-\nabla_{\frc_n}}e^{g/\hbar}$ 
and $R_n/\hbar =\frac{1}{2}[\nabla_w,\frac{g}{\hbar}]  
+\frac{1}{12}[\frac{g}{\hbar}, [\nabla_w,\frac{g}{\hbar}]]+\cdots$ 
is the remainder. 
Because the space 
$\frac{1}{\hbar} \C[\hbar,\smnabla^\hbar][\![\lambda^{-1},t]\!]
[\nabla^\hbar][\![\Qla]\!]$ is closed under the commutator 
$[\cdot,\cdot]$, $R_n/\hbar$ belongs to this space.  
Since $w^i(Q,t)\in \sum_{|d|>n} \Qla^d \qmO_B$ 
and $g|_{\Qla=0}=0$, it follows that 
$R_n\in \sum_{|d|>n+1} \Qla^d
\C[\hbar, \smnabla^\hbar][\![\lambda^{-1},t]\!]
[\nabla^\hbar][\![\Qla]\!]$. 
Now $\frc_{n+1}^i$ and 
$h_{n+1} := h_n - \nabla^\hbar_w - R_n$ satisfy 
(\ref{eq:factor_n}) and (\ref{eq:vanish_upton}). 
By construction, $\frc_n$ and $h_n$ converges in 
$\Qla$-adic topology and the conclusion follows. 
\end{proof}
By (\ref{eq:exp_diff_1}) and 
Lemma \ref{lem:exp_diff_fac}, there exist 
$\frc^i(Q,\tau)\in \qmO_B$ 
and $v^i(Q,\tau)\in \qmO_B$  
such that $\frc^i(0,\tau)=0$, $v^i(0,\tau) =\delta_{i0}$, 
\begin{align*}
e^{\sum_{i=0}^s \modif(c^i(Q)\partial_i)} 
\Itw_\mV(Q,\tau,\hbar) 
&= e^{\sum_{i=0}^s \frc^i(Q,\tau) \partial_i} 
\left(\sum_{i=0}^s v^i(Q,\tau) 
\hbar\partial_i \Itw_\mV(Q,\tau,\hbar)
\right)  \\
&= \sum_{i=0}^s v^i(Q,\tau^*) 
(\hbar\partial_i \Itw_\mV)(Q,\tau^*,\hbar) 
\end{align*}  
where $\tau^* = \tau^*(Q,\tau) = 
e^{\sum_{i=0}^s \frc^i(Q,\tau)\partial_i}\tau$. 
From what we have discussed, the 
$\qmD_B$-module generated by $\Itw_\mV$ 
is considered to be a big AQDM over the ground ring 
$K=\C[\![\lambda^{-1}]\!]$ (with $Q$ replaced with $\Qla$) 
and is isomorphic to a $\qmD_B$-submodule of $QDM_\be(X,\mV)$. 
Under this isomorphism, $\Itw_\mV(\tau^*(Q,0),Q,\hbar)$ 
corresponds to a section $v'(Q,\hbar)$ of the QDM on a locus 
$\tau'=\tau'(Q)\in H^*(X)\otimes \C[\![\lambda^{-1},\Qla]\!]$ 
such that $v'(0,\hbar)=1$ and 
$\tau'(0)=0$. Thus by Proposition \ref{prop:Jgenerates_AQDM}, 
there exists a differential operator 
$f'(Q,\hbar\bpartial+p,\hbar) 
\in \C[\hbar][\![\lambda^{-1},\Qla]\!][\hbar\bpartial+p]$ 
such that $f'(0,\hbar\bpartial+p,\hbar)=1$ 
and 
\[
\sum_{i=0}^s v^i(Q,\tau^*) (\hbar\partial_i \Itw_\mV)
(Q,\tau^*,\hbar)\Bigl|_{\tau=0} 
=f'(Q,\hbar\bpartial+p,\hbar)
\Itw_\mV(Q,\tau^*(Q,0),\hbar). 
\] 
Then by (\ref{eq:IXV_modif_IV}), for 
$f''= \modif(f) f'$, 
\begin{equation}
\label{eq:mod_mod}
I_{X,\mV}(Q,\hbar) = 
f''(Q,\hbar\bpartial+p,\hbar) 
\Itw_\mV(Q,\tau^*(Q,0),\hbar) 
\end{equation} 
The differential operator $f''$ 
belongs to 
$\hmD=\C[\hbar,\hbar\bpartial+p][\![\lambda^{-1},\Qla]\!]$  
which is considered to be a certain completion of 
$\mD=\C[\hbar][\![\lambda^{-1},\Qla]\!]\langle \hbar\bpartial+p\rangle$. 
It is easy to show that a small AQDM 
(over $K=\C[\![\lambda^{-1}]\!]$) has the structure 
of a $\hmD$-module. 
Thus, by Corollary \ref{cor:modif_generates}, 
the above formula means that $I_{X,\mV}$ corresponds 
to a section $v''$ of $QDM_\be(M,\mV)$ on the locus 
$\tau'(Q)$ such that $v''|_{Q=0}=1$. 
Thus by Proposition \ref{prop:Jgenerates_AQDM}, the 
$\mD$-module generated by $I_{X,\mV}(Q,\hbar)$ is isomorphic 
to the restriction of $QDM_\be(X,\mV)$ to 
the subspace $\{\tau'=\tau'(Q)\}$. 
Now the conclusion follows from 
that $FH_{T^2}(L_{X/\mV})$ is generated by $I_{X/\mV}$ and 
the uniqueness of the reconstruction 
in Theorem \ref{thm:reconstruction}. 
A routine argument shows that $\Loc_0^{-1}$ and 
$[\Delta]$ corresponds to the standard trivialization of 
the QDM and the unit $1\in H^*(X)$.

\subsubsection{When $\mV=0$} 
\label{subsubsec:Vzero} 
By Condition \ref{cond:ambienttoric}, 
we embed $X$ into $Y$ as a complete intersection 
of nef toric divisors $D_1,\dots, D_k$. 
Without loss of generality, 
we can assume that $X$ and $Y$ are given as 
the GIT quotients $\C^N/\!/\T_\C$ and $\C^{N+k}/\!/\T_\C$
respectively and that the embedding $i\colon X\hookrightarrow Y$ 
is induced from the inclusion of the co-ordinate subspace 
$\C^N \subset \C^{N+k}$ consisting of first $N$ factors. 
Let $z_1,\dots, z_{N+k}$ be the standard co-ordinates on $\C^{N+k}$  
and $u_i\in H^2(Y)$ be the class of the toric divisor $\{z_i=0\}$. 
Then we have $u_{N+i}=[D_i]$ for $1\le i\le k$. 
In particular, $u_{N+1},\dots,u_{N+k}$ are nef. 
Because the cohomology rings of $X$ and $Y$ are generated by 
toric divisor classes, 
the restriction map $i^*\colon H^*(X)\rightarrow H^*(Y)$ 
is surjective. 
Thus, the push-forward 
$i_* \colon H_*(X)\rightarrow H_*(Y)$ is injective. 
Let $p_1,\dots,p_{r+l}$ be a nef integral basis of $H^2(Y)$ 
such that $i^*p_{r+1}=\cdots=i^*p_{r+l}=0$ and 
that $i^*p_1,\dots, i^*p_r$ form a nef integral basis of $H^2(X)$. 
Introduce Novikov variables $Q^1,\dots,Q^{r+l}$ dual to 
$p_1,\dots,p_{r+l}$. We use 
\begin{align*} 
Q^d &= (Q^1)^{\pair{p_1}{d}} \cdots (Q^{r+l})^{\pair{p_{r+l}}{d}} 
 &\text{for } d\in H_2(Y), \\
Q^d &= Q^{i_*d} = (Q^1)^{\pair{i^*p_1}{d}} \cdots 
(Q^r)^{\pair{i^*p_r}{d}} 
 &\text{for } d\in H_2(X). 
\end{align*}
Let $\mW$ denote the sum of nef line bundles on $Y$: 
\[
\mW = \mO(D_1) \oplus \cdots \oplus \mO(D_k). 
\]

\begin{lemma}
Let $J_{Y,\mW}(Q,\tau,\hbar)$ 
be the $J$-function of 
the $(\mW, \be)$-twisted quantum cohomology of $Y$. 
Let $J_X(Q,\tau,\hbar)$ be the $J$-function of the 
quantum cohomology of $X$. 
Then we have 
\label{lem:Lefschetz}
\[ 
\lim_{\lambda \to 0}
\be(\mW)\cup 
J_{Y,\mW}(Q,\tau,\hbar)
 =i_*J_X(Q,i^*\tau,\hbar). 
\] 
\end{lemma}
\begin{proof}
This follows from Theorem \ref{thm:kkp}.  
\end{proof} 

Let $J_Y(Q,\tau,\hbar)$ be the $J$-function 
of $Y$ and $\Itw_{\mW}(Q,\tau,\hbar)$ 
be the hypergeometric modification 
(\ref{eq:hypergeometric_modif}) of $J_Y$ with respect to $\mW$.  
Let $I_{Y,\mW}(Q,\hbar)$ be the $I$-function 
(\ref{eq:I-funct}) of $(Y,\mW)$. 
Since $c_1(Y)$ is nef, by Givental's mirror theorem 
\cite{givental-mirrorthm-toric}, 
there exist  $f_i(Q)\in \C[\![Q]\!]$, $1\le i\le r+l$ 
such that $f_i(0)=0$ and 
\[
J_Y(Q, c(Q), \hbar)= I_Y(Q,\hbar), \quad 
c(Q)=\sum_{a=1}^{r+l} f_a(Q) p_a. 
\] 
Since $D_1,\dots,D_k$ are nef, 
the modification factor with respect to $\mW$ 
is polynomial in $\lambda$. 
We repeat the same argument as in Section 
\ref{subsubsec:Vzero} 
by replacing the rings  
\[
\C[\hbar][\![\lambda^{-1},\Qla,t]\!], \quad  
\C(\!(\hbar^{-1})\!)[\![\lambda^{-1},\Qla,t]\!], \quad 
\C[\hbar,\hbar\bpartial+p]
[\![\lambda^{-1},t]\!][\hbar\partial][\![\Qla]\!]
\]
appearing there with 
\[
\C[\hbar,\lambda][\![Q,t]\!], \quad 
\C(\!(\hbar^{-1})\!)[\lambda][\![Q,t]\!] \quad 
\C[\hbar, \hbar\bpartial+p,\lambda]
[\![t]\!][\hbar\partial][\![Q]\!]. 
\]
Then we can show that 
the relation (\ref{eq:mod_mod}) holds 
when $I_{X,\mV}$ and $\Itw_\mV$ replaced with 
$I_{Y,\mW}$ and $\Itw_\mW$ 
and for some 
$f''\in \hmD:= \C[\hbar,\hbar\bpartial+p,\lambda][\![Q]\!]$
and $\tau^*(Q,0) \in \C[\lambda][\![Q]\!]$. 
This again shows (by Proposition \ref{prop:Jgenerates_AQDM}) 
that 
$\mD:=\C[\lambda][\![Q]\!]\langle \hbar\bpartial+p \rangle$-module 
generated by $I_{Y,\mW}(Q,\hbar)$ equals 
the $\mD$-module generated by $\Itw_{\mW}(Q,\tau^*(Q,0),\hbar)$ 
as a submodule of 
$H^*(Y)\otimes \C(\!(\hbar^{-1})\!)[\lambda][\![Q]\!]$.  
In particular, there exists a differential operator 
$f_1(Q,\hbar\bpartial+p,\hbar) \in \mD$ such that 
$f_1(0,\hbar\bpartial+p,\hbar)=1$ and 
\begin{equation}
\label{eq:Itw_IYW}
\Itw_{\mW}(Q,\tau^*(Q,0),\hbar)= f_1(Q,\hbar\bpartial + p,\hbar) 
I_{Y,\mW}(Q,\hbar).  
\end{equation}
A routine argument using 
Corollary \ref{cor:modif_generates} and Proposition 
\ref{prop:Jgenerates_AQDM} shows that 
the $\mD$-module generated by 
$\Itw_\mW(Q,\tau^*(Q,0),\hbar)$ 
is isomorphic to the restriction of the 
$QDM_\be(Y,\mW)$ to some 
$\htau(Q)\in H^*(Y)\otimes \C[\lambda][\![Q]\!]$. 
It is easy to see that 
$\Itw_\mW(Q,\tau^*(Q,0),\hbar)$ is a $J$-function 
of the small AQDM generated by 
$\Itw_\mW(Q,\tau^*(Q,0),\hbar)$ for some frame. 
By the relationships  
(\ref{eq:fromItoJ_diff}), (\ref{eq:fromJtohJ}) among $J$-functions, 
there exists a differential operator 
$f_2(Q,\hbar\bpartial+p,\hbar)\in \mD$ 
such that $f_2(0,\hbar\bpartial+p,\hbar)=1$ and 
\begin{align}
\nonumber 
J_{Y,\mW}(Q,\htau(Q),\hbar) &= 
f_2(Q,\hbar \bpartial+p,\hbar) 
\Itw_{\mW}(Q,\tau^*(Q,0),\hbar)  \\
\label{eq:JYW_IYW}
&= V(Q,\hbar \bpartial + p, \hbar) I_{Y,\mW} (Q,\hbar).   
\end{align} 
where in the second line we used (\ref{eq:Itw_IYW}) 
and put $V = f_2 f_1$.  
On the other hand, using 
\begin{align}
I_{Y,\mW}(Q,\hbar) = \sum_{d\in\Lambda_Y} Q^d 
\prod_{i=1}^N 
\frac{\prod_{\nu=-\infty}^{0} (u_i+\nu\hbar)}
     {\prod_{\nu=-\infty}^{\pair{u_i}{d}}(u_i+\nu\hbar)} 
\cup \prod_{j=N+1}^{N+k} \prod_{\nu=1}^{\pair{u_j}{d}} 
\frac{u_j+\lambda+\nu\hbar}{u_j+\nu\hbar}
\end{align}
and 
$i_* i^*(\alpha) = \alpha \cup \prod_{j=N+1}^{N+k} u_j$ 
for $\alpha\in H^*(Y)$, we can easily check that 
\begin{equation}
\label{eq:Itw_IX}
\lim_{\lambda \to 0}
\be(\mW)\cup I_{Y,\mW}(Q,\hbar)
=i_*I_{X}(Q,\hbar)
\end{equation} 
By Lemma \ref{lem:Lefschetz} and Equations (\ref{eq:JYW_IYW}), 
(\ref{eq:Itw_IX}), we have 
\[
i_*J_X(Q,\lim_{\lambda\to 0} i^*\htau(Q),\hbar) = 
\left(\lim_{\lambda\to 0} V(Q,\hbar \bpartial+p, \hbar)\right) 
i_* I_X(Q,\hbar).   
\]
Note that the limit $\lambda\to 0$ exists 
since everything is defined over $K=\C[\lambda]$. 
Set $V_0(Q,\Qp,\hbar)=\lim_{\lambda\to 0} V(Q,\Qp,\hbar)$ 
and $\htau_0(Q)=\lim_{\lambda \to 0}i^* \htau(Q)$. 
By the injectivity of $i_*$, 
\begin{equation}
\label{eq:JX_IX}
J_X(Q,\htau_0(Q),\hbar) = V_0(Q,\hbar \bpartial+i^*p, \hbar) 
I_X(Q,\hbar).  
\end{equation} 
Since $I_X$ and $J_X$ does not depend on $Q^{r+1},\dots,Q^{r+l}$, 
by putting $Q^{r+1}=\cdots=Q^{r+l}=0$ if necessary, 
we can assume $\htau_0(Q)$ and $V_0(Q,\Qp,\hbar)$ depend only on 
$Q^1,\dots,Q^r$ and $\Qp_1,\dots,\Qp_r$ and $\hbar$.   
Then (\ref{eq:JX_IX}) implies that $J_X(Q,\htau_0(Q),\hbar)$ 
belongs to the $\mD$-module $FH_{S^1}^*(L_X)$ 
generated by $I_X(Q,\hbar)$. 
Because $V_0|_{Q=0}=1$, by Proposition \ref{prop:Jgenerates_AQDM}, 
$J_X(Q,\htau_0(Q),\hbar)$ also generates $FH_{S^1}^*(L_X)$. 
Therefore, the restriction of $QDM(X)$ to $\tau=\htau_0(Q)$ 
(which is generated by $J_X(Q,\htau_0(Q),\hbar)$ by Proposition 
\ref{prop:Jgenerates_AQDM}) 
is isomorphic to $FH_{S^1}^*(L_X)$. 
Now the conclusion follows from the uniqueness of 
the reconstruction in Theorem \ref{thm:reconstruction}.




\section{Example}
Let $M_N^k$ be a degree $k$ hypersurface of $\Proj^{N-1}$. 
This variety is Fano if $k<N$, Calabi-Yau if $k=N$ and 
of general type if $k>N$. 
We compute the big quantum cohomology of a 
general type hypersurface $M_8^9$ . 
More precisely, we will compute the Euler-twisted 
quantum cohomology $QH_{\Euler}^*(\Proj^7,\mO(9))$. 
We follow the procedure of the generalized mirror transformation 
in Definition \ref{def:GMT}. 

Since $\mO(9)$ is ample, the $T^2$-equivariant 
Floer cohomology $FH_{T^2}^*(L_{\Proj^7/\mO(9)})$ is 
defined over $\C[\lambda]$. 
We will denote the specialization of 
$FH_{T^2}^*(L_{\Proj^7/\mO(9)})$ to $\lambda=0$ 
by $FH_{S^1}^*(L_{\Proj^7/\mO(9)})$. 
This is  
generated by the Floer fundamental class $\Delta$ 
over the Heisenberg algebra $\mD$ 
with the following relation: 
\begin{equation}
\label{eq:PicardFucksofM89}
P^8\Delta=Q(9P+9\hbar)(9P+8\hbar)\cdots(9P+2\hbar)(9P+\hbar)\Delta,
\end{equation}
where $[P,Q]=\hbar Q$. 
By Theorem \ref{thm:freenessofefc}, 
a frame $\Phi$ of the Floer cohomology is given by 
\[
\Phi(1)=\Delta,\quad \Phi(p)=P\Delta,\quad 
\dots, \quad \Phi(p^7)=P^7\Delta, 
\] 
where we choose $p$ to be a positive generator of $H^2(\Proj^7,\Z)$. 
The connection matrix $A=(A_{kj})$ associated with $\Phi$ 
is defined by 
\[
P \Phi(p^j)=\sum_{k=0}^7 \Phi(p^k) A_{kj}.
\] 
The matrix $A$ is of the form
\[
A(Q,\hbar)=
\begin{bmatrix}
0 & 0 & & &  0 & C_0(Q,\hbar) \\
1 & 0 & & &  0 & C_1(Q,\hbar)  \\
0 & 1 & & &  0 & C_2(Q,\hbar)  \\
  &  & \ddots&&& \vdots   \\
  &   && &&                \\
0 & & & 0& 1 & C_7(Q,\hbar) 
\end{bmatrix},
\]
where $P^8\Delta=\sum_{k=0}^7 C_k(Q,\hbar)P^k\Delta$. 
The last column $P^8\Delta$ is calculated 
in the following way. By expanding 
(\ref{eq:PicardFucksofM89}) in $P$, we get 
\[
P^8 \Delta = Q (9^9 P + 9^8\cdot 45 \hbar) P^8 \Delta + O(P^7)\Delta. 
\] 
The term $P^8\Delta$ in the right-hand side is again 
replaced with the right-hand side of (\ref{eq:PicardFucksofM89}). 
Repeating the substitution for infinitely many times, 
we will arrive at $P^8\Delta=\sum_{k=0}^7 C_k(Q,\hbar)P^k\Delta$. 
Note that the remainder in each step will go to zero in the 
$Q$-adic topology. 
First few terms of the (in fact non-convergent) power series 
$C_k(Q,\hbar)$ are given by 
\tiny
\begin{align*}
C_0(Q,\hbar) & 
= 362880{\hbar}^{9}Q+843522882289920{\hbar}^{10}{Q}^{2}
+2872595183309735497205760{{\hbar}}^{11}{Q}^{3}+\cdots, \\
C_1(Q,\hbar)&=9239184\hbar^{8}Q+21617282246494176\hbar^{9}{Q}^{2}+
73846387657103705389012608{{\hbar}}^{10}{Q}^{3}+\cdots, \\
C_2(Q,\hbar)&=94988700\hbar^{7}Q+224382860804086776\hbar^{8}{Q}^{2}+
770022503217483472097175312{{\hbar}}^{9}{Q}^{3}+\cdots, \\
C_3(Q,\hbar)&=527562720\hbar^6Q+1263132210366894780\hbar^7Q^2+
4362972010749555043532127804{{\hbar}}^{8}{Q}^{3}+\cdots, \\
C_4(Q,\hbar)&=1767041325\hbar^{5}Q+4311916692248817630\hbar^{6}{Q}^{2}+
15031733439971730690200607660\hbar^{7}{Q}^{3}+\cdots, \\ 
C_5(Q,\hbar)&=3736207377\hbar^{4}Q+9369487748231192043\hbar^{5}{Q}^{2}+
33103288447539778489031223849\hbar^{6}{Q}^{3}+\cdots, \\
C_6(Q,\hbar)&=5022117450\hbar^3Q+13121510478769345653\hbar^4Q^2+
47311019540125905135150100746{{\hbar}}^{5}{Q}^{3}+\cdots, \\ 
C_7(Q,\hbar)&=4161183030\hbar^{2}Q+11618436584101043070\hbar^3Q^2+
43300442548663832211730173027{{\hbar}}^{4}Q^3+\cdots. 
\end{align*}
\normalsize
Next we perform the Birkhoff factorization of 
the fundamental solution $L$ (\ref{eq:fundsol_A}) 
of the flat connection 
$\smnabla^\hbar=\hbar Q(\partial/\partial Q)+A(Q,\hbar)$. 
By $L^{-1}=\Loc\circ \Phi$, 
$L^{-1}(1)$ is the $I$-function (\ref{eq:I-funct}) for 
$(\Proj^7,\mO(9))$ with $\lambda=0$ 
and $L^{-1}(p^k)=\Loc(P^k\Delta)$ is given by 
the derivatives of the $I$-function. 
\begin{align*} 
L^{-1} &=
\begin{bmatrix}
\vert & &\vert  \\
I_0 &  
\cdots & I_7 \\
\vert&  &\vert  
\end{bmatrix}, \quad 
I_k(Q,\hbar)= 
(\hbar Q\parfrac{}{Q}+p)^k 
\sum_{d=0}^\infty Q^d
\frac{\prod_{m=1}^{9d}(9p+m\hbar)}{\prod_{m=1}^d(p+m\hbar)^8}. 
\end{align*} 
Here, we regard $I_k$ as a column vector 
by expanding it in a basis $\{1,p,\dots,p^7\}$. 
Let $L=L_+L_-$ be the Birkhoff factorization 
(\ref{eq:Birkhoff}) of $L$. 
Let $\pi_+$ denote the projection 
$\C(\!(\hbar^{-1})\!) = \C[\hbar]\oplus \hbar^{-1}\C[\![\hbar^{-1}]\!]
\rightarrow\C[\hbar]$. 
By $\pi_+(L^{-1}L_+)=\pi_+(L_{-}^{-1})=\id$, 
we have the following recursive formula for $L_{+,d}$, where 
$L_+=\sum_{d=0}^\infty L_{+,d} Q^d$.   
\[ 
L_{+,0}=\id, \quad 
L_{+,d}=-\sum_{k=1}^d \pi_+(T_k L_{+,d-k}), \quad
\text{where } L^{-1}=\id+\sum_{d=1}^\infty T_d Q^d. 
\] 
On the other hand, for the reason of the 
grading $\deg Q=\deg \hbar^{-1} =-2$, 
$L_-$ becomes a finite series in $Q$ and $\hbar^{-1}$. 
For example, the first column of $L_-^{-1}$ 
($J$-function for a canonical frame) is given by 
\tiny
\begin{align*}
L_-^{-1}(1) = & 
1+(34138908Q/\hbar)p^2+
(56718144{Q}/{\hbar^2}+8404934443598718{Q^2}/{\hbar})p^3 \\
&
+(-22818915{Q}/{\hbar^3}+64923366053493693{Q^2}/{8\hbar^2} 
+3815933053700462506215462{Q^3}/{\hbar})p^4 \\ 
&
+(-44979543{Q}/{\hbar^4}-41161611741786333{Q^2}/{16\hbar^3}
+1568163327547517306411844{Q^3}/{\hbar^2}  \\
& \qquad 
+219544798390763529724114822821260793{Q^4}/{128\hbar})p^5 \\ 
&
+(89959086{Q}/{\hbar^5}-2387486769247188{Q^2}/{\hbar^4} 
-1841411178101141933423191{Q^3}/{2\hbar^3} \\
& \qquad 
+165593248955035194721662391017258{Q^4}/{\hbar^2} \\ 
& \qquad 
+7727272362231749241168150195184170620342513631{Q^5}/
  {12500\hbar})p^6 \\
& +(-83567214{Q}/{\hbar^6}
 +128193071703568551{Q^2}/{32\hbar^5} 
 +2536603825689258986824613{Q^3}/{12\hbar^4} \\ 
& 
\qquad -198293209598115335601311499223555059{Q^4}/{1024\hbar^3} \\
& 
\qquad -13718052706792335194606021984159356468758455727{Q^5}/
  {500000\hbar^2} \\
&  
\qquad +3542419384285237175282517996282946380767283552791120571{Q^6}/
 {20000\hbar})
p^7. 
\end{align*}
\normalsize 
The connection matrix $\A$ associated with a canonical 
frame $\Phi_{\rm can}= \Phi\circ L_+$ is independent 
of $\hbar$ and is given by 
\begin{align*} 
\A &= L_+^{-1} A L_+ + \hbar L_+^{-1} Q \parfrac{L_+}{Q}
=(L_+|_{\hbar=0})^{-1} A|_{\hbar=0} L_+|_{\hbar=0}  \\
& =
\begin{bmatrix}
0&0&0&0&0&0&0&0 \\
1&0&0&0&0&0&0&0 \\
\alpha Q&1&0&0&0&0&0&0\\
\beta{Q}^{2}&\gamma Q&1&0&0&0&0&0\\
\delta{Q}^{3}&\epsilon{Q}^{2}&\phi Q&1&0&0&0&0\\
\rho{Q}^{4}& \xi{Q}^{3}&\epsilon{Q}^{2}&\gamma Q&1&0&0&0\\
\eta{Q}^{5}& \rho{Q}^{4}&\delta {Q}^{3}& \beta{Q}^{2} &
\alpha Q & 1 & 0 & 0\\
\omega{Q}^{6}& \nu{Q}^{5}& \lambda{Q}^{4}
&\pi{Q}^{3}&\mu{Q}^{2}&\sigma Q& 1 & 0 
\end{bmatrix},
\end{align*} 
where $\alpha,\beta,\gamma,\dots$ are the following 
constants: 
\tiny
\begin{align*}
&\alpha=34138908, \quad \beta=16809868887197436,\quad
\gamma=90857052,\quad   
\delta=11447799161101387518646386, \\ 
&\epsilon=81506931029963973/2, \quad  
\phi=124756281,\quad 
\rho=219544798390763529724114822821260793/32, \\
&\xi=18892465499391490557425853, \quad 
\eta=7727272362231749241168150195184170620342513631/2500,\\
&\omega=10627258152855711525847553988848839142301850658373361713/10000,\\
&\nu=2411335276367964113374706805471621675307861731/1250,\quad 
\lambda=81865678061602904275032886226470995/32,\\
&\pi=2727763447102590732569280,\quad
\mu=2985296281746390,\quad
\sigma=5973264.
\end{align*}
\normalsize
Now we perform the reconstruction.  
Let $t^0,\dots,t^7$ be a co-ordinate system on $H^*(\Proj^7)$
dual to a basis $1,p,\dots,p^7$. 
We solve for $\hbar$-independent connection matrices 
$\A(Q,t)$, $\Omega_0(Q,t)$,
$\Omega_{1}(Q,t),\dots,\Omega_{7}(Q,t)$ 
such that they satisfy the 
conditions (\ref{eq:commutingconn}), (\ref{eq:integrableconn}) 
of flatness and  
\[
\A(Q,0) = \A(Q), \quad 
\Omega_i(Q,t) e_0 = e_i 
\]
Here, $\{e_0,\dots,e_7\}$ denotes the standard basis of 
$\C^8 \cong H^*(\Proj^7)$ corresponding to $\{1,p,\dots,p^7\}$. 
In view of the string and the divisor equation,  
it suffices to compute the deformation 
in parameters $t^2,\dots,t^7$. 
Let $\A^{(n)}$ and $\Omega_k^{(n)}$ be the degree $n$ 
part of $\A(Q,t)|_{t^0=t^1=0}$ and $\Omega(Q,t)|_{t^0=t^1=0}$ 
with respect to 
the variables $t^2,\dots,t^7$. 
Put $\A^{\le n} = \sum_{j=0}^n \A^{(j)}$ and 
$\Omega^{\le n} = \sum_{j=0}^n \Omega^{(j)}$. 
Because $t^2,t^3,\dots,t^7$ have negative degrees 
$-2,-4,\dots,-12$, it turns out that 
$\A(Q,t)|_{t^0=t^1=0}$ and $\Omega(Q,t)|_{t^0=t^1=0}$ 
are polynomial in $t^2,\dots,t^7$ and 
$\A(Q,t)|_{t^0=t^1=0}=\A^{\le 5}$, 
$\Omega_k(Q,t)|_{t^0=t^1=0}=\Omega_{k}^{\le 5}$. 
Assume inductively that we know $\A^{\le n}$.   
Because $\C^8\otimes\C[\![Q,t]\!]$ is 
generated by $e_0$ as a $\C[\![Q,t]\!][\A^{\le n}]$-module, 
this admits a unique $\C[\![Q,t]\!][\A^{\le n}]$-algebra structure 
such that $e_0$ is a unit. 
Define $\Omega_{k}^{\le n}$ to be the multiplication matrix by 
$e_k$ in this ring. 
This is calculated in the following way.    
\[
\Omega_{k}^{\le n} = 
\sum_{j=0}^7 B_{kj}(\A^{\le n})^j,
\quad \text{where $B_{kj}$ is determined by } 
e_k=\sum_{j=0}^7 B_{kj}(\A^{\le n})^j e_0.
\]
By the relation 
$(Q\partial/\partial Q) \Omega_{k}^{\le n}
=(\partial/\partial t^k) \A^{\le n+1}$, 
we compute $\A^{\le n+1}$ as 
\[
\A^{\le n+1} = \A^{\le 0}
+ \int_0^{(t^2,\dots,t^7)}
\sum_{k=2}^7 Q\parfrac{\Omega_{k}^{\le n}}{Q} dt^k.
\]
Suppose that we have obtained $\A(Q,t)|_{t^0=t^1=0}$ and 
$\Omega(Q,t)|_{t^0=t^1=0}$ in the above way. 
We will take flat co-ordinates $\hatt^0,\hatt^1,\dots,\hatt^7$ 
of the form $\hatt^k= t^k + g^k(Q)$, $g^k(0)=0$. 
By (\ref{eq:connchange_coordchange}), 
the connection matrices $\hAA$, $\Omega_{\hat{k}}$ associated 
with $\hatt^k$ are given by 
\[
\hAA(Q,t) = \A(Q,t) - \sum_{k=0}^7 Q\parfrac{g^k(Q)}{Q} 
\Omega_k(Q,t), \quad 
\Omega_{\hat{k}}(Q,t) = \Omega_k(Q,t).  
\] 
The functions $g^k$ are determined by the condition 
$e_1=\hAA(Q,0)e_0=\A(Q)e_0 - 
\sum_{k=0}^7 Q(\partial g^k/\partial Q)e_k$. 
We find (only from $\A(Q)$) 
\begin{align*}
&\hatt^0=t^0, \quad \hatt^1= t^1, \quad 
\hatt^2=t^2+\alpha Q,\quad 
\hatt^3=t^3+\frac{1}{2}\beta Q^2, \quad 
\hatt^4=t^4+\frac{1}{3}\delta Q^3, \\ 
&\hatt^5=t^5+\frac{1}{4}\rho Q^4,\quad 
\hatt^6=t^6+\frac{1}{5}\eta Q^5,\quad
\hatt^7=t^7+\frac{1}{6}\omega Q^6.
\end{align*}
By the string and the divisor equation in 
Proposition \ref{prop:divisor_AQDM}, 
we have $\hAA=\Omega_{\hat{1}}$ and  
\begin{align*}
\Omega_{\hat{k}}(Q,t(Q,\hatt)) &=  
\Omega_{\hat{k}}(Qe^{\hatt^1}, t(Q, 0,0,\hatt^2,\dots,\hatt^7)) \\
& = \Omega_{k}(Qe^{\hatt^1}, 0, 0, \hatt^2 - g^2(Q),\dots, \hatt^7-g^7(Q)) 
\end{align*} 
The matrix-valued function 
$(Q,\hatt)\mapsto \Omega_{\hat{k}}(Q,t(Q,\hatt))$
represents the multiplication by $p^k$ in 
the big quantum cohomology $QH^*_{\Euler}(\Proj^7,\mathcal{O}(9))$
with respect to the basis $\{1,p,\dots,p^7\}$.  
For simplicity, we present a $(7\times 7)$ submatrix 
$(\Omega_{\hat{1},jk})_{0\le j,k\le 6}$
of $\Omega_{\hat{1}}=(\Omega_{\hat{1},jk})_{0\le j,k\le 7}$. 
\begin{gather*}
\text{$(7\times 7)$-submatrix of } 
\Omega_{\hat{1}}=
\begin{bmatrix}
0&0&0&0&0&0&0\\
1&0&0&0&0&0&0\\
0&1&0&0&0&0&0\\
0& A & 1 & 0 & 0 & 0 & 0 \\
0& B & D & 1 & 0 & 0 & 0\\
0& C & B & A & 1 & 0 & 0\\
0& 0 & 0 & 0 & 0 & 1 &0
\end{bmatrix}, 
\end{gather*}
where the functions $A,B,C,D$ are given by 
\begin{align*} 
A&=(\gamma-\alpha)Qe^{\hatt^1}=56718144 Qe^{\hatt^1},\\
D&=(\phi-\alpha)Qe^{\hatt^1}=90617373Qe^{\hatt^1},\\
B&=(\epsilon+2\alpha(\alpha-\phi)-\beta)(Qe^{\hatt^1})^{2}+(\phi-\alpha)\hatt^2Qe^{\hatt^1} \\
&=\frac{35512880615374365}{2}(Qe^{\hatt^1})^2+90617373\, 
\hatt^2Qe^{\hatt^1}, \\
C&=(\frac{9}{2}\alpha^2(\phi-\alpha)+
\frac{3}{2}\beta(3\alpha-\gamma)-3\epsilon\alpha-\delta+\xi)
(Qe^{\hatt^1})^3 \\
&\quad +(4\alpha(\alpha-\phi)+2(\epsilon-\beta))\hatt^2 
(Qe^{\hatt^1})^2
+(\frac{\phi-\alpha}{2}(\hatt^2)^2
+(\gamma-\alpha) {\hatt^3})Qe^{\hatt^1} \\
&=4037555975532386945225553 (Qe^{\hatt^1})^3 + 35512880615374365\, \hatt^2(Qe^{\hatt^1})^2 \\ 
&\quad +(\frac{90617373}{2}(\hatt^2)^2+56718144\,\hatt^3)
Qe^{\hatt^1}.
\end{align*}
\normalsize
By Theorem \ref{thm:kkp}, this submatrix 
is related to the product by $p$ in $QH^*(M_8^9)$ as 
\[
\frac{1}{9} \int_{M_8^9} (p* p^i) \cup p^{6-j} 
= \Omega_{\hat{1},ji}, \quad 0\le i,j\le 6.  
\]
They agree with the Jinzenji's calculations 
\cite[Section 6]{jinzenji1}. 
In particular, 
the genus 0 Gromov-Witten potential 
$F_{M_8^9}$ of $M_8^9$ restricted to 
the image of $H^*(\Proj^7)\to H^*(M_8^9)$ 
is determined by 
$(\partial/\partial \hatt^1)^3F_{M_8^9}=9C$ 
and the classical part. 

\bibliographystyle{amsplain}

\end{document}